\documentclass[12pt]{article}
\usepackage[latin1]{inputenc}
\usepackage[OT2,T1]{fontenc}
\usepackage[francais]{babel}
\usepackage{lmodern}
\usepackage{amsmath}
\usepackage{amssymb}
\usepackage{mathrsfs}
\usepackage[a4paper,vscale=0.7048]{geometry}

\usepackage[ps,dvips,arrow,matrix,tips]{xy}
\entrymodifiers={+!!<0pt,\the\fontdimen22\textfont2>}
\SelectTips{cm}{10}

\makeatletter
\def\myrightarrow{{\setbox\z@\hbox{$\rightarrow$}\dimen0\ht\z@\multiply\dimen0 6\divide\dimen0 10\ht\z@\dimen0\box\z@}}

\def\myrightarrowfill@{\arrowfill@\relbar\relbar\myrightarrow}
\newcommand{\myxrightarrow}[2][]{\ext@arrow 0359\myrightarrowfill@{#1}{#2}}
\makeatother

{\setbox0\hbox{$ $}}
\mathcode`A="7041 \mathcode`B="7042 \mathcode`C="7043 \mathcode`D="7044
\mathcode`E="7045 \mathcode`F="7046 \mathcode`G="7047 \mathcode`H="7048
\mathcode`I="7049 \mathcode`J="704A \mathcode`K="704B \mathcode`L="704C
\mathcode`M="704D \mathcode`N="704E \mathcode`O="704F \mathcode`P="7050
\mathcode`Q="7051 \mathcode`R="7052 \mathcode`S="7053 \mathcode`T="7054
\mathcode`U="7055 \mathcode`V="7056 \mathcode`W="7057 \mathcode`X="7058
\mathcode`Y="7059 \mathcode`Z="705A

\lineskip=\lineskiplimit

\newcommand{\fp}[2]{#1_{\mkern-1.5mu#2}}

\newcounter{enonce}[section]
\makeatletter
\def\theenonce{\thesection.\@arabic\c@enonce}
\makeatother

\newenvironment{enonce}[1]{\noindent{\textbf{#1}}---\!\, \begin{itshape}}{\end{itshape}}
\newenvironment{proposition}[1][]{\refstepcounter{enonce}\begin{enonce}{Proposition \theenonce{}#1 }}{\end{enonce}}
\newenvironment{definition}[1][]{\refstepcounter{enonce}\begin{enonce}{Définition \theenonce{}#1 }}{\end{enonce}}
\newenvironment{theoreme}[1][]{\refstepcounter{enonce}\begin{enonce}{Théorème \theenonce{}#1 }}{\end{enonce}}
\newenvironment{corollaire}[1][]{\refstepcounter{enonce}\begin{enonce}{Corollaire \theenonce{}#1 }}{\end{enonce}}
\newenvironment{lemme}[1][]{\refstepcounter{enonce}\begin{enonce}{Lemme \theenonce{}#1 }}{\end{enonce}}
\newenvironment{conjecture}[1][]{\refstepcounter{enonce}\begin{enonce}{Conjecture \theenonce{}#1 }}{\end{enonce}}
\newenvironment{question}[1][]{\refstepcounter{enonce}\begin{enonce}{Question \theenonce{}#1 }}{\end{enonce}}
\newenvironment{questions}[1][]{\refstepcounter{enonce}\begin{enonce}{Questions \theenonce{}#1 }}{\end{enonce}}
\newenvironment{exemples*}[1][]{\refstepcounter{enonce}\noindent{\textbf{Exemples #1}}---\!\, }{}
\newenvironment{exemple*}[1][]{\refstepcounter{enonce}\noindent{\textbf{Exemple #1}}---\!\, }{}
\newenvironment{remarque}[1][]{\noindent{\textbf{Remarque{}#1 }}---\!\, }{}
\newenvironment{remarques}[1][]{\noindent{\textbf{Remarques{}#1 }}---\!\, }{}
\newenvironment{demo}[1][]{\noindent{\it{}Démonstration{}#1 }--- }{\hfill $\square$}
\newenvironment{esquissedemo}[1][]{\noindent{\it{}Esquisse de démonstration{}#1 }--- }{\hfill $\square$}

\renewcommand{\phi}{\varphi}
\renewcommand{\leq}{\leqslant}
\renewcommand{\geq}{\geqslant}
\renewcommand{\emptyset}{\varnothing}
\newcommand{\cl}{{\mathrm{cl}}}
\newcommand{\CH}{{\mathrm{CH}}}
\newcommand{\CHz}{{\mathrm{A}_0}}
\newcommand{\et}{{\text{ét}}}
\newcommand{\Qbar}{\overline \Q}
\newcommand{\Zbbarre}{\overline \Z}
\newcommand{\Hom}{{\mathrm{Hom}}}
\newcommand{\Br}{{\mathrm{Br}}}
\newcommand{\Hilb}{{\mathrm{Hilb}}}
\newcommand{\Brond}{\mathscr{B}}
\newcommand{\Crond}{\mathscr{C}}
\newcommand{\Frond}{\mathscr{F}}
\newcommand{\Irond}{\mathscr{I}}
\newcommand{\Orond}{\mathscr{O}}
\newcommand{\Srond}{\mathscr{S}}
\newcommand{\congru}[3]{#1 \equiv #2  \left(\!\!\!\!\mod #3\right)}
\newcommand{\Xrond}{\mathscr{X}}
\newcommand{\Yrond}{\mathscr{Y}}
\newcommand{\C}{\mathbf{C}}
\newcommand{\R}{\mathbf{R}}
\newcommand{\Q}{\mathbf{Q}}
\newcommand{\Z}{\mathbf{Z}}
\newcommand{\F}{\mathbf{F}}
\newcommand{\N}{\mathbf{N}}
\renewcommand{\P}{\mathbf{P}}
\newcommand{\A}{\mathbf{A}}
\newcommand{\Nrd}{\mathrm{Nrd}}

\newcommand{\uplet}[2]{#1, \mskip2.5mu \ldots \mskip-1mu, \mskip2.5mu #2}
\newcommand{\Gal}{\mathrm{Gal}}
\newcommand{\Pic}{\mathrm{Pic}}
\newcommand{\Gm}{\mathbf{G}_\mathrm{m}}
\newcommand{\Qp}{\mathbf{Q}_p}
\newcommand{\Fp}{\mathbf{F}_p}
\newcommand{\Ql}{{\mathbf{Q}_{\ell}}}
\newcommand{\Zp}{\mathbf{Z}_p}
\newcommand{\isoto}{\myxrightarrow{\,\sim\,}}
\newcommand{\kbar}{{\mkern1mu\overline{\mkern-1mu{}k\mkern-1mu}\mkern1mu}}
\newcommand{\kalg}{\kbar}
\newcommand{\Abarre}{{\mkern1mu\overline{\mkern-1mu{}A\mkern-1.5mu}\mkern1.5mu}}
\newcommand{\Bbarre}{{\mkern.4mu\overline{\mkern-.4mu{}B\mkern-1.3mu}\mkern1.3mu}}
\newcommand{\Xbarre}{{\mkern.4mu\overline{\mkern-.4mu{}X\mkern-1mu}\mkern1mu}}
\newcommand{\Ybarre}{{\mkern.15mu\overline{\mkern-.15mu{}Y\mkern-1mu}\mkern1mu}}
\newcommand{\Zbarre}{{\mkern1mu\overline{\mkern-1mu{}Z\mkern-1mu}\mkern1mu}}
\newcommand{\Ubarre}{{\mkern.77mu\overline{\mkern-.77mu{}U\mkern-1mu}\mkern1mu}}
\newcommand{\xtilde}{{\widetilde{x}}}
\newcommand{\ytilde}{{\widetilde{y}}}
\DeclareMathOperator{\Spec}{Spec}
\DeclareMathOperator{\inv}{inv}
\DeclareMathOperator{\Card}{Card}
\DeclareMathOperator{\Ker}{Ker}
\newcommand{\xsouligne}{\mkern.2mu\underline{\mkern-.2mu{}x\mkern-1mu}\mkern1mu}
\newcommand{\alphasouligne}{{\mkern.2mu\underline{\mkern-.2mu{}\alpha\mkern-1mu}\mkern1mu}}
\numberwithin{equation}{subsection}

\title{La connexité rationnelle en arithmétique}
\date{}\author{Olivier Wittenberg}
\begin{document}
\maketitle
\insert\footins{\footnotesize{}Version du 28 septembre 2008, révisée le 30 juin 2010, mise à jour le 24 septembre 2010.}
Une variété algébrique propre et lisse~$X$ sur un corps~$k$ de caractéristique~$0$
est dite \emph{rationnellement connexe} si pour toute extension algébriquement close~$K$ de~$k$,
par tout couple de $K$\nobreakdash-points de~$X$ passe une courbe rationnelle (définie sur~$K$).
Cette notion, introduite au début des années 1990 par Kollár, Miyaoka et Mori, et indépendamment par Campana, a d'abord joué un rôle important dans l'étude
de la géométrie des variétés complexes.  Le développement des techniques géométriques propres aux variétés rationnellement connexes s'est ensuite répercuté
en arithmétique. Ainsi l'article fondateur de Kollár~\cite{kollarloc} établissait-il, pour toutes les variétés rationnellement connexes définies sur
un corps $p$\nobreakdash-adique, la finitude de la $R$\nobreakdash-équivalence --- une propriété de nature arithmétique
qui jusque-là n'était connue que dans des cas très particuliers et qui n'avait même pu être envisagée
dans cette généralité, faute  de disposer de la notion de connexité rationnelle.  Depuis une dizaine d'années, plusieurs autres résultats concernant
l'arithmétique des variétés rationnellement connexes ont vu le jour.  C'est sur ces résultats que nous nous proposons de faire le point dans le présent rapport.

Le premier chapitre introduit brièvement les principales questions qui se posent dans l'étude de l'arithmétique des variétés rationnellement connexes.
Dans chacun des chapitres suivants,
un théorème général est discuté et démontré.  Le corps de base est un corps $p$\nobreakdash-adique (ou fertile)
au second chapitre, un corps pseudo-algébriquement clos (puis fini ou $p$\nobreakdash-adique) au troisième chapitre,
un corps fini au quatrième chapitre.
On ne connaît à ce jour aucun résultat qui s'applique à toutes les variétés rationnellement connexes définies sur~$\Q$, exception faite du corollaire~\ref{corkolszcdn}
ci-dessous. Les corps de nombres ne joueront donc qu'un rôle mineur dans ce rapport.
Il~existe toutefois une abondante littérature consacrée à l'arithmétique de diverses classes
de variétés rationnellement connexes sur les corps de nombres: surfaces rationnelles, espaces homogènes de groupes
linéaires, intersections d'hypersurfaces de bas degré dans l'espace projectif, fibrations en des variétés de l'un de ces trois types.
Nous ne tenterons pas de la survoler ici.

Les chapitres~\ref{seckoll} et~\ref{seckollsz} font appel à des techniques de déformation de courbes rationnelles.  Leur lecture présuppose un minimum de familiarité
avec la plus simple d'entre elles dans le cas où le corps de base est algébriquement clos, telle qu'elle est présentée dans~\cite[§4]{bonavero}.

\bigskip
\noindent{}\textbf{Remerciements.}\;  Ces notes ont constitué le support d'une série de cinq exposés
à la session SMF États de la Recherche «~Variétés rationnellement connexes: aspects géométriques et arithmétiques~» tenue à Strasbourg en mai~2008,
dont je remercie les organisateurs.
Je tiens d'autre part à remercier Antoine Chambert-Loir, Jean-Louis
Colliot-Thélène, Olivier Debarre, Hélène Esnault, Bruno Kahn et János Kollár
pour leurs réponses à mes questions, pour leurs remarques et pour de nombreuses discussions éclairantes,
et enfin le rapporteur pour sa lecture vigilante.

\bigskip
\noindent{}\textbf{Conventions.}\;  Une \emph{variété} est un schéma de type fini sur un corps.
Soit~$X$ une variété sur un corps~$k$.  Un \emph{point rationnel} de~$X$ est un $k$\nobreakdash-point de~$X$.  L'ensemble des points rationnels est noté $X(k)$.
Soit~$K$ un corps algébriquement clos non dénombrable contenant~$k$.  On dit que la variété~$X$ est \emph{rationnelle} (resp.~\emph{unirationnelle})
si $X \otimes_k K$ l'est, c'est-à-dire si $X \otimes_k K$ est birationnellement équivalente à un espace projectif (resp.~s'il existe une application rationnelle dominante
d'un espace projectif vers $X \otimes_k K$).  Si~$X$ est propre sur~$k$, on dit que~$X$ est \emph{rationnellement connexe}
(resp.~\emph{rationnellement connexe par chaînes}, \emph{séparablement rationnellement connexe}) si $X \otimes_k K$ l'est, c'est-à-dire si $X \otimes_k K$ vérifie les définitions données
dans~\cite{bonavero} ou dans~\cite{kollrational}.  Nous convenons que les qualificatifs \emph{rationnelle}, \emph{unirationnelle}, \emph{rationnellement connexe} et \emph{séparablement rationnellement connexe}
sous-entendent que $X\otimes_k K$ est irréductible (en particulier non vide).
On dira que $X$ est \emph{$k$\nobreakdash-rationnelle} (resp.~\emph{$k$\nobreakdash-unirationnelle})
s'il existe une application rationnelle $\P^n_k \dashrightarrow X$ qui soit birationnelle (resp.~dominante).
Une \emph{conique} (ou \emph{conique projective}) est une courbe projective plane de degré~$2$; rappelons que toute courbe rationnelle propre et lisse
est isomorphe à une conique.

\tableofcontents

\section{Un aperçu de quelques problèmes concernant l'arithmétique des variétés rationnellement connexes}
\label{secundeux}

\subsection{Corps~$(C_i)$}
\label{subseccorpsci}
\medskip

Le paragraphe~\ref{subseccorpsci} est consacré aux corps~$(C_i)$; les variétés rationnellement connexes n'y apparaissent qu'implicitement,
sous la forme d'hypersurfaces projectives de bas degré ou d'intersections d'icelles.
La propriété~$(C_1)$, introduite par Artin~\cite{artincollectedpapers}, est notamment liée aux questions d'existence de points rationnels sur
les variétés de Severi--Brauer ainsi que sur certaines variétés de Fano. (Les unes comme les autres sont des exemples de variétés rationnellement connexes (par chaînes), d'après
les travaux de Campana, Kollár, Miyaoka et Mori.)

\subsubsection{Définition, exemples et théorèmes de transition}

\begin{definition}[ (Artin, Lang)]%
Soient~$k$ un corps et $i \geq 0$ un entier.  On dit que $k$ est \emph{un corps~$(C_i)$} si pour tout~$n$ et tout~$d$, toute hypersurface de~$\P^n_k$ de degré~$d$
avec $n \geq d^i$ admet un point $k$\nobreakdash-rationnel.
\end{definition}

\bigskip
Autrement dit, le corps~$k$ est~$(C_i)$ si les équations de la forme
$$
f(\uplet{x_0}{x_n})=0 \rlap{\text{,}}
$$
où $f \in k[\uplet{x_0}{x_n}]$ est un polynôme homogène de degré $d>0$, admettent une solution
$(\uplet{x_0}{x_n}) \in k^{n+1} \setminus \{(\uplet{0}{0})\}$ dès que $n \geq d^i$.

Les corps~$(C_0)$ sont les corps algébriquement clos.
Chevalley~\cite{chevalleyfinisqac} a démontré
que les corps finis sont~$(C_1)$:

\bigskip
\begin{theoreme}[ (Chevalley--Warning)]%
Soit~$k$ un corps fini de caractéristique~$p$.
Soit $H \subset \P^n_k$ une hypersurface de degré~$d$ avec $n \geq d$.
Alors $\congru{\Card H(k)}{1}{p}$, et en particulier $H(k) \neq \emptyset$.
\end{theoreme}

\bigskip
\begin{demo}[ (due à Ax)]%
Notons~$q$ le cardinal de~$k$ et $f \in k[\uplet{x_0}{x_n}]$ un polynôme homogène de degré~$d$ s'annulant sur~$H$.  Soit~$N$ le nombre de solutions
dans $k^{n+1}$ de l'équation $f=0$.  Comme
$\Card H(k) = (N-1)/(q-1)$, il suffit d'établir la congruence $\congru{N}{0}{p}$.

Posons $F=1-f^{q-1}$.  Le polynôme~$F$ ne prend sur $k^{n+1}$ que les valeurs~$0$ et~$1$; il vaut~$1$ précisément sur les solutions de l'équation $f=0$.
D'où l'égalité
\begin{equation}
\label{eqchevwar}
N=\sum_{\xsouligne \in k^{n+1}} F(\xsouligne) \text{.}
\end{equation}
Pour tout entier $\alpha$ tel que $0 \leq \alpha < q-1$, on a $\sum_{x \in k} x^\alpha=0$ (si $\alpha=0$ c'est clair, sinon c'est un petit exercice).
Par conséquent, pour $\uplet{\alpha_0}{\alpha_n} \in \N$, on a
$$
\sum_{(\uplet{x_0}{x_n}) \in k^{n+1}} x_0^{\alpha_0} \dots x_n^{\alpha_n}=0
$$
dès que $\min(\alpha_i)<q-1$, en particulier dès que $\alpha_0 + \dots + \alpha_n < (n+1)(q-1)$.  Il~s'ensuit que
$\sum_{\xsouligne \in k^{n+1}} G(\xsouligne)=0$
pour tout $G \in k[\uplet{x_0}{x_n}]$ de degré strictement inférieur à $(n+1)(q-1)$.  Appliquons cela au polynôme~$F$, qui est de degré $d(q-1)$, et combinons
l'égalité obtenue avec~(\ref{eqchevwar}): on trouve que~$N$ s'annule dans~$k$, autrement dit $\congru{N}{0}{p}$.
\end{demo}

\bigskip
\begin{remarques}
(i) L'hypothèse que le polynôme~$f$ est homogène n'a pas été utilisée, ainsi le résultat prouvé est quelque peu plus général que celui énoncé.

(ii) Toujours sous les hypothèses du théorème de Chevalley--Warning, Ax~\cite{axzeroes} a démontré que l'on a même $\congru{\Card H(k)}{1}{q}$,
où~$q$ désigne le cardinal de~$k$.
\end{remarques}

\bigskip
À partir des corps algébriquement clos et des corps finis, il est facile de fabriquer des exemples de corps~$(C_i)$ pour $i>1$ grâce à la propriété de transitivité suivante:

\bigskip
\begin{theoreme}[ (Tsen--Lang--Nagata)]%
\label{tsenlangnagata}
Soit~$k'/k$ une extension de corps de degré de transcendance $d < \infty$, et soit~$i$ un entier naturel.
Si~$k$ est un corps~$(C_i)$, alors $k'$ est un corps $(C_{i+d})$.
\end{theoreme}

\bigskip
En particulier, le corps $\C(t)$, et plus généralement le corps des fonctions d'une courbe sur un corps algébriquement clos, est~$(C_1)$ (c'est le théorème de Tsen),
et le corps des fonctions de toute variété algébrique intègre de dimension~$i$ sur un corps algébriquement clos (resp.~fini) est un exemple de corps~$(C_i)$
(resp.~$(C_{i+1})$, d'après le théorème de Chevalley).

\bigskip
\begin{demo}
Nous nous contentons ici de démontrer le cas particulier le plus important du théorème~\ref{tsenlangnagata},
c'est-à-dire celui où $k'=k(t)$ et où~$k$ est algébriquement clos.
Pour le cas général, le principe est le même, mais il vaut mieux commencer
par établir le théorème d'Artin--Lang--Nagata dont il est question au~§\ref{subsecintersectionscompletes} ci-dessous.
Nous renvoyons à \cite[Chapter~3]{greenberglecturesonforms} pour une démonstration complète.

Soient donc~$k$ un corps algébriquement clos et $f \in k(t)[\uplet{x_0}{x_n}]$ un polynôme homogène de degré $d \leq n$ à coefficients dans~$k(t)$.
Nous voulons prouver que l'équation $f=0$ admet une solution dans $k(t)^{n+1} \setminus \{(\uplet{0}{0})\}$.  Écrivons~$f$ comme
\begin{equation}
\label{eqtsen}
f = \sum_{\alphasouligne=(\uplet{\alpha_0}{\alpha_n}) \in \N^{n+1}} c_\alphasouligne x_0^{\alpha_0}\dots x_n^{\alpha_n}
\end{equation}
avec $c_\alphasouligne \in k(t)$.  Quitte à multiplier~$f$ par un scalaire,
on peut supposer que les $c_\alphasouligne$ sont dans~$k[t]$.  Soit~$N$ un entier assez grand, à préciser.
Posons $$x_i=y_{i0} + t y_{i1}+\dots + t^N y_{iN}$$ pour chaque $i \in \{\uplet{0}{n}\}$,
où les $y_{ij}$ sont des indéterminées.  Réécrivons~(\ref{eqtsen}) en termes des $y_{ij}$, développons, et rassemblons les monômes obtenus selon les puissances
de~$t$; on aboutit à $f=\sum_{m \geq 0} t^m \phi_m$ où les~$\phi_m$ sont des polynômes en les $y_{ij}$ et à coefficients dans~$k$.  Notant~$\delta$ le maximum des degrés
des polynômes $c_\alphasouligne$, on a $\phi_m=0$ pour $m>Nd+\delta$.  Le système
$\phi_0 = \phi_1 = \dots = \phi_{Nd+\delta} = 0$
est un système de $Nd+\delta+1$ équations polynomiales homogènes en $(N+1)(n+1)$ variables, les $y_{ij}$.  Comme $d < n+1$, ce système admettra une solution
dans $k^{(N+1)(n+1)} \setminus \{(\uplet{0}{0})\}$ si l'on choisit~$N$ assez grand, puisque~$k$ est algébriquement clos. Une telle solution est
précisément ce que l'on cherche.
\end{demo}

\bigskip
Soit maintenant~$A$ un anneau de valuation discrète complet, de corps résiduel~$\kappa$, de corps des fractions~$K$.
Si~$\kappa$ est~$(C_i)$, peut-on conclure que~$K$ est~$(C_{i+1})$~?  En général non,
comme il résulte du contre-exemple de Terjanian à la conjecture
d'Artin (voir~\textsection\ref{parterjanian}).  Néanmoins la réponse est oui dans deux cas notables:
le cas d'égale caractéristique (Greenberg), et le cas où~$i=0$ (Lang).
Pour ces deux résultats il n'est même pas nécessaire de supposer~$A$ complet; il suffit qu'il soit hensélien et que le complété de~$K$ soit une extension
séparable de~$K$ (ce qui bien sûr est automatique si~$K$ est de caractéristique~$0$).
Nous renvoyons à~\cite{terjanianprogresrecents} pour un résumé des résultats de Greenberg et à~\cite{greenbergihes} pour les démonstrations.
Les deux corollaires les plus importants sont:

\bigskip
\begin{theoreme}[ (Greenberg)]%
\label{thgreenberg}
Si~$k$ est un corps~$(C_i)$, le corps~$k((t))$ des séries formelles en une variable à coefficients dans~$k$ est $(C_{i+1})$.
\end{theoreme}

\bigskip
\begin{theoreme}[ (Lang)]%
Soit~$p$ un nombre premier. L'extension non ramifiée maximale de~$\Qp$ (ou plus généralement d'un corps $p$\nobreakdash-adique) est un corps~$(C_1)$.
\end{theoreme}

\subsubsection{Le cas de~$\Qp$}
\label{parterjanian}

Le corps~$\Qp$ lui-même n'est pas~$(C_1)$: il est facile d'exhiber des coniques dépourvues de point rationnel
sur~$\Qp$. Si $p\neq 2$ et si $a \in \Zp^\star$ n'est pas un carré,
la courbe d'équation homogène $x^2 - a y^2 + p z^2 = 0$ en est un exemple,
tandis que
si $p=2$, la courbe d'équation homogène $x^2+y^2+z^2=0$ en est un.

Artin avait conjecturé que~$\Qp$ est~$(C_2)$.  Il est vrai que toute hypersurface de degré~$d$ dans~$\P^n_{\Qp}$
avec $n \geq d^2$ admet un point $\Qp$\nobreakdash-rationnel si $d=2$ (Hasse~\cite{hassequadriques}) ou si $d=3$ (Demjanov~\cite{demjanov} pour $p\neq 3$ et Lewis~\cite{lewis} en général).
Néanmoins, Terjanian~\cite{terjaniancontreexemple} a exhibé un contre-exemple
à la conjecture d'Artin, avec $d=4$ et $p=2$.  Ce contre-exemple fut ensuite généralisé
par divers auteurs pour aboutir finalement au

\bigskip
\begin{theoreme}[ (Arkhipov--Karatsuba~\cite{arkhipovkaratsuba}, Alemu~\cite{alemu})]%
Soit~$p$ un nombre premier.  Soit~$k$ un corps $p$\nobreakdash-adique, c'est-à-dire une extension finie de~$\Qp$.
Alors~$k$ n'est~$(C_i)$ pour aucun $i \geq 0$.
\label{alemu}
\end{theoreme}

\bigskip
La question suivante n'est pas encore résolue (cf.~\cite[p.~6]{terjanianprogresrecents} et~\cite[p.~67]{wooleysurv}): une hypersurface de~$\P^n_k$ de
degré~$d$, avec $n \geq d^2$, admet-elle nécessairement un point rationnel si~$k$
est $p$\nobreakdash-adique \emph{et si~$d$ est impair}~?

Dans le sens d'une réponse positive à la conjecture d'Artin, signalons le
célèbre théorème d'Ax--Kochen~\cite{axkochen} selon lequel pour~$d$ fixé, les
corps~$\Qp$ sont «~$(C_2)$ en degré~$d$~» sauf pour un nombre fini
de~$p$. (L'ensemble des mauvais~$p$ dépend de~$d$ de manière non explicite
mais calculable --- du moins en théorie, cf.~\cite[Theorem~9]{axkochen3}).
On trouvera dans~\cite{terjanianprogresrecents} une variante portant sur tous les corps $p$\nobreakdash-adiques: pour tout~$d$, il existe~$p_0$ tel que pour tout $p \geq p_0$,
tout corps $p$\nobreakdash-adique est «~$(C_2)$ en degré~$d$~».  Dans le cas où $d=5$, Leep et Yeomans~\cite{leepyeomans} établissent le résultat
plus précis suivant: pour tout corps $p$\nobreakdash-adique~$k$ dont le corps résiduel possède au moins~$47$ éléments,
toute hypersurface de~$\P^n_k$ de degré~$d=5$ avec $n \geq d^2=25$ admet un point $k$\nobreakdash-rationnel.  De même,
d'après Wooley~\cite{wooleyartin}, la conjecture vaut pour $d=7$ (resp.~pour~$d=11$) sur un corps $p$\nobreakdash-adique~$k$ dès que le corps
résiduel de~$k$ possède au moins $887$ (resp.~$8059$) éléments.

La démonstration d'Ax et Kochen du théorème d'Ax--Kochen fait appel à la théorie des modèles.  Denef a récemment annoncé une preuve géométrique de ce théorème, ainsi que d'une
généralisation conjecturée par Colliot-Thélène~\cite[§3]{ctcras}.
Celle-ci étend le théorème d'Ax--Kochen à des variétés qui ne sont pas nécessairement des hypersurfaces de~$\P^n_k$.

Si~$k$ est un corps $p$\nobreakdash-adique et si $H \subset \P^n_k$ est une hypersurface de degré~$d$ avec~\mbox{$n \geq d^2$}, Kato et Kuzumaki~\cite{katokuzumaki} conjecturent que~$H$
contient un $0$\nobreakdash-cycle de degré~$1$ (autrement dit: les degrés des points fermés de~$H$ sont premiers entre eux dans leur ensemble),
et établissent cette conjecture dans le cas où~$d$ est premier.  Une version géométrique de cette conjecture est suggérée dans~\cite[Remarque~2]{ctcras}.

\subsubsection{Autres corps; quelques questions ouvertes}
\label{subsubsecautrescorps}

Il est très facile de voir que le corps~$\R$ des nombres réels n'est lui non plus~$(C_i)$ pour aucun $i \geq 0$.  Cela implique (exercice~!) que si~$X$ est une variété intègre lisse sur~$\R$,
et si $X(\R)\neq\emptyset$, alors le corps des fonctions de~$X$ n'est~$(C_i)$ pour aucun $i \geq 0$.  Lorsque $X(\R)=\emptyset$, la situation est moins claire:

\bigskip
\begin{conjecture}[ (Lang~\cite{langrealplaces})]%
\label{conjlangreel}
Soit~$X$ une variété intègre sur~$\R$ (ou, plus généralement, sur un corps réel clos).  Soit $i=\dim(X)$.
Si $X(\R)=\emptyset$, alors le corps des fonctions de~$X$ est~$(C_i)$.
\end{conjecture}

\bigskip
Cette conjecture audacieuse est ouverte même dans le cas où~$X$ est la conique réelle sans point réel, c'est-à-dire la courbe affine d'équation $x^2+y^2=-1$.
Lorsque~$X$ est de dimension~$2$, on ne sait même pas si toute hypersurface quadrique dans~$\P^4_{\R(X)}$ admet un point rationnel.

Si maintenant~$k$ est un corps de nombres, il résulte du théorème~\ref{alemu} que~$k$ n'est~$(C_i)$ pour aucun $i\geq 0$. Cependant, le problème suivant est ouvert:

\bigskip
\newcommand{\refartin}{\cite[p.~x]{artincollectedpapers}}
\begin{conjecture}[ (Artin~\refartin)]%
Le corps $\Q^{\mathrm{ab}}$ est~$(C_1)$.
\end{conjecture}

\bigskip
Rappelons que $\Q^{\mathrm{ab}}$, qui est par définition l'extension abélienne maximale de~$\Q$, est aussi le sous-corps de~$\Qbar$ engendré par les racines de l'unité (théorème de Kronecker--Weber).

Si~$k$ est un corps algébriquement clos, le corps des séries formelles itérées $k((x))((y))$ est~$(C_2)$ d'après le théorème~\ref{thgreenberg}.
La question de savoir si le corps des fractions $k((x,y))$ de l'anneau $k[[x,y]]$ est lui aussi~$(C_2)$ est en revanche ouverte,
même si $k=\C$.  Dans cette direction on dispose seulement d'un théorème dû à Choi, Dai, Lam et Reznick~\cite{cdlr} selon lequel $k((x,y))$ est «~$(C_2)$ pour les hypersurfaces
diagonales~».  Plus généralement, on peut se demander si $k((\uplet{x_1}{x_n}))$ est~$(C_n)$.

Citons pour terminer une question soulevée par Ax sur la propriété~$(C_1)$ pour les corps pseudo-algébriquement clos.

\bigskip
\begin{definition}[ (Ax, Frey, Jarden)]%
\label{defpac}
Un corps~$k$ est dit \emph{pseudo-algébriquement clos} si toute variété géométriquement intègre sur~$k$ admet un point $k$\nobreakdash-rationnel.
\end{definition}

\bigskip
Les corps algébriquement clos et les extensions algébriques infinies de corps finis sont des exemples de corps
pseudo-algébriquement clos (cf.~\cite[Corollary~11.2.4]{friedjarden}).

\bigskip
\newcommand{\refax}{\cite[p.~270, Problem~3]{axconj}}
\begin{question}[ (Ax~\refax)]%
\label{questionax}
Tout corps parfait et pseudo-algébriquement clos est-il~$(C_1)$~?
\end{question}

\bigskip
Cette question est équivalente à la question suivante (cf.~\cite[Corollary~21.3.3]{friedjarden} ou
une variante du lemme~\ref{lemmepac} ci-dessous):
pour tout corps~$k$ et tout entier $n \geq 1$, toute hypersurface $H \subset \P^n_k$
de degré au plus~$n$ contient-elle une sous-variété géométriquement irréductible~?

Si~$k$ est un corps parfait pseudo-algébriquement clos, on sait que~$k$ est~$(C_1)$
dès que l'une au moins des trois hypothèses suivantes est vérifiée:
\begin{itemize}
\item le groupe de Galois absolu de~$k$ est
abélien (Ax~\cite[Theorem~D]{axconj});
\item le corps~$k$ contient un corps algébriquement clos (Denef--Jarden--Lewis~\cite{denefjardenlewis});
\item le corps~$k$ est
de caractéristique nulle (Kollár~\cite{kollarax}).
\end{itemize}
Kollár démontre en réalité
dans~\cite{kollarax} l'énoncé suivant, d'une portée bien plus large: sur un
corps pseudo-algébriquement clos de caractéristique nulle, toute variété
obtenue en faisant dégénérer une variété de Fano (qui n'est pas nécessairement
une hypersurface ni même une intersection d'hypersurfaces de bas degré dans un
espace projectif) admet un point rationnel.  Hogadi et Xu~\cite{hogadixu} ont depuis généralisé ce résultat
aux variétés rationnellement connexes (qui ne sont pas nécessairement de Fano).

\subsubsection{Intersections d'hypersurfaces de bas degré dans l'espace projectif}
\label{subsecintersectionscompletes}

Soient~$i,n,r \geq 1$ des entiers, soit~$k$ un corps~$(C_i)$ et soient $\uplet{H_1}{H_r} \subset \P^n_k$ des hypersurfaces.
Notons $X = H_1 \cap \dots \cap H_r$ et supposons que \mbox{$n \geq d_1^i + \cdots + d_r^i$}.
Par définition de la propriété~$(C_i)$, la variété~$X$ contient un point rationnel si~$r=1$.
Pour~$r$ quelconque, a-t-on encore nécessairement $X(k)\neq\emptyset$~?
En toute généralité, cette question est ouverte.  On sait néanmoins qu'elle
admet une réponse affirmative si $d_1=\dots=d_r$ (Artin--Lang--Nagata~\cite{nagata}) ou si~$k$
possède une forme normique de niveau~$i$ et de degré~$d$ pour tout entier
$d>1$ (Lang~\cite[Theorem~4]{langonquasi}).
Une \emph{forme normique de niveau~$i$ et de degré~$d$} est un polynôme homogène~$f$ à coefficients dans~$k$,
de degré~$d$, en $d^i$ variables, tel que l'hypersurface projective d'équation \mbox{$f=0$} n'admette pas de point rationnel.
La terminologie vient de ce que si $\ell/k$ est une extension finie, la norme de~$\ell$ à~$k$, exprimée dans une base quelconque du $k$\nobreakdash-espace vectoriel~$\ell$, est une forme normique
de niveau~$1$ et de degré \mbox{$[\ell:k]$}.  En particulier, un corps fini, et plus généralement tout corps possédant une extension finie de chaque degré, admet une forme normique
de niveau~$1$ et de degré~$d$ pour tout~$d>1$.
Pour la plupart des corps~$(C_i)$ que l'on rencontre en géométrie arithmétique,
l'existence de formes normiques de tout degré est assurée par la propriété suivante (cf.~\cite[p.~377]{langonquasi}): s'il existe sur~$k$ une valuation
discrète (de rang~$1$) dont le corps résiduel
possède une forme normique de niveau~$i-1$ et de degré~$d$, alors~$k$ possède
une forme normique de niveau~$i$ et de degré~$d$. (On prend pour
convention qu'un corps algébriquement clos possède une forme normique de
niveau~$0$ et de degré~$d$ pour tout $d>1$.)  Ainsi, notamment, le corps des fonctions d'une variété intègre de dimension~$i$ (resp.~$i-1$) sur un corps algébriquement clos (resp.~sur
un corps fini ou sur l'extension non ramifiée maximale d'un corps $p$\nobreakdash-adique)
possède une forme normique de niveau~$i$ et de degré~$d$ pour tout $d>1$. Si~$k$ est un tel corps, on a donc bien $X(k)\neq\emptyset$.

\subsection{Interlude: surfaces rationnelles et corps~$(C_1)$}
\label{subsecinterlude}

\medskip
Jusqu'ici nous avons constaté que la propriété~$(C_1)$ pour le corps~$k$ exerce une influence sur les points rationnels des hypersurfaces de degré~$d$ dans~$\P^n_k$ si \mbox{$n \geq d$}
(par définition) et plus généralement sur les points rationnels
des intersections d'hypersurfaces de degrés $\uplet{d_1}{d_r}$ dans~$\P^n_k$ si \mbox{$n \geq \sum_{i=1}^r d_i$}
(d'après le~§\ref{subsecintersectionscompletes}).  Lorsqu'elles sont lisses, ces variétés sont rationnellement connexes,
au moins si~$k$ est de caractéristique nulle.  Inversement, si~$X$ est une variété propre, lisse, rationnellement connexe, sur un corps~$k$ de caractéristique nulle, mais que~$X$
n'est pas donnée comme une intersection d'hypersurfaces de bas degré
dans un espace projectif, il n'y a aucune raison \emph{a priori} de s'attendre à ce que la propriété~$(C_1)$ pour~$k$ ait un quelconque
impact sur
l'ensemble $X(k)$.  C'est pourtant ce qui se produit pour plusieurs classes de variétés rationnellement connexes, comme l'illustrent
l'exemple des variétés de Severi--Brauer (cf.~§\ref{subsecgroupebrauer} ci-après), le théorème de Springer--Steinberg (selon lequel si~$k$ est un corps parfait~$(C_1)$,
tout espace homogène sous un groupe algébrique linéaire connexe sur~$k$ admet un point rationnel; cf.~\cite[III.2.4]{serrecg}), et le

\bigskip
\newcommand{\refmanin}{\cite[§4]{maninihes}}
\begin{theoreme}[ (Manin~\refmanin, Colliot-Thélène~\cite{cticm})]%
\label{thmmct}
Soit~$X$ une surface propre, lisse, rationnelle, sur un corps~$k$.
Si~$k$ est~$(C_1)$ alors $X(k)\neq\emptyset$.
\end{theoreme}

\bigskip
La démonstration s'appuie sur la classification des surfaces rationnelles sur
un corps arbitraire, établie par Enriques, Manin et Iskovskikh~\cite{iskoizv}
et qui s'énonce ainsi.  Soit~$k$ un corps.
Toute surface sur~$k$ propre, lisse et rationnelle s'obtient par un nombre fini d'éclatements de centre lisse
à partir d'une surface de l'un des deux types suivants:
\begin{enumerate}
\item Les surfaces propres, lisses et géométriquement connexes dont le faisceau anti-canonique est ample.
\item Les surfaces~$X$ propres, lisses et géométriquement connexes telles qu'il existe une courbe~$C$ lisse, connexe, rationnelle et un morphisme $\pi \colon X \rightarrow C$
dont la fibre générique est lisse et dont les fibres sont toutes isomorphes à des coniques irréductibles géométriquement réduites.
\end{enumerate}
Les surfaces de type~1 (resp.~2) sont dites \emph{de del Pezzo} (resp.~\emph{fibrées en coniques au-dessus d'une conique}; la courbe~$C$ étant propre, lisse, rationnelle, est en effet
isomorphe à une conique projective lisse \emph{via} le plongement anti-canonique).  Elles sont toujours rationnelles.

\bigskip
\begin{esquissedemo}[ du théorème~\ref{thmmct}]%
(Voir \cite[IV.6.8]{kollrational} pour les détails.)

Soient~$k$ un corps~$(C_1)$ et~$X$ une surface rationnelle, propre et lisse sur~$k$.  Pour établir le théorème, on peut supposer que~$X$ est de l'un des deux types standard
d'Enriques--Manin--Iskovskikh.

Si $\pi \colon X \rightarrow C$ est une surface rationnelle fibrée en coniques, on a $C(k) \neq \emptyset$ puisque~$C$ est isomorphe à une conique
et que~$k$ est~$(C_1)$. Fixons $c \in C(k)$. La fibre de~$\pi$ en~$c$ admet à son tour un point rationnel, étant elle aussi isomorphe à une conique.
D'où $X(k)\neq\emptyset$.

Supposons maintenant~$X$ de del Pezzo.  On appelle \emph{degré} de~$X$ le nombre d'auto-intersection de son faisceau canonique.
Il est compris entre~$1$ et~$9$.  Notons-le~$d$.  Pour $d \geq 3$, le faisceau anti-canonique est très ample et permet de voir~$X$ comme une sous-variété
de degré~$d$ dans~$\P^d_k$.
Si $d=3$, la surface~$X$ s'identifie ainsi à une surface cubique dans $\P^3_k$,
de sorte que $X(k)\neq\emptyset$ par définition de la propriété~$(C_1)$.  Si $d=4$, la surface~$X$ est une intersection lisse de deux quadriques
dans~$\P^4_k$.  D'après Artin--Lang--Nagata, cela entraîne de nouveau que $X(k)\neq\emptyset$ (cf.~§\ref{subsecintersectionscompletes}).
Swinnerton-Dyer a montré que les surfaces de del Pezzo de degré~$5$ ont toujours un point rationnel, sans hypothèse sur le corps de base.
Il en va de même pour les surfaces de del Pezzo de degrés~$1$ et~$7$, pour des raisons triviales (si $d=1$, le système linéaire anti-canonique admet un unique point base; si $d=7$,
la surface est isomorphe à l'éclaté de~$\P^2_k$ en deux points rationnels ou en un point fermé de degré~$2$).
Les surfaces de del Pezzo de degré~$9$ sont les variétés de Severi--Brauer de dimension~$2$, \emph{i.e.} les surfaces
qui après une extension finie des scalaires deviennent isomorphes au plan projectif.  Nous montrerons au~§\ref{subsecgroupebrauer} qu'elles admettent un point rationnel
sur tout corps~$(C_1)$.  Un argument de même nature permet de traiter plus généralement les surfaces de del Pezzo de degré~$\geq 6$ (on se ramène à la propriété:
sur un corps~$(C_1)$, tout espace principal homogène sous un tore
admet un point rationnel~\cite[X.7]{serrecorpslocaux}).
Restent les surfaces de del Pezzo de degré~$2$.  Supposons pour simplifier que~$k$ soit de caractéristique $\neq 2$.
Les surfaces de del Pezzo de degré~$2$ sur~$k$ sont alors les revêtements doubles de~$\P^2_k$ ramifiés le long d'une courbe quartique lisse.  Elles sont donc définies
par une équation «~homogène~» de la
forme $t^2=f(x,y,z)$, où~$f$ est un polynôme homogène de degré~$4$ et où les
variables $x$, $y$, $z$ sont de poids~$1$ et~$t$ est de poids~$2$.
Si tout élément de~$k$ est un carré, on a $X(k)\neq\emptyset$ de façon évidente.
Sinon, il existe un polynôme homogène $g \in k[u,v]$ de degré~$2$, tel que l'équation $g(u,v)=0$ n'ait pas de solution
dans $k^2 \setminus \{(0,0)\}$.  L'équation $g(u,v)^2=f(x,y,z)$ définit alors une hypersurface de degré~$4$ dans $\P^4_k$.
Comme~$k$ est $(C_1)$, cette hypersurface admet un $k$\nobreakdash-point, d'où $X(k)\neq\emptyset$.
\end{esquissedemo}

\bigskip
Ces diverses considérations nous amènent à poser la question suivante, qui fait partie du folklore:

\bigskip
\begin{question}
\label{laquestionc1}
Soit~$X$ une variété propre, lisse, séparablement\footnote{Une variété~$X$ géométriquement connexe sur~$k$ est dite \emph{séparablement rationnellement connexe}
s'il existe
un morphisme $f \colon \P^1_\kalg \rightarrow X \otimes_k \kalg$
tel que $f^\star T_X$ soit ample, où~$\kalg$ désigne une clôture algébrique de~$k$.  Si~$k$ est de caractéristique~$0$, il revient au même de demander que~$X$ soit rationnellement connexe.
En caractéristique $p>0$ la condition «~séparablement rationnellement connexe~» est strictement plus forte. Elle permet d'écarter certaines pathologies comme les variétés de type
général qui sont en même temps rationnellement connexes~\cite{shiodaexample}.} rationnellement connexe, sur un corps~$k$.
Si~$k$ est~$(C_1)$, a-t-on nécessairement $X(k)\neq\emptyset$~?
\end{question}

\bigskip
On sait que la réponse est affirmative si~$k$ est le corps des fonctions d'une courbe sur un corps algébriquement clos (d'après Graber, Harris, Starr et de~Jong, cf.~\cite{graberharrisstarr} et~\cite{dejongstarr}) ou
si~$k$ est un corps fini (d'après Esnault, cf.~§\ref{secesnault} ci-dessous).  Le~cas particulier de la question~\ref{laquestionc1} où~$k$
est l'extension non ramifiée maximale d'un corps $p$\nobreakdash-adique mérite d'être énoncé séparément:

\bigskip
\begin{question}
\label{questionnonram}
Toute variété propre, lisse et rationnellement connexe sur un corps $p$\nobreakdash-adique acquiert-elle
un point rationnel après une extension finie non ramifiée des scalaires~?
\end{question}

\bigskip
La question~\ref{questionnonram} est parallèle au théorème de Graber, Harris, Starr et de~Jong qui vient d'être mentionné (cf.~\cite[Théorème~7.7, Remarque~7.8]{ctcime}).

\bigskip
Nous venons de rencontrer avec la question~\ref{questionnonram} un premier problème ouvert concernant l'arithmétique des variétés
rationnellement connexes.  Les questions~\ref{laquestionc1} et~\ref{questionnonram} sont cependant en elles-mêmes d'une portée très limitée, cela autant à cause de
leurs hypothèses (beaucoup de corps ne sont pas~$(C_1)$, surtout parmi les corps qui présentent un intérêt du point de vue de l'arithmétique)
qu'à cause de la propriété qu'elles considèrent  (savoir que $X(k) \neq\emptyset$ ne renseigne en rien sur la structure de l'ensemble~$X(k)$,
sur les courbes rationnelles géométriquement intègres contenues dans~$X$, sur les groupes de Chow de~$X$, etc).

L'un des buts de la seconde moitié du chapitre~\ref{secundeux} sera de parvenir à la formulation de quelques grandes questions ouvertes au sujet de l'arithmétique des variétés rationnellement connexes.
Au~§\ref{subsecgroupebrauer} nous introduisons la notion de groupe de Brauer. Elle intervient à plusieurs endroits dans la suite du texte.
Au~§\ref{subsecobelem} nous définissons l'obstruction élémentaire, qui est une obstruction à l'existence de points rationnels sur un corps arbitraire,
et l'utilisons pour formuler, suivant de Jong et Starr, un analogue de la question~\ref{laquestionc1} sur le corps $\C(x,y)$.
Au~§\ref{subsecobbm}, consacré à l'obstruction de Brauer--Manin, nous arrivons à une première question ouverte sur les variétés rationnellement connexes définies sur un corps de nombres.
La $R$\nobreakdash-équivalence et les groupes de Chow de $0$\nobreakdash-cycles forment le sujet du~§\ref{subsecreq}.
Enfin, le~§\ref{subsecautre} rassemble quelques questions ouvertes qui n'ont pas trouvé leur place dans les paragraphes précédents.

\subsection{Groupe de Brauer}
\label{subsecgroupebrauer}

\medskip
La notion de groupe de Brauer sera utilisée aux paragraphes~\ref{subsecobelem} et~\ref{subsecobbm}
ainsi que dans~\cite{starrcevolume}.
Elle est étroitement liée à la propriété~$(C_1)$ ainsi qu'à l'étude des variétés de Severi--Brauer, qui constituent l'exemple
le plus simple de variétés rationnellement connexes.
Nous nous bornons ici à donner une définition et quelques exemples. 
Dans ce paragraphe, la plupart des preuves sont omises.
Le lecteur désireux d'en savoir plus consultera avec profit Platonov et
Yanchevski\u{\i}~\cite{platonovyanchevskii}.

Soient~$k$ un corps et $\kbar$ une clôture séparable de~$k$.  Le groupe de Brauer de~$k$ est un groupe abélien de torsion défini comme suit.

Une \emph{algèbre centrale simple} sur~$k$ est une $k$\nobreakdash-algèbre (unitaire) $A$ telle que la $\kbar$\nobreakdash-algèbre
$A \otimes_k \kbar$ soit isomorphe à l'algèbre de matrices $M_n(\kbar)$ pour un $n \geq 1$.  L'entier~$n$ est alors appelé le \emph{degré} de~$A$.
Les corps gauches de centre~$k$ et de dimension finie sur~$k$ sont des algèbres centrales simples sur~$k$.
Si~$A$ est une algèbre centrale simple sur~$k$, il existe un entier $r \geq 1$ et un corps gauche~$D$ de centre~$k$
tels que~$A$ soit isomorphe à la $k$\nobreakdash-algèbre~$D \otimes_k M_r(k)$.  Le corps gauche~$D$ est déterminé par~$A$ à un $k$\nobreakdash-isomorphisme près, on peut donc parler du corps gauche
associé à~$A$.
Deux algèbres centrales simples sur~$k$ sont \emph{semblables} si les corps gauches qui leur sont associés sont $k$\nobreakdash-isomorphes.

La relation $M_n(\kbar) \otimes M_m(\kbar) \simeq M_{nm}(\kbar)$ montre que les algèbres centrales simples sont stables par produit tensoriel.
De plus, si $A$, $A'$, $B$ et $B'$ sont des algèbres centrales simples sur~$k$, si~$A$ est semblable à~$A'$ et si~$B$ est semblable à~$B'$,
alors $A \otimes_k B$ est semblable à $A' \otimes_k B'$.  Par conséquent le produit tensoriel définit une opération sur les classes de similitude d'algèbres centrales simples sur~$k$.

Le \emph{groupe de Brauer} de~$k$, noté~$\Br(k)$, est par définition l'ensemble des classes de similitude d'algèbres centrales simples sur~$k$, muni du produit tensoriel.
L'élément neutre est représenté par la $k$\nobreakdash-algèbre~$k$.  Si~$A$ est une algèbre centrale simple sur~$k$, l'inverse de la classe de~$A$
dans $\Br(k)$ est représenté par l'algèbre opposée $A^\circ$ (par définition $A^\circ=A$ en tant que $k$\nobreakdash-espace vectoriel,
mais la multiplication de $A^\circ$ envoie $x,y \in A$ sur $yx \in A$); en effet $A \otimes_k A^\circ$ s'identifie à l'algèbre des endomorphismes du $k$\nobreakdash-espace vectoriel~$A$,
par l'application qui à $a \otimes b$ associe l'endomorphisme $x \mapsto axb$.

\bigskip
\begin{exemples*}
(i) Si~$k$ est séparablement clos, alors $\Br(k)=0$.

(ii) Si~$k$ est un corps fini, alors $\Br(k)=0$.  C'est le théorème de Wedderburn, habituellement énoncé ainsi: «~tout corps fini est commutatif~».

(iii) Les quaternions de Hamilton forment une algèbre centrale simple sur~$\R$.  Plus généralement, si~$k$ est un corps de caractéristique $\neq 2$ et si $a,b \in k^\star$,
la $k$\nobreakdash-algèbre $k \oplus ke \oplus kf \oplus kef$ définie par $e^2=a$, $f^2=b$, $fe=-ef$ est une algèbre centrale simple sur~$k$, de degré~$2$,
encore appelée \emph{algèbre de quaternions}.
C'est un corps gauche si et seulement si la conique d'équation homogène $ax^2 + by^2 = z^2$ n'admet pas de point $k$\nobreakdash-rationnel.

Notons $(a,b)_k \in \Br(k)$ la classe de l'algèbre de quaternions associée à $a,b \in k^\star$.
L'application $k^\star \times k^\star \rightarrow \Br(k)$, $(a,b)\mapsto(a,b)_k$ est $\Z$\nobreakdash-bilinéaire
symétrique et se factorise par $k^\star/k^{\star 2} \times k^\star/k^{\star 2}$.

(iv) Si~$\R$ est le corps des réels, on a $\Br(\R)=\Z/2$, engendré par la classe $(-1,-1)_\R$ des quaternions de Hamilton.
Cette propriété fut prouvée à l'origine par Frobenius mais elle résulte immédiatement de l'interprétation cohomologique du groupe de Brauer
(voir plus bas).

(v) Si~$k$ est un corps $p$\nobreakdash-adique, la théorie du corps de classes locale fournit un isomorphisme canonique $\Br(k) \isoto \Q/\Z$.

(vi) Si~$k$ est un corps de nombres, la théorie du corps de classes globale fournit une suite exacte
\begin{equation}
\label{suiteexcorpsglo}
\xymatrix@C=7ex{
0 \ar[r] & \Br(k) \ar[r] & \displaystyle\bigoplus_{v \in \Omega} \Br(k_v) \ar[r]^(.57){\sum \inv_v} & \Q/\Z \ar[r] & 0
}
\end{equation}
où~$\Omega$ désigne l'ensemble des places de~$k$ et~$k_v$ est le complété de~$k$ en~$v$, et où $\inv_v \colon \Br(k_v) \rightarrow \Q/\Z$
est l'isomorphisme évoqué en~(v) si~$v$ est une place finie et est l'inclusion évidente résultant de~(iv) ou de~(i) si~$v$ est infinie.
La flèche de gauche envoie la classe d'une $k$\nobreakdash-algèbre centrale simple~$A$ sur la famille des classes
des $k_v$\nobreakdash-algèbres centrales simples $A \otimes_k k_v$.
\end{exemples*}

\bigskip
Avec l'exemple des algèbres de quaternions, on voit que le groupe de Brauer contrôle l'existence de points rationnels sur les coniques.  Le phénomène est plus général.
On dit qu'une variété~$X$ sur~$k$ est une \emph{variété de Severi--Brauer} si $X \otimes_k \kbar$ est $\kbar$\nobreakdash-isomorphe à l'espace projectif
$\P^n_\kbar$ pour un $n \geq 0$.
Il y a une correspondance biunivoque entre les algèbres centrales simples sur~$k$ de degré~$n$ (à isomorphisme près)
et les variétés de Severi--Brauer sur~$k$ de dimension $n-1$ (à isomorphisme près).  Les algèbres de quaternions correspondent aux coniques projectives lisses;
l'algèbre de matrices $M_n(k)$ correspond à l'espace projectif $\P^{n-1}_k$.   On peut montrer qu'une variété de Severi--Brauer admettant un point rationnel
est automatiquement isomorphe à l'espace projectif.
Ainsi, à toute variété de Severi--Brauer~$X$ sur~$k$ est associée une classe de~$\Br(k)$, et cette classe est nulle
si et seulement si $X(k)\neq\emptyset$.

Bien entendu, les variétés de Severi--Brauer sont séparablement rationnellement connexes.  Pour que la question~\ref{laquestionc1}
soit raisonnable, il faudrait donc que toute variété de Severi--Brauer sur un corps~$(C_1)$ admette un point
rationnel.  C'est bien le cas:

\bigskip
\begin{proposition}
\label{brauercunnul}
Le groupe de Brauer d'un corps~$(C_1)$ est nul.
\end{proposition}

\bigskip
\begin{demo}
Soient~$k$ un corps~$(C_1)$ et~$D$ un corps gauche de centre~$k$ et de dimension finie sur~$k$.

Choisissons une clôture séparable~$\kbar$ de~$k$ et
un isomorphisme de $\kbar$\nobreakdash-algèbres $\phi \colon D \otimes_k \kbar \isoto M_n(\kbar)$.
L'application $D \otimes_k \kbar \rightarrow \kbar$, $x \mapsto \det(\phi(x))$
est polynomiale, au sens où dans une base du $k$\nobreakdash-espace vectoriel~$D$ elle s'écrit comme un polynôme en les coordonnées.
Il résulte du théorème de Skolem--Noether, selon lequel tout automorphisme d'une algèbre centrale simple est intérieur, qu'elle
ne dépend pas du choix de~$\phi$. Cela entraîne qu'elle est $\Gal(\kbar/k)$\nobreakdash-équivariante --- en effet, la conjuguer par $\sigma \in \Gal(\kbar/k)$ revient à remplacer $\phi$
par $\sigma \phi \sigma^{-1}$, c'est-à-dire à changer le choix de~$\phi$.
Elle provient donc par extension des scalaires d'une unique application polynomiale $D \rightarrow k$, appelée \emph{norme réduite}
et notée $\Nrd$.

Dans une base de~$D$ sur~$k$, la norme réduite est un polynôme homogène de degré~$n$ en~$n^2$ variables.  
Ce polynôme ne s'annule qu'en $0 \in D$ puisque~$D$ est un corps gauche et que $\Nrd(xy)=\Nrd(x)\Nrd(y)$ pour tous $x,y\in D$.
Le corps~$k$ étant~$(C_1)$, il s'ensuit que $n \geq n^2$, autrement dit $n=1$ et donc $D=k$.
\end{demo}

\bigskip
La proposition~\ref{brauercunnul} implique en particulier que le groupe de Brauer d'un corps fini ou de $\C(t)$ est nul.

Le groupe de Brauer $\Br(k)$ admet une interprétation cohomologique: il s'identifie au groupe de cohomologie galoisienne $H^2(k,\kbar^\star)$
(cohomologie du groupe profini $\Gal(\kbar/k)$ à valeurs dans le module galoisien discret~$\kbar^\star$).  Cela entraîne qu'il est de torsion.

Jointe à cette interprétation cohomologique, la proposition~\ref{brauercunnul} constitue le point de départ de la théorie de la dimension cohomologique des corps,
due à Tate (cf.~\cite[Ch.~II]{serrecg}).

\bigskip
Grothendieck~\cite{grothbrauer} a défini le groupe de Brauer (cohomologique) $\Br(X)$ d'un schéma~$X$ par la formule $\Br(X)=H^2_\et(X,\Gm)$.
Si~$k$ est un corps et que $X=\Spec(k)$, on retrouve le groupe $\Br(k)$ considéré ci-dessus.   Sous des hypothèses assez générales il existe une interprétation
des éléments de $\Br(X)$ en termes de familles d'algèbres centrales simples paramétrées par~$X$; nous n'en aurons cependant pas l'utilité.
Nous nous contentons de noter que si~$X$ est une variété intègre lisse sur un corps~$k$, alors $\Br(X)$ est naturellement
un sous-groupe de $\Br(k(X))$, où~$k(X)$ désigne le corps des fonctions de~$X$ (cf.~\cite[II, §1]{grothbrauer}).
Ainsi, par exemple, le groupe de Brauer d'une courbe lisse sur un corps algébriquement clos est nul, d'après le théorème~\ref{tsenlangnagata}
et la proposition~\ref{brauercunnul}.

\subsection{Obstruction élémentaire}
\label{subsecobelem}

\medskip
Dans ce paragraphe nous définissons l'obstruction élémentaire; c'est l'une des
rares obstructions à l'existence d'un point rationnel qui aient un sens pour des
variétés lisses arbitraires sur un corps arbitraire.
Cette obstruction intervient
dans la formulation des questions pertinentes concernant les
variétés «~rationnellement simplement connexes~» sur le corps $\C(x,y)$.

Soit~$X$ une variété lisse et géométriquement intègre sur un corps~$k$.  Soient $\kbar$ une clôture séparable de~$k$ et $\kbar(X)$ le corps des fonctions de $X \otimes_k \kbar$.

\bigskip
\begin{proposition}
Si $X(k)\neq\emptyset$, le morphisme de groupes $\kbar^\star \hookrightarrow \kbar(X)^\star$ admet une rétraction $\Gal(\kbar/k)$\nobreakdash-équivariante.
\end{proposition}

\bigskip
\begin{demo}
Soit $x \in X(k)$.  Notons $n=\dim(X)$.
Comme~$X$ est régulier, le complété de l'anneau local $\Orond_{X,x}$ est $k$\nobreakdash-isomorphe à $k[[\uplet{t_1}{t_n}]]$.  Le corps~$\kbar(X)$ se plonge donc
de manière $\Gal(\kbar/k)$\nobreakdash-équivariante
dans le corps de séries formelles itérées $\kbar((t_1))\cdots((t_n))$.
Pour conclure on remarque que le morphisme de groupes $\kbar^\star \hookrightarrow \kbar((t_1))\cdots((t_n))^\star$ admet une rétraction équivariante évidente:
celle obtenue en composant les~$n$ morphismes $\kbar((t_1))\cdots((t_i))^\star \rightarrow \kbar((t_1))\cdots((t_{i-1}))^\star$
qui à une série formelle en~$t_i$ et à coefficients dans $\kbar((t_1))\cdots((t_{i-1}))$ associent son coefficient non nul de plus bas degré.
\end{demo}

\bigskip
Si le morphisme de groupes $\kbar^\star \hookrightarrow \kbar(X)^\star$ n'admet pas de rétraction équivariante,
on dit, suivant Colliot-Thélène et Sansuc~\cite{ctsanobelem}, qu'il y a une \emph{obstruction élémentaire} à l'existence d'un $k$\nobreakdash-point sur~$X$.
Il ne peut y avoir d'obstruction élémentaire si~$k$ est un corps~$(C_1)$ (pour les variétés rationnellement connexes cela résulte de~\cite{ctsanobelem};
voir~\cite[Th.~3.4.1]{wittobelem} pour le cas général).

L'obstruction élémentaire est liée au problème de la descente de faisceaux inversibles sur $X \otimes_k \kbar$ qui, à isomorphisme près, sont invariants par $\Gal(\kbar/k)$:

\bigskip
\begin{proposition}
\label{proppicdescente}
Supposons que toute fonction inversible sur~$X \otimes_k \kbar$ soit constante (par exemple~$X$ propre).  Les flèches naturelles $\Pic(X) \rightarrow \Pic(X \otimes_k \kbar)$ et $\Br(k) \rightarrow \Br(X)$ s'inscrivent dans une suite exacte
canonique
$$
\xymatrix{
0 \ar[r] & \Pic(X) \ar[r] & \Pic(X \otimes_k \kbar)^{\Gal(\kbar/k)} \ar[r]^(0.65)\delta & \Br(k) \ar[r] & \Br(X) \rlap{\text{.}}
}
$$
Si l'obstruction élémentaire à l'existence d'un $k$\nobreakdash-point sur~$X$ s'évanouit,
alors $\Br(k) \rightarrow \Br(X)$ est injective, de sorte que $\Pic(X) \isoto \Pic(X\otimes_k \kbar)^{\Gal(\kbar/k)}$.
\end{proposition}

\bigskip
\begin{demo}
La première assertion résulte formellement de la suite spectrale de Hochschild--Serre
$H^p(k,H^q_\et(X \otimes_k \kbar,\Gm))\Rightarrow H^{p+q}_\et(X,\Gm)$ et du théorème~90 de Hilbert (qui affirme que $H^1(k,\kbar^\star)=0$).
Si l'inclusion $\kbar^\star \hookrightarrow \kbar(X)^\star$ admet une rétraction
$\Gal(\kbar/k)$\nobreakdash-équivariante, la flèche $H^2(k,\kbar^\star) \rightarrow H^2(k,\kbar(X)^\star)$ qu'elle induit en cohomologie admet aussi une rétraction et est donc injective.
Or d'une part $H^2(k,\kbar^\star)=\Br(k)$, et d'autre part $H^2(k,\kbar(X)^\star)$ s'injecte dans $\Br(k(X))$, comme il résulte de nouveau de la suite spectrale
de Hochschild--Serre et du théorème~90 de Hilbert.  Cela démontre la proposition puisque $\Br(X) \subset \Br(k(X))$.
\end{demo}

\bigskip
\begin{exemple*}
Pour comprendre la proposition~\ref{proppicdescente}, prenons pour~$X$ une conique projective lisse sur~$k$.
Dans ce cas $X \otimes_k \kbar$ est isomorphe à $\P^1_{\kbar}$
et l'application «~degré~» détermine donc un isomorphisme $\Pic(X \otimes_k \kbar) \isoto \Z$.  L'action de $\Gal(\kbar/k)$ sur $\Pic(X \otimes_k \kbar)$ est triviale.
L'inclusion $\Pic(X) \subset \Pic(X \otimes_k\kbar)^{\Gal(\kbar/k)}$ s'identifie à l'inclusion, dans~$\Z$, du sous-groupe engendré par les degrés des points fermés de~$X$; donc
à $\Z\subset \Z$ si $X(k)\neq\emptyset$ et à $2\Z \subset \Z$ sinon.  L'image de $1 \in \Z=\Pic(X\otimes_k\kbar)^{\Gal(\kbar/k)}$
par~$\delta$ est la classe dans $\Br(k)$ de la conique~$X$ vue comme variété de Severi--Brauer.
\end{exemple*}

\bigskip
C'est \emph{via} la proposition~\ref{proppicdescente} que l'obstruction
élémentaire joue un rôle dans l'étude des variétés «~rationnellement simplement
connexes~».
La définition précise de la simple connexité rationnelle est en cours
d'élaboration (cf.~\cite[\textsection1]{dejonghestarr}).  Cette notion
est à la connexité rationnelle ce qu'en topologie la simple connexité est à la
connexité.  L'espoir formulé par de~Jong et Starr dans~\cite{dejonghestarr} est que
sur le corps $\C(x,y)$, toute variété «~rationnellement simplement connexe~» pour laquelle
l'obstruction élémentaire s'évanouit admette un point rationnel.
Un tel énoncé serait l'analogue du théorème de Graber, Harris et Starr~\cite{graberharrisstarr}
selon lequel sur le corps~$\C(x)$, toute variété rationnellement connexe admet un point rationnel.
Sur le corps~$\C(x)$, qui est~$(C_1)$, l'obstruction élémentaire ne joue pas de rôle; mais sur $\C(x,y)$
elle ne peut être ignorée puisque même dans le cas des variétés dont la géométrie est la plus simple
(les variétés de Severi--Brauer), elle ne s'annule pas en général.

Citons pour terminer l'étude systématique
entreprise par Borovoi, Colliot-Thélène et Skorobogatov~\cite{boctsko} de
l'obstruction élémentaire pour les espaces homogènes de groupes algébriques.

\subsection{Obstruction de Brauer--Manin}
\label{subsecobbm}

\medskip
Les phénomènes globaux qui rendent l'arithmétique des corps de nombres plus complexe que celle des corps
finis, des corps $p$\nobreakdash-adiques, ou que celle des corps de fonctions
d'une variable sur un corps algébriquement clos, jouent bien sûr aussi un rôle important dans l'arithmétique
des variétés rationnellement connexes définies sur un corps de nombres.

Ces phénomènes compliquent la théorie mais ils constituent en même temps des outils pour la comprendre.

Ainsi, en combinant la notion de groupe de Brauer d'une variété (introduite par Grothendieck) avec la théorie du corps de classes globale
(et plus précisément la loi de réciprocité), Manin~\cite{maninicm} a défini, pour toute variété~$X$ sur un corps de nombres~$k$,
un ensemble qui joue le rôle d'une «~première approximation~» de l'ensemble~$X(k)$ des points rationnels de~$X$.
Plus précisément,
supposons la variété~$X$ propre et lisse et notons $X(\A_k)$ l'espace $\prod_{v \in \Omega} X(k_v)$, où~$\Omega$ désigne l'ensemble
des places de~$k$ et~$k_v$ le complété de~$k$ en~$v$.  Manin considère le sous-ensemble $X(\A_k)^{\Br} \subset X(\A_k)$ constitué des
familles $(P_v)_{v \in \Omega} \in X(\A_k)$ telles que pour tout $A \in \Br(X)$, la somme $\sum_{v \in \Omega} \inv_v A(P_v)$
soit égale à $0 \in \Q/\Z$.  (On montre que cette somme est finie.)  Ici $A(P_v) \in \Br(k_v)$ désigne l'évaluation de~$A$ en~$P_v$ et $\inv_v \colon \Br(k_v) \hookrightarrow \Q/\Z$
est l'application apparaissant dans la suite exacte~(\ref{suiteexcorpsglo}).  Comme~(\ref{suiteexcorpsglo}) est un complexe, l'ensemble $X(k)$, qui est naturellement un sous-ensemble
de $X(\A_k)$, est même inclus dans $X(\A_k)^\Br$.  En particulier, si $X(\A_k)^\Br=\emptyset$, alors $X(k)=\emptyset$; on dit dans ce cas qu'il y a une \emph{obstruction de Brauer--Manin à
l'existence d'un $k$\nobreakdash-point sur~$X$}.

Munissons~$X(\A_k)$ de la topologie produit des topologies $v$\nobreakdash-adiques.
Le sous-ensemble $X(\A_k)^{\Br}$ est fermé dans $X(\A_k)$.  Par conséquent, si $X(\A_k)^{\Br} \neq X(\A_k)$, l'ensemble $X(k)$ ne peut être dense dans $X(\A_k)$.
On dit dans ce cas qu'il y a une \emph{obstruction de Brauer--Manin à l'approximation faible sur~$X$}.

Y a-t-il d'autres phénomènes que l'obstruction de Brauer--Manin qui peuvent empêcher qu'une variété rationnellement connexe sur un corps de nombres possède un point
rationnel ou encore satisfasse à l'approximation faible~?  La question est ouverte:

\bigskip
\newcommand{\refct}{\cite[p.~174]{ctbudapest}}
\begin{question}[ (Colliot-Thélène~\refct)]%
\label{questionct}
Soient~$k$ un corps de nombres et~$X$ une variété propre, lisse, rationnellement connexe, sur~$k$.  L'ensemble $X(\A_k)^\Br$ coïncide-t-il nécessairement avec l'adhérence
de $X(k)$ dans $X(\A_k)$~?
\end{question}

\bigskip
Que cette question admette une réponse affirmative dans le cas des surfaces est une conjecture de Colliot-Thélène et Sansuc.
De très nombreux travaux ont été consacrés à ce problème, tant pour les surfaces rationnelles que pour d'autres types de variétés rationnellement connexes (notamment
les compactifications lisses d'espaces homogènes de groupes algébriques linéaires ou encore les intersections complètes de bas degré dans~$\P^n$).
Nous renvoyons à~\cite{peyresembour} pour un survol et pour de nombreuses références supplémentaires à la littérature,
et à~\cite[§4]{cttoulouse} pour une conjecture plus précise dans le cas des variétés rationnelles.

Notons qu'une réponse affirmative à la question~\ref{questionct} impliquerait que tout groupe fini est groupe de Galois sur~$\Q$ (cf.~\cite[§3.5]{serretopics}).

\subsection{Groupe de Chow, $R$-équivalence}
\label{subsecreq}

\medskip
Soit~$X$ une variété sur un corps~$k$.

Deux points rationnels $x,y \in X(k)$ sont \emph{directement $R$\nobreakdash-équivalents} s'il existe une application
rationnelle $\phi \colon \P^1_k \dashrightarrow X$ définie en~$0$ et en~$\infty$ telle que $\phi(0)=x$ et $\phi(\infty)=y$.
La $R$\nobreakdash-équivalence est par définition la plus petite relation d'équivalence sur~$X(k)$ contenant la $R$\nobreakdash-équivalence directe.
Cette notion fut introduite par Manin~\cite{manincubicforms}.
L'ensemble des classes de $R$\nobreakdash-équivalence est noté $X(k)/R$.

Un \emph{$0$\nobreakdash-cycle} sur~$X$ est un élément du $\Z$\nobreakdash-module libre de base l'ensemble des points fermés de~$X$.
Deux $0$\nobreakdash-cycles sont \emph{rationnellement équivalents} si leur différence est une combinaison linéaire
de cycles de la forme $f_\star(z)$, où $f \colon C \rightarrow X$ est un morphisme propre d'une courbe normale connexe
vers~$X$ et~$z$ est un diviseur principal sur~$C$.
Le \emph{groupe de Chow des $0$\nobreakdash-cycles sur~$X$} est le quotient $\CH_0(X)$ du groupe des $0$\nobreakdash-cycles sur~$X$ par le sous-groupe
des $0$\nobreakdash-cycles rationnellement équivalents à~$0$.  Si~$X$ est propre sur~$k$, il y a une flèche «~degré~» de $\CH_0(X)$ vers~$\Z$; son noyau est le \emph{groupe de Chow
des $0$\nobreakdash-cycles de degré~$0$ sur~$X$}, noté~$\CHz(X)$.

Deux points $R$\nobreakdash-équivalents sont rationnellement équivalents.  En conséquence,
l'étude du quotient $X(k)/R$ est liée à celle du groupe $\CHz(X)$.

Il existe plusieurs grandes questions ouvertes concernant l'ensemble $X(k)/R$ et le groupe $\CHz(X)$ pour les variétés séparablement rationnellement connexes.
Nous énonçons dans ce paragraphe les plus significatives d'entre elles et renvoyons le lecteur au rapport de Colliot-Thélène~\cite[§10--§11]{ctcime}
pour une discussion plus approfondie de ces questions, de questions connexes, et des résultats connus.

\bigskip
\begin{question}[ (Colliot-Thélène~\cite{ctcime})]%
\label{questionctqac}
Soient~$k$ un corps~$(C_1)$ et~$X$ une variété propre, lisse, séparablement rationnellement connexe, sur~$k$.
A-t-on nécessairement $\Card(X(k)/R)\leq 1$ et $\CHz(X)=0$~?
\end{question}

\bigskip
Comme Colliot-Thélène le remarque dans~\cite{ctcime}, si l'inégalité
$\Card(X(k)/R)\leq 1$ était vérifiée pour toute surface rationnelle~$X$ sur le corps $k=\C(t)$, les variétés complexes de dimension~$3$ fibrées en coniques
sur~$\P^2$ seraient toutes unirationnelles (ce qui est peu vraisemblable).
La question mérite néanmoins d'être posée, au vu de tous les résultats connus allant dans le sens d'une réponse positive.
Entre autres, on a bien $\CHz(X)=0$ si~$k$ est un corps fini (Kato--Saito) ou si~$X$ est une surface rationnelle (Colliot-Thélène), et l'on a bien
$\Card(X(k)/R)\leq 1$ si~$k$ est fini et de cardinal assez grand en un sens précis (Kollár--Szabó; ce théorème fera l'objet du~§\ref{seckollsz}).
Voir~\cite{ctcime} pour d'autres résultats et pour des références à la littérature.

Si maintenant~$k$ est un corps $p$\nobreakdash-adique ou un corps de nombres, la mauvaise réduction a tendance à s'opposer à ce que $\Card(X(k)/R)\leq 1$ ou $\CHz(X)=0$.
Kollár a néanmoins démontré que $X(k)/R$ est un ensemble fini si~$X$ est une variété lisse rationnellement connexe sur un corps $p$\nobreakdash-adique~$k$.
Ce théorème fera l'objet du~§\ref{seckoll}.  On peut de même espérer que les questions suivantes admettent une réponse positive:

\bigskip
\begin{questions}[ (Colliot-Thélène~\cite{cttoulouse}, \cite{ctbordeaux}, \cite{ctcime})]%
\label{questionctfinitude}
Le groupe $\CHz(X)$ est-il fini si~$X$ est une variété propre, lisse, rationnellement connexe sur un corps $p$\nobreakdash-adique~?
L'ensemble $X(k)/R$ et le groupe $\CHz(X)$ sont-ils finis si~$X$ est une variété propre, lisse, rationnellement connexe sur un corps de nombres~?
\end{questions}

\bigskip
Nous renvoyons de nouveau à~\cite{ctcime} pour une discussion des cas connus (voir
aussi~\cite[§5]{ctfinitudechow}).  Si~$k$ est un corps de nombres, la question de la finitude de $X(k)/R$ est ouverte même pour les surfaces
rationnelles.

Si~$k$ est un corps de type fini sur son sous-corps premier et~$X$ est une variété propre, lisse, séparablement rationnellement connexe, sur~$k$,
on ne sait pas si $\CHz(X)$ est toujours fini ou non (cf.~\cite[§11]{ctcime}).  En revanche Kollár~\cite{kollsp} a fabriqué des exemples montrant que $X(k)/R$ peut être infini,
avec $k=\Q(t)$.

\subsection{Quelques autres questions sur l'arithmétique des variétés rationnellement connexes}
\label{subsecautre}

\medskip
La conjecture suivante est un cas particulier d'une conjecture de Campana qui caractérise de façon purement géométrique, parmi les variétés définies sur un corps de nombres,
celles dont les points rationnels sont «~potentiellement denses~».

\bigskip
\newcommand{\campanaref}{\cite[Conjecture~9.20]{campanaconj}}
\begin{conjecture}[ (Campana~\campanaref)]%
\label{conjcampana}
Soient~$k$ un corps de nombres et $X$ une variété propre, lisse, rationnellement connexe, sur~$k$.
Il existe une extension finie $\ell/k$ telle que $X(\ell)$ soit dense dans $X \otimes_k \ell$ pour la topologie de Zariski.
\end{conjecture}

\bigskip
On ignore même:

\bigskip
\begin{question}
\label{questiondensite}
Soit~$X$ une variété propre, lisse, séparablement rationnellement connexe sur un corps infini~$k$.
Si l'ensemble $X(k)$ n'est pas vide, est-il dense dans~$X$ pour la topologie de Zariski~?
\end{question}

\bigskip
Il n'y a pas de raison d'espérer que dans cette généralité, la réponse à cette question soit positive.
Cependant, il en va ainsi lorsque~$k$ est le corps des fonctions d'une courbe sur un corps algébriquement clos,
du moins si la variété~$X$ est projective (Kollár--Miyaoka--Mori~\cite[IV.6.10]{kollrational}).
Pour $k=\Q$, la question~\ref{questiondensite} est ouverte même dans le cas des surfaces rationnelles.

Rosenlicht~\cite[p.~46]{rosenlicht} donne des exemples de corps~$k$ infinis et de courbes affines lisses rationnelles~$X$
telles que $X(k)$ soit fini et non vide.
(Bien entendu, le corps~$k$ est imparfait et la compactification régulière de~$X$ n'est pas lisse sur~$k$.)
L'hypothèse de propreté est donc essentielle dans la question~\ref{questiondensite}.

Une réponse affirmative à la question~\ref{questionct} entraînerait une réponse affirmative à la question~\ref{questiondensite} pour les corps de nombres,
et donc impliquerait la validité de la conjecture~\ref{conjcampana}.  La conjecture~\ref{conjcampana} est connue pour toutes les variétés de Fano lisses de dimension~$3$
à l'exception des revêtements doubles de~$\P^3$ ramifiés le long d'une surface sextique (Bogomolov, Harris, Tschinkel; cf.~\cite{bogotschinkel}) mais elle est ouverte
pour les variétés de dimension~$3$ fibrées en coniques sur~$\P^2$.

Jusqu'ici nous avons seulement considéré des questions de nature qualitative
concernant les points rationnels et les $0$\nobreakdash-cycles sur les variétés
rationnellement connexes. 
Sur un corps de nombres, les~aspects quantitatifs donnent
également lieu à des problèmes
très intéressants (estimation asymptotique du nombre de points rationnels de hauteur bornée, 
interprétation des constantes qui apparaissent).
Leur étude fut entamée par Manin à la fin des années~1980.
Nous renvoyons à Peyre~\cite{peyrebourmanin} pour un survol des conjectures et résultats connus en~2001.

Signalons enfin deux problèmes qui, à défaut d'être de nature arithmétique,
sont directement inspirés de préoccupations arithmétiques: la question de
l'approximation faible pour les variétés rationnellement connexes définies sur
le corps des fonctions d'une courbe complexe, traitée en détail dans~\cite{hassettcevolume} et~\cite{hassetttschinkel};
et un problème analogue pour les variétés
rationnellement connexes sur le corps des fonctions d'une courbe réelle.  Par
exemple on peut demander: si~$C$ est une courbe propre et lisse sur le
corps~$\R$ des réels, si $f \colon X \rightarrow C$ est un morphisme propre de
fibre générique géométrique lisse, rationnellement connexe (donc connexe), et si
l'application $X(\R) \rightarrow C(\R)$ induite par~$f$ admet une section
continue à valeurs dans l'ouvert de lissité de~$f$, alors le morphisme~$f$
admet-il une section~?  Ce type de questions, qui trouve son origine dans les
travaux de Witt, est abordé par Ducros~\cite{ducros}, auquel nous renvoyons le
lecteur pour davantage de références à la littérature.

\section{Variétés rationnellement connexes sur les corps fertiles}
\label{seckoll}

Dans ce chapitre nous exposons les résultats obtenus par Kollár~\cite{kollarloc} en~1999 au sujet de l'arithmétique des variétés rationnellement connexes sur
des corps tels que~$\Qp$ ou~$\R$, puis nous en évoquons quelques raffinements ultérieurs (dus à Kollár~\cite{kollsp}).
Ici entrent en jeu les techniques de déformation de courbes rationnelles; pour la première fois elles sont appliquées
dans des situations où le corps de base n'est pas algébriquement clos.

\subsection{Énoncé du théorème principal; conséquences}

\medskip
Soient~$k$ un corps, $\kbar$ une clôture algébrique de~$k$ et~$X$ une variété projective, lisse et séparablement rationnellement connexe, sur~$k$.

Que la variété~$X$ soit séparablement rationnellement connexe signifie que $X \otimes_k \kbar$ contient une courbe rationnelle très libre,
c'est-à-dire qu'il existe un morphisme $f \colon \P^1_\kbar \rightarrow X \otimes_k \kbar$ tel que $f^\star T_X$ soit ample\footnote{Rappelons
que $f^\star T_X$ est isomorphe (comme tout faisceau de modules
localement libre de type fini sur~$\P^1$) à une somme directe $\bigoplus_{i=1}^r
\Orond(a_i)$ pour des $a_i \in \Z$; il est ample si et seulement si $a_i>0$ pour tout~$i$.}.
D'après Kollár, Miyaoka et Mori, pour tout point de $X \otimes_k \kbar$ il existe une courbe rationnelle très libre sur $X \otimes_k \kbar$ passant par ce point.
Lorsque~$k$ n'est pas algébriquement clos, on peut se demander si la même chose est vraie sans étendre les scalaires de~$k$ à~$\kbar$:

\bigskip
\begin{question}
\label{questionexisteample}
Soit $x \in X(k)$.
Existe-t-il un morphisme $f \colon \P^1_k \rightarrow X$ tel que $f^\star T_X$ soit ample et que $f(0)=x$~?
\end{question}

\bigskip
C'est là une question ouverte (cf.~\cite[§9]{araujokollar}).  Il est peu vraisemblable qu'elle admette une réponse affirmative en toute généralité (notamment
pour~$k$ fini ou pour $k=\C(t)$).
Le théorème principal de ce chapitre (théorème~\ref{corpslocauxkollarth} ci-dessous) affirme néanmoins que tel est le cas lorsque le corps~$k$ est \emph{fertile}.

\bigskip
\begin{definition}[ (Pop~\cite{pop})]%
\label{deffertile}
Le corps~$k$ est dit \emph{fertile} («~large~» en anglais) si pour toute variété~$X$ connexe, lisse sur~$k$ et telle que $X(k)\neq\emptyset$,
l'ensemble $X(k)$ est dense dans~$X$ pour la topologie de Zariski.
\end{definition}

\bigskip
Dans la définition~\ref{deffertile}, la variété~$X$ n'est pas supposée rationnellement connexe.
On peut vérifier que le corps~$k$ est fertile si et seulement si pour toute courbe~$C$ lisse sur~$k$ telle que $C(k)\neq\emptyset$,
l'ensemble~$C(k)$ est infini.

\bigskip
\begin{exemples*}
(i) Tout corps local usuel (c'est-à-dire les corps~$\Qp$, $\Fp((t))$, $\R$ et leurs extensions finies) est fertile.  Cela résulte du théorème d'inversion
locale pour les variétés analytiques~\cite[Part~II, Th.~III.9.2]{serrelie}.

(ii) Tout corps réel clos et tout corps $p$\nobreakdash-adiquement clos (\emph{i.e.} tout corps dont le groupe de Galois absolu est isomorphe à celui de~$\R$ ou d'un corps $p$\nobreakdash-adique;
par exemple la fermeture algébrique de~$\Q$ dans~$\Qp$) est fertile (cf.~\cite[Theorem~7.8]{prestelroquette}).

(iii) Tout corps pseudo-algébriquement clos (voir la définition~\ref{defpac}) est fertile.

(iv) Tout corps dont le groupe de Galois absolu est un pro-$p$\nobreakdash-groupe est fertile\footnote{Lorsque le corps est parfait, cette remarque est due à Colliot-Thélène~\cite{ctann}.
Il s'avère que l'argument de~\cite{ctann} s'applique aussi lorsque~$k$ n'est pas parfait, compte tenu que pour toute courbe lisse~$X$ sur~$k$ et tout ouvert dense $U \subset X$,
tout diviseur sur~$U$ est linéairement équivalent, sur~$X$, à un diviseur dont le support est lisse et inclus dans~$U$ (cf.~\cite[Lemma~3.16]{voevodsky}).}.

(v) Le corps des fractions de tout anneau local hensélien intègre est fertile (cf.~\cite[Theorem~1.1]{poplarge}).
En particulier, quels que soient le corps~$k$ et l'entier $n \geq 1$, le corps $k((\uplet{x_1}{x_n}))$ est fertile.

(vi) Si~$k$ est un corps global et~$S$ un ensemble fini de places de~$k$, la plus grande extension algébrique de~$k$ dans laquelle les places de~$S$ sont totalement décomposées
est un corps fertile (cf.~\cite[2.4.3]{mbconstruction}).  En particulier le corps des nombres algébriques totalement réels est fertile.
\end{exemples*}

\bigskip
\begin{theoreme}[ (Kollár~\cite{kollarloc})]%
\label{corpslocauxkollarth}
Soient~$k$ un corps fertile et~$X$ une variété projective, lisse et séparablement rationnellement connexe, sur~$k$.  Soit $x \in X(k)$.
Il existe un morphisme $f \colon \P^1_k \rightarrow X$ tel que $f(0)=x$
et que $f^\star T_X$ soit ample.
\end{theoreme}

\bigskip
Nous consacrerons le §\ref{seckollocpreuve} à la preuve du théorème~\ref{corpslocauxkollarth}.
Son corollaire le plus frappant (que nous prouverons au §\ref{seckolloccorpreuve}) concerne la $R$\nobreakdash-équivalence (cf.~§\ref{subsecreq}):

\bigskip
\begin{corollaire}
\label{kolloccorreq}
Soit~$k$ un corps $p$\nobreakdash-adique ou le corps~$\R$ des réels.  Soit~$X$ une variété projective, lisse et rationnellement connexe, sur~$k$.
L'ensemble $X(k)/R$ est fini et les classes de $R$\nobreakdash-équivalence sont des parties ouvertes et fermées de~$X(k)$.  En particulier,
si $k=\R$, les classes de $R$\nobreakdash-équivalence sont exactement les composantes connexes de~$X(\R)$.
\end{corollaire}

\bigskip
La finitude de $X(k)/R$ lorsque~$k$ est un corps local
n'était auparavant connue que dans un certain nombre de cas particuliers (les surfaces rationnelles, à l'exception des surfaces de del Pezzo de degré~$1$ ou~$2$;
les intersections lisses de deux quadriques dans~$\P^n_k$ pour $n\geq 4$; les hypersurfaces cubiques lisses; les compactifications de groupes algébriques
linéaires).

En revanche le théorème~\ref{corpslocauxkollarth} ne dit rien sur la finitude de $\CHz(X)$.

Un autre corollaire notable du théorème~\ref{corpslocauxkollarth} concerne les variétés unirationnelles.  Une variété~$X$ unirationnelle et lisse sur un corps~$k$ est-elle $k$\nobreakdash-unirationnelle dès que $X(k)\neq\emptyset$~?
Cette question est ouverte.
Il semble probable que la réponse soit négative en général; néanmoins
le théorème~\ref{corpslocauxkollarth} permet d'apporter une réponse positive pour certaines variétés~$X$ lorsque~$k$ est fertile.
Le cas le plus simple est celui des surfaces fibrées en coniques (cf.~§\ref{subsecinterlude}), qui au moins sur les corps $p$\nobreakdash-adiques avait déjà été établi par Yanchevski\u{\i}~\cite{yanchevskiigeneve}:

\bigskip
\begin{corollaire}[ (Yanchevski\u{\i})]%
\label{coryanch}
Soit~$X$ une surface projective, lisse, rationnelle, fibrée en coniques, sur un corps~$k$.
Supposons~$k$ fertile. Alors~$X$ est $k$\nobreakdash-unirationnelle si et seulement si $X(k)\neq\emptyset$.
\end{corollaire}

\bigskip
\begin{demo}
Si~$X$ est $k$\nobreakdash-unirationnelle alors $X(k)\neq\emptyset$ de façon évidente.  Réciproquement, supposons~$X(k)$ non vide.
Notons $\pi \colon X \rightarrow \P^1_k$ une fibration en coniques sur~$X$.
D'après le~théorème~\ref{corpslocauxkollarth}, il existe un morphisme $f \colon \P^1_k \rightarrow X$ tel que $f^\star T_X$ soit ample.
Le~composé $\pi \circ f \colon \P^1_k \rightarrow \P^1_k$ est fini (sinon l'image de~$f$ serait une courbe rationnelle contenue dans une fibre de~$\pi$; elle ne pourrait
donc pas être très libre).
Soit $X'$ le produit fibré de $\pi \colon X \rightarrow \P^1_k$ avec $\pi \circ f \colon \P^1_k \rightarrow \P^1_k$.
Le morphisme~$f$ induit une section de
la seconde projection $\pi' \colon X' \rightarrow \P^1_k$;
la fibre générique de~$\pi'$, qui est une conique lisse, admet donc un point rationnel.
Par conséquent elle est $k(\P^1_k)$\nobreakdash-rationnelle et la surface~$X'$ est $k$\nobreakdash-rationnelle.
Comme~$X'$ domine~$X$, la $k$\nobreakdash-unirationalité de~$X$ s'ensuit.
\end{demo}

\bigskip
Le même argument permet d'établir la conclusion du corollaire~\ref{coryanch} pour quelques autres
classes de variétés unirationnelles que les surfaces fibrées en coniques au-dessus d'une conique.
Voir \cite[Corollary~1.8]{kollarloc}.

Pour~$k=\Q$, la question de la $k$\nobreakdash-unirationalité des surfaces rationnelles~$X$ fibrées en coniques telles que $X(k)\neq\emptyset$
est une question ouverte.  Pour ces surfaces, même la densité de $X(k)$ dans~$X$ pour la topologie de Zariski est inconnue.

\bigskip
\begin{remarque}
Soit~$X$ une variété projective, lisse et rationnellement connexe sur~$\R$.
Soit~$C$ la conique réelle sans point réel.  Si $X(\R)\neq\emptyset$, le
théorème~\ref{corpslocauxkollarth} entraîne l'existence d'un morphisme non constant
$C \rightarrow X$ (composer le morphisme donné par le théorème~\ref{corpslocauxkollarth}
avec un morphisme dominant $C \to \P^1_\R$).  La question de savoir si un tel morphisme existe toujours est
ouverte (cf.~\cite[Remarks~20]{araujokollar}).  Une réponse affirmative à cette
question résulterait de la conjecture~\ref{conjlangreel}
et d'une réponse affirmative à la question~\ref{laquestionc1}.
\end{remarque}

\subsection{Preuve du théorème~\ref{corpslocauxkollarth}}
\label{seckollocpreuve}

\medskip
\subsubsection{Esquisse de l'argument}

Le principe de la preuve est le suivant.

Soit~$X$ une variété projective, lisse, séparablement rationnellement connexe et de dimension~$\geq 1$, sur un corps fertile~$k$.
Soit $x \in X(k)$.  Soit~$\kbar$ une clôture séparable de~$k$.  Au moins si~$k$ est parfait, on sait qu'il existe une
courbe rationnelle (irréductible) très libre $B \subset X \otimes_k \kbar$
contenant~$x$.  On aimerait en exhiber une qui soit définie sur~$k$.
Pour cela faisons agir le groupe de Galois $\Gal(\kbar/k)$ sur~$B$; on obtient ainsi un nombre fini de courbes conjuguées $\uplet{B_1}{B_n}$.
Leur réunion $C = \bigcup_{i=1}^n B_i$ est une courbe connexe contenue dans~$X$ et passant par~$x$.  Elle définit donc un point $k$\nobreakdash-rationnel $[C]$
du schéma de Hilbert $\Hilb(X,x)$ paramétrant les sous-$k$\nobreakdash-schémas fermés de~$X$ contenant~$x$.

Supposons un instant que $\Hilb(X,x)$ soit lisse sur~$k$ en~$[C]$ et que le point générique de la composante irréductible~$I$ de $\Hilb(X,x)$ contenant~$[C]$
corresponde à une courbe rationnelle sur $X \otimes_k k(I)$.  Comme~$k$ est fertile et que la variété~$I$ possède un point rationnel lisse
(à savoir~$[C]$), l'ensemble $I(k)$ est alors dense dans~$I$ pour la topologie de Zariski.  Un point de $I(k)$ situé en dehors d'un certain fermé strict de~$I$ fournit ainsi
une courbe rationnelle $D \subset X$ contenant~$x$; un petit raffinement de cet argument permet d'imposer à~$D$ d'être très libre.

Cette esquisse est approximative mais elle contient néanmoins l'idée centrale de la démonstration: si l'on dispose d'un espace de modules~$M$ paramétrant un certain type d'objets
sur un corps fertile~$k$ (par exemple des courbes connexes de genre arithmétique~$0$) et si l'on cherche à fabriquer un point rationnel
d'un ouvert $M^0 \subset M$ (par exemple~$M^0$ pourrait être l'ensemble des points de~$M$ qui correspondent à des courbes irréductibles), il suffit de fabriquer un point rationnel lisse de~$M$ situé
sur la frontière de~$M^0$.

Il y a deux raisons pour lesquelles l'approche \emph{via} $\Hilb(X,x)$ fonctionne mal.
La~première est que l'étude infinitésimale de $\Hilb(X,x)$ en~$[C]$
est d'autant plus compliquée que la courbe~$C$ est singulière en~$x$; or on ne dispose d'aucun contrôle sur la singularité de~$C$ en~$x$ (même si les courbes~$B_i$ sont lisses,
leurs intersections deux à deux en~$x$ ne sont pas nécessairement transverses).
Vérifier la lissité de $\Hilb(X,x)$ en~$[C]$ ne saurait donc être une formalité.
La seconde est que la courbe~$C$ n'est pas toujours de genre arithmétique nul
(par exemple les~$B_i$ pourraient se rencontrer ailleurs qu'en~$x$); il n'y a donc pas de raison d'espérer qu'elle se déforme
(dans une famille plate) en une courbe rationnelle.

La manière la plus simple de contourner ces deux obstacles est de chercher à déformer non pas la courbe~$C$ mais une courbe réductible plus simple (un peigne),
dont les singularités et le genre arithmétique sont contrôlés.  La contrepartie est que la courbe n'est plus plongée dans~$X$: elle est seulement munie d'un morphisme (non
fini) vers~$X$.
Concrètement, lorsqu'on dispose d'un peigne et qu'on souhaite le déformer, la méthode la plus souple est de considérer l'espace de modules des courbes stables de Kontsevich
(cf.~\cite[§8]{araujokollar}).  Néanmoins, pour prouver le théorème~\ref{corpslocauxkollarth}, une construction un peu artificielle mais techniquement beaucoup plus économique sera suffisante
(comparer les arguments de~\cite[§4]{araujokollar} et de~\cite[§10]{araujokollar}).
C'est la voie que
nous suivrons (tout comme Kollár~\cite{kollarloc} mais contrairement à Araujo et Kollár~\cite[§9]{araujokollar}).

\subsubsection{La preuve proprement dite}
\label{chap2preuveproprementdite}

Soient~$k$ un corps fertile et~$X$ une variété projective et lisse sur~$k$,
séparablement rationnellement connexe, munie d'un point $x \in X(k)$.
Notons~$\kbar$ une clôture séparable de~$k$.

\bigskip
\begin{definition}
\label{defpeigne}
Un \emph{peigne} sur~$k$ est une courbe~$C$ sur~$k$, propre et réduite, munie d'une composante irréductible distinguée $C_0 \subset C$,
telle que les conditions suivantes soient remplies:
\begin{itemize}
\item la courbe~$C_0$ (dite \emph{manche} de~$C$) est lisse et géométriquement connexe sur~$k$;
\item les composantes irréductibles de $C \otimes_k \kbar$ autres que $C_0 \otimes_k \kbar$ (dites \emph{dents} de~$C$)
sont des courbes rationnelles lisses deux à deux disjointes qui rencontrent chacune
$C_0 \otimes_k \kbar$ en un unique point, l'intersection en ce point étant de plus transverse.
\end{itemize}
\end{definition}

\bigskip
Nous allons construire un peigne à partir d'une courbe définie sur une extension finie galoisienne de~$k$.

\bigskip
\begin{lemme}
Il existe une extension finie galoisienne $\ell/k$ et un $\ell$\nobreakdash-morphisme très libre $f_\ell \colon \P^1_\ell \rightarrow X \otimes_k \ell$ tel que $f_\ell(0)=x$.
\end{lemme}

\bigskip
\begin{demo}
L'espace de modules
$\Hom^0(\P^1_k,X\;\!;0\mapsto x)$
des morphismes très libres de~$\P^1_k$ vers~$X$ envoyant~$0$ sur~$x$ est un $k$\nobreakdash-schéma lisse
(cf.~\cite[2.6]{debarrelivre})
et non vide (puisqu'il admet un point dans une clôture algébrique de~$k$,
cf.~\cite[IV.3.9]{kollrational}).
Par conséquent il admet un point dans une extension finie séparable de~$k$; d'où le lemme.
\end{demo}

\bigskip
Fixons~$\ell$ et~$f_\ell$ comme dans le lemme.  Soit $M \in \A^1_k$ un point fermé de corps résiduel~$\ell$.
Soit $C \subset \P^1_k \times \P^1_k$ la réunion de $\{0\} \times \P^1_k$ et de~$\P^1_k \times M$.
La courbe~$C$ est un peigne sur~$k$, de manche $C_0 = \{0\} \times \P^1_k$.  L'action de $\Gal(\ell/k)$
sur l'ensemble des dents de $C \otimes_k \kbar$ est simplement transitive.

Le morphisme $\{0\} \times \P^1_k \rightarrow X$ constant égal à~$x$ et le morphisme \mbox{$f_\ell \colon \P^1_k \times M \rightarrow X$} coïncident
sur $\{0\} \times M$ et se recollent donc pour former un $k$\nobreakdash-morphisme $f \colon C \rightarrow X$.
C'est ce peigne~$C$, muni de ce morphisme $f \colon C \rightarrow X$, que l'on souhaite maintenant déformer en un morphisme de~$\P^1_k$ vers~$X$.
Pour cela on commence par déformer abstraitement le peigne~$C$ en la courbe~$\P^1_k$, de façon \emph{ad hoc}; l'étude infinitésimale des schémas de Hilbert
permettra ensuite de déformer le morphisme~$f$.

Soit~$T$ une courbe connexe lisse sur~$k$, munie d'un point rationnel $t \in T(k)$ (par exemple on peut prendre $T=\P^1_k$ et $t=0$).
Soit~$\Crond$ la surface lisse sur~$k$ obtenue en éclatant le point fermé $t \times M$ dans $T \times \P^1_k$.  Notons $\pi \colon \Crond \rightarrow T$
la composée de l'éclatement et de la première projection.  La courbe $\Crond_t=\pi^{-1}(t)$ est alors un peigne sur~$k$ (dont le manche est le transformé strict
de $t \times \P^1_k$ et dont les dents sont les diviseurs exceptionnels).  Ce peigne est isomorphe à~$C$, de sorte
qu'en posant $\Xrond = T\times X$ et en notant $\Xrond_t = t \times X = X$, le morphisme $f \colon C \rightarrow X$ induit un
morphisme $f_t \colon \Crond_t \rightarrow \Xrond_t$.

Notons $\Hom_T(\Crond,\Xrond)$ le schéma de Hilbert des $T$\nobreakdash-morphismes de~$\Crond$ vers~$\Xrond$.  C'est un $T$\nobreakdash-schéma.
Le morphisme $f_t$ en définit un $k$\nobreakdash-point situé au-dessus de $t \in T(k)$.  Tout
$k$\nobreakdash-point de $\Hom_T(\Crond,\Xrond)$ dont l'image dans~$T$ est \emph{différente} de~$t$ fournit un morphisme $\P^1_k \rightarrow X$;
mais l'image de ce morphisme peut ne pas contenir~$x$.
Pour y remédier, considérons l'image
réciproque, par l'éclatement $\Crond \rightarrow T \times \P^1_k$, de la courbe $T \times \infty$.  C'est une section de $\pi \colon \Crond \rightarrow T$.
Notons-la $\infty_T \subset \Crond$ et posons $x_T = T \times x \subset \Xrond$.  Le morphisme $f_t$ définit alors un $k$\nobreakdash-point $[f_t]$
du schéma de Hilbert $\Hom_T(\Crond,\Xrond;\infty_T \mapsto x_T)$ des $T$\nobreakdash-morphismes de~$\Crond$ vers~$\Xrond$ qui envoient~$\infty_T$ sur~$x_T$.

\bigskip
\begin{proposition}
\label{kollcproplisse}
Le $T$\nobreakdash-schéma $\Hom_T(\Crond,\Xrond;\infty_T \mapsto x_T)$ est lisse au point $[f_t]$.
\end{proposition}

\bigskip
Admettons un instant cette proposition.  Elle implique, grâce à la fertilité de~$k$,
que le schéma $\Hom_T(\Crond,\Xrond;\infty_T \mapsto x_T)$ admet un $k$\nobreakdash-point dont l'image dans~$T$ appartient à $T \setminus \{t\}$;
autrement dit, qu'il existe un morphisme $\P^1_k \rightarrow X$ envoyant $\infty \in \P^1(k)$ sur $x \in X(k)$.
Pour trouver un morphisme très libre l'argument sera un peu plus élaboré.

\bigskip
\begin{demo}[ de la proposition~\ref{kollcproplisse}]%
L'étude infinitésimale des schémas de Hilbert fournit le critère de lissité suivant:

\bigskip
\begin{proposition}[ (critère de lissité pour les schémas Hom)]%
\label{criterelissite}
Soit~$T$ un schéma noethérien irréductible muni d'un point $t \in T$.
Soit~$\Crond$ un $T$\nobreakdash-schéma projectif et plat de dimension relative~$1$, à fibres géométriquement réduites.
Soit~$\Xrond$ un $T$\nobreakdash-schéma lisse et quasi-projectif.
Soit $M \subset \Crond$ un sous-schéma fermé fini et plat sur~$T$
tel que~$\Crond$ soit lisse sur~$T$ aux points de~$M$.
Soient enfin $g \colon M \rightarrow \Xrond$ un $T$\nobreakdash-morphisme
et $f_t \colon \Crond_t \rightarrow \Xrond_t$ un $t$\nobreakdash-morphisme
dont la restriction à~$M_t$ coïncide avec la restriction de~$g$.
Pour que le $T$\nobreakdash-schéma
$\Hom_T(\Crond,\Xrond;g)$
paramétrant les $T$\nobreakdash-morphismes de~$\Crond$ vers~$\Xrond$ dont la restriction
à~$M$ coïncide avec~$g$ soit lisse au point~$[f_t]$, il suffit que l'on ait
$$H^1(\Crond_t, (f_t^\star T_{\Xrond_t}) \otimes \Irond_{M_t})=0$$ où $T_{\Xrond_t}$ désigne le faisceau tangent de $\Xrond_t$ sur~$t$
et où $\Irond_{M_t} \subset \Orond_{\Crond_t}$ désigne le faisceau (localement libre)
d'idéaux du sous-schéma fermé $M_t \subset \Crond_t$.
\end{proposition}

\bigskip
\begin{demo}
Voir \cite[II.1.7]{kollrational} et \cite[2.11]{debarrelivre}.
\end{demo}

\bigskip
Pour établir la proposition~\ref{kollcproplisse} il suffit
donc de vérifier que l'on a
$$H^1(C, f^\star T_X \otimes_{\Orond_C} \Irond_\infty)=0\text{,}$$
où $\Irond_\infty \subset \Orond_C$ désigne le faisceau d'idéaux
défini par le point $(\{0\} \times \infty) \in C_0 \subset C$.
Le $\Orond_C$\nobreakdash-module $f^\star T_X \otimes_{\Orond_C} \Irond_\infty$ est localement libre.
Sa restriction à $C_0 \subset C$ est isomorphe à $\Orond(-1)^{\dim(X)}$.  Sa restriction
aux dents de~$C$ est ample puisqu'elle coïncide avec la restriction de~$f^\star T_X$.  Le lemme général suivant permet donc de conclure.
\end{demo}

\bigskip
\begin{lemme}
\label{lemmepeigne}
Soit~$C$ un peigne sur un corps algébriquement clos.  Supposons que le manche de~$C$ soit une courbe rationnelle.
Soit~$E$ un $\Orond_C$\nobreakdash-module localement libre dont la restriction à chaque
dent est isomorphe à une somme directe $\bigoplus_{i=1}^r \Orond(a_i)$ pour des $a_i \geq 0$ et dont la restriction au manche est
isomorphe à une somme directe $\bigoplus_{i=1}^r \Orond(a_i)$ pour des $a_i \geq -1$.
Alors $H^1(C,E)=0$.
\end{lemme}

\bigskip
\begin{demo}
Notons~$M$ le manche de~$C$ et~$D$ la réunion des dents de~$C$, de sorte que $C = M \cup D$.
Notons $i_D \colon D \hookrightarrow C$ et $i_M \colon M \hookrightarrow C$ les inclusions.
On a alors une suite exacte de $\Orond_C$\nobreakdash-modules
$$
\xymatrix{
0 \ar[r] & \Irond_M \ar[r] & \Orond_C \ar[r] & i_{M\star}\Orond_M \ar[r] & 0
}
$$
d'où, après tensorisation par~$E$, une suite exacte
$$
\xymatrix{
0 \ar[r] & i_{D\star}i_D^\star(E \otimes \Irond_M) \ar[r] & E \ar[r] & i_{M\star}i^\star_ME \ar[r] & 0
}
$$
(compte tenu de l'égalité $\Irond_M = i_{D\star}i_D^\star\Irond_M$).  On en déduit une suite exacte
$$
\xymatrix{
H^1(D, E|_D \otimes \Irond_{M \cap D}) \ar[r] & H^1(C,E) \ar[r] & H^1(M,E|_M)\text{.}
}
$$
L'hypothèse sur $E|_M$ entraîne tout de suite que $H^1(M,E|_M)=0$.  D'autre part, le $\Orond_D$\nobreakdash-module $\Irond_{M \cap D}$ est localement libre
et sa restriction à chaque dent est isomorphe à~$\Orond(-1)$.  Vu l'hypothèse sur $E|_D$, il en résulte que $H^1(D,E|_D \otimes \Irond_{M \cap D})=0$.
D'où le lemme.
\end{demo}

\bigskip
Terminons maintenant la preuve du théorème~\ref{corpslocauxkollarth}.
D'après la proposition~\ref{kollcproplisse}, il existe une courbe connexe
$B \subset \Hom_T(\Crond,\Xrond;\infty_T \mapsto x_T)$
lisse sur~$k$, passant par le point $[f_t]$ et dominant~$T$.
La propriété universelle du schéma de Hilbert fournit
un $B$\nobreakdash-morphisme $$\Crond \times_T B \longrightarrow \Xrond \times_T B$$ envoyant $\infty_T \times_T B$ sur $x_T \times_T B$,
ou de façon équivalente un $k$\nobreakdash-morphisme
$$
\phi \colon \Crond \times_T B \longrightarrow X
$$
tel que $\phi(\infty_T \times_T B)=\{x\}$.

Comme $[f_t]\in B(k)$, on a $B(k)\neq\emptyset$. L'ensemble~$B(k)$
est donc dense dans~$B$ puisque~$k$ est fertile.  Notons $B^0 = B \times_T (T \setminus \{t\})$.
Soit $b \in B^0(k)$.
La restriction de~$\phi$ à $\Crond \times_T b$ est un morphisme $\phi_b \colon \P^1_k \rightarrow X$ envoyant~$\infty$ sur~$x$.  Que $\phi_b^\star T_X$ soit ample équivaut à ce que
$H^1(\P^1_k,\phi_b^\star T_X(-2))=0$ (en effet $H^1(\P^1_k,\Orond(d))=0$ si et seulement si $d\geq -1$).  Par conséquent, compte tenu de la densité de $B^0(k)$ dans~$B^0$ et
du théorème de semi-continuité supérieure,
une condition suffisante pour qu'il existe un $b \in B^0(k)$ tel que $\phi_b^\star T_X$ soit ample est l'existence d'un
diviseur effectif $D \subset \Crond$, de degré~$\geq 2$ sur les fibres de $\pi \colon \Crond \rightarrow T$, tel que,
si l'on pose $\Frond=\phi^\star T_X \otimes \Irond_{D \times_T B}$,
on ait $H^1(\Crond \times_T [f_t],\Frond|_{\Crond \times_T [f_t]})=0$.  Ici $\Irond_{D \times_T B}$ désigne le faisceau (localement libre) d'idéaux de $\Orond_{\Crond \times_T B}$ défini par
$D \times_T B$.

Prenons pour $D \subset \Crond$ la réunion de~$\infty_T$ et du transformé strict de $T \times M$ par l'éclatement $\Crond \rightarrow T \times \P^1_k$.
Le diviseur $D_t \otimes_k \kbar$ sur le peigne $\Crond_t \otimes_k \kbar = (\Crond \times_T [f_t]) \otimes_k \kbar$ est de degré~$1$ sur le manche et de degré~$1$ sur chaque dent.
Par conséquent la restriction de $\Frond|_{\Crond \times_T [f_t]}$ au manche est isomorphe à $\Orond(-1)^{\dim(X)}$ et sa restriction
à chaque dent est isomorphe à une somme directe $\bigoplus_{i=1}^r \Orond(a_i)$ pour des $a_i \geq -1$.
Le lemme~\ref{lemmepeigne} permet donc de conclure.

\subsection{Preuve du corollaire~\ref{kolloccorreq}}
\label{seckolloccorpreuve}

\medskip
Le corps~$k$ est maintenant un corps $p$\nobreakdash-adique ou est le corps des réels.

Pour établir le corollaire il suffit de montrer que les classes de
$R$\nobreakdash-équivalence sont ouvertes.  En effet elles seront alors
automatiquement fermées, puisqu'elles forment une partition de $X(k)$.  De plus,
comme le corps~$k$ est localement compact et que la variété~$X$ est propre, l'espace topologique $X(k)$ est compact,
de sorte que la finitude de $X(k)/R$ s'ensuivra également.

Soit donc $x \in X(k)$; nous allons prouver que la classe de $R$\nobreakdash-équivalence de~$x$ dans~$X(k)$ est un voisinage de~$x$.
On peut supposer que $\dim(X)>0$.
Fixons un morphisme $f \colon \P^1_k \rightarrow X$ tel que $f(0)=x$ et que $f^\star T_X$ soit ample.

Que $f^\star T_X$ soit très ample signifie que la courbe $f(\P^1_k)$ se déforme dans toute direction donnée, en fixant le point~$x$.  On peut fabriquer ainsi
de nombreux points $R$\nobreakdash-équivalents à~$x$, arbitrairement proches de~$x$ et distincts de~$x$.  Mais cela ne signifie pas pour autant
que tout point suffisamment proche de~$x$ est $R$\nobreakdash-équivalent à~$x$. Une astuce permettant d'aboutir à cette conclusion consiste à déformer la courbe $f(\P^1_k)$ en fixant
un \emph{autre} point que~$x$, disons~$y$.  Les déformations du point~$x$ couvriront alors un voisinage de~$x$ et seront toutes $R$\nobreakdash-équivalentes à~$y$,
donc à~$x$ puisque $y$ est lui-même $R$\nobreakdash-équivalent à~$x$.  Plus formellement:

Notons $H = \Hom^0(\P^1_k,X\;\!\!;\infty\mapsto f(\infty))$ l'espace de modules des morphismes très libres de~$\P^1_k$ vers~$X$ qui envoient~$\infty$ sur $f(\infty)$.
C'est un $k$\nobreakdash-schéma localement de type fini (cf.~\cite[II.3.5.4]{kollrational}).
Le morphisme d'évaluation $H \times \P^1_k \rightarrow X$ 
est lisse sur $H \times \A^1_k$ (cf.~\cite[II.3.5.3]{kollrational} ou~\cite[4.8]{debarrelivre}).
L'application induite $H(k) \times \A^1(k) \rightarrow X(k)$ est donc ouverte (théorème d'inversion locale).
Ainsi son image est une partie ouverte de~$X(k)$ contenant~$x$ (puisque $f(0)=x$)
et incluse dans la classe de $R$\nobreakdash-équivalence de~$x$ (car incluse dans la classe de $R$\nobreakdash-équivalence de~$f(\infty)$, qui est $R$\nobreakdash-équivalent à~$f(0)=x$).
La classe de $R$\nobreakdash-équivalence de~$x$ est donc bien un voisinage de~$x$.

\subsection{$R$-équivalence et $R$-équivalence directe sur les corps fertiles}
\label{reqreqdirecte}

\medskip
Soient~$k$ un corps et~$X$ une variété projective, lisse et séparablement rationnellement connexe, sur~$k$.  Nous venons de voir (théorème~\ref{corpslocauxkollarth})
que si~$k$ est fertile, par tout point rationnel de~$X$ passe une courbe rationnelle très libre.  Si~$k$ est algébriquement clos, on sait que l'on peut même trouver une courbe rationnelle
aussi libre que l'on veut passant par tout ensemble fini donné de points rationnels de~$X$.  Cette propriété ne saurait être satisfaite si le corps~$k$ est seulement supposé fertile:
par exemple, si $k=\R$, deux points rationnels de~$X$ appartenant à des composantes connexes distinctes de $X(\R)$ ne peuvent
certainement pas être joints par une courbe rationnelle définie sur~$\R$.  Il se pose néanmoins la question de savoir s'il existe toujours une courbe rationnelle aussi libre que l'on veut
 passant par un ensemble fini donné
de points rationnels tous $R$\nobreakdash-équivalents entre eux.  Kollár~\cite{kollsp} a montré que tel est bien le cas si~$k$ est fertile, à l'aide d'une technique
empruntée à Graber, Harris et Starr~\cite{graberharrisstarr}:

\bigskip
\begin{theoreme}[ (Kollár~\cite{kollsp})]%
\label{kospec}
Soient~$k$ un corps fertile et~$X$ une variété projective et lisse sur~$k$, séparablement rationnellement connexe,  de dimension~\mbox{$\geq 3$}.
Soient $n \geq 1$ et $\uplet{x_1}{x_n} \in X(k)$.  Supposons les~$x_i$ deux à deux $R$\nobreakdash-équivalents.
Alors il existe une immersion fermée $f \colon \P^1_k \rightarrow X$ telle que $(f^\star T_X)(-n)$ soit ample et que $\{\uplet{x_1}{x_n}\} \subset f(\P^1_k)$.
\end{theoreme}

\bigskip
L'hypothèse $\dim(X) \geq 3$ n'est pas essentielle. Elle permet d'exiger que~$f$ soit une immersion fermée, ce qui simplifie l'énoncé du théorème.
De toute manière, si $\dim(X)=2$, on peut appliquer le théorème à la variété $X \times \P^1_k$ et obtenir une conclusion presque aussi bonne.
En particulier, le théorème~\ref{kospec} généralise\footnote{La preuve du théorème~\ref{kospec} ne fait pas appel au théorème~\ref{corpslocauxkollarth}.  En fait,
le théorème~\ref{corpslocauxkollarth} découle même du théorème~\ref{thkodireq} (appliqué à $x=y$).}
le théorème~\ref{corpslocauxkollarth}.

Moret-Bailly~\cite[Corollaire~2]{mbcomplement} avait auparavant démontré un résultat très proche du théorème~\ref{kospec}
dans la situation où~$X$ est le quotient de $\mathrm{GL}_n$ par un sous-groupe fini étale.

Le théorème~\ref{kospec} est une conséquence facile du théorème ci-dessous, qui en constitue donc véritablement le c\oe{}ur:

\bigskip
\begin{theoreme}
\label{thkodireq}
Soient~$k$ un corps fertile et~$X$ une variété projective, lisse et séparablement rationnellement connexe, sur~$k$.
Soient $x,y \in X(k)$.  Si~$x$ et~$y$ sont directement $R$\nobreakdash-équivalents (cf.~§\ref{subsecreq}), il existe un morphisme $f \colon \P^1_k \rightarrow X$
tel que $f^\star T_X$ soit ample et  que $f(0)=x$, $f(\infty)=y$.
\end{theoreme}

\bigskip
\begin{esquissedemo}[ du théorème~\ref{kospec} à partir du théorème~\ref{thkodireq}]%
Tout d'abord on déduit du théorème~\ref{thkodireq} que la $R$\nobreakdash-équivalence directe est une relation transitive sur~$X(k)$ (de sorte qu'il n'y a pas de différence
entre les notions de $R$\nobreakdash-équivalence et de $R$\nobreakdash-équivalence directe sur~$X(k)$).

Pour cela, fixons $x,y,z \in X(k)$ tels que~$x$ soit directement $R$\nobreakdash-équivalent à~$y$ et que~$y$ soit directement $R$\nobreakdash-équivalent à~$z$.
D'après le théorème~\ref{thkodireq}, il existe des morphismes très libres $f \colon \P^1_k \rightarrow X$ et $g \colon \P^1_k \rightarrow X$
tels que $f(0)=x$, $f(\infty)=y$, $g(0)=y$ et $g(\infty)=z$.
Comme~$f$ et~$g$ sont très libres, la technique de déformation employée au~§\ref{chap2preuveproprementdite} montre que
le morphisme du peigne à une dent vers~$X$ obtenu en recollant~$f$ et~$g$
se déforme en un morphisme $h \colon \P^1_k \rightarrow X$
tel que $h(0)=x$ et $h(\infty)=z$.  Ainsi~$x$ et~$z$ sont bien directement $R$\nobreakdash-équivalents.

Prouvons maintenant le théorème~\ref{kospec}.  Soient $\uplet{x_1}{x_n} \in X(k)$ des points deux à deux $R$\nobreakdash-équivalents.
Fixons un point $y \in X(k)$ appartenant à la classe de $R$\nobreakdash-équivalence des~$x_i$ (par exemple $y=x_1$).
Comme~$y$ est $R$\nobreakdash-équivalent aux~$x_i$, il leur est même directement $R$\nobreakdash-équivalent;
d'après le théorème~\ref{thkodireq}, il existe donc des morphismes très libres $f_i \colon \P^1_k \rightarrow X$ tels que $f_i(0)=y$ et $f_i(\infty)=x_i$.
On peut assembler les~$f_i$ en un morphisme d'un peigne vers~$X$, le manche du peigne étant entièrement envoyé sur le point~$y$.
La technique du~§\ref{chap2preuveproprementdite} montre que ce peigne et ce morphisme se déforment en un morphisme $f \colon \P^1_k \rightarrow X$
tel que $(f^\star T_X)(-n)$ soit ample et qu'il existe $\uplet{y_1}{y_n} \in \P^1(k)$ vérifiant $f(y_i)=x_i$ pour tout~$i$.
Qu'en déformant~$f$ on puisse supposer que~$f$ est une immersion fermée résulte enfin de~\cite[II.3.14]{kollrational}.
\end{esquissedemo}

\bigskip
Voici l'idée de la démonstration du théorème~\ref{thkodireq}.  Tout d'abord, quitte à remplacer~$X$ par $X \times \P^1_k$ et $x$, $y$ par $x\times 0$ et $y\times \infty$,
on peut supposer que $x\neq y$.  Comme~$x$ et~$y$ sont par hypothèse directement $R$\nobreakdash-équivalents, il existe $g \colon \P^1_k \rightarrow X$
tel que $g(0)=x$ et $g(\infty)=y$.   Rappelons que si~$k$ est algébriquement clos, Kollár, Miyaoka et Mori ont montré qu'il suffit, pour trouver un~$g$ tel que~$g^\star T_X$ soit de plus ample,
de choisir un nombre suffisamment élevé de courbes rationnelles libres sur~$X$ rencontrant $g(\P^1_k)$ en des points deux à deux distincts; ces courbes rationnelles libres forment les dents
d'un peigne dont~$g$ est le manche; un \emph{sous-peigne} du peigne obtenu se déforme alors en un morphisme  $g' \colon \P^1_k \rightarrow X$  très libre vérifiant $g'(0)=x$, $g'(\infty)=y$.
(Voir~\cite[§4.2]{bonavero} pour tout cela.) 
Lorsque~$k$ n'est pas algébriquement clos, cette approche ne fonctionne plus.  En effet, il y a en général une
obstruction au problème de déformation considéré, \emph{i.e.} l'espace de modules associé n'est pas lisse\footnote{Si~$C$ est un peigne de manche~$C_0$ et si $f \colon C \rightarrow X$ est un
morphisme, pour que $H^1(C,f^\star T_X)=0$ il faut que $H^1(C_0,f_0^\star T_X)=0$, où~$f_0$ désigne la restriction de~$f$ à~$C_0$. Or cette condition n'a rien à voir avec les dents du peigne.
Si elle n'est pas satisfaite, on ne peut donc pas la forcer en ajoutant des dents à~$C$.};
or la fertilité de~$k$ ne fournit pas de nouveaux points rationnels sur une variété dont on sait seulement qu'elle admet un point rationnel singulier (ainsi par exemple la surface affine normale réelle d'équation $x^2+y^2+z^2=0$ admet
un unique point à coordonnées réelles).
Ce qui permet à l'argument d'aboutir lorsque~$k$ est algébriquement clos, ce n'est pas un calcul de lissité mais une comparaison des dimensions des espaces de modules
correspondant aux déformations de tous les sous-peignes possibles (cf.~\cite[§5]{araujokollar}).

Graber, Harris et Starr~\cite{graberharrisstarr} se sont rendu compte que la situation est meilleure si l'on cherche à déformer non pas un morphisme $f \colon C \rightarrow X$
d'un peigne vers~$X$ mais le sous-schéma fermé $f(C) \subset X$ lui-même (du moins si~$f$ est une immersion fermée, ce qui, au cours de la preuve du théorème~\ref{thkodireq}
est une condition que l'on peut toujours supposer satisfaite).
Pour ce nouveau problème de déformation,
l'obstruction s'évanouit si l'on ajoute au peigne un nombre élevé de dents libres suffisamment générales\footnote{L'espace de modules en jeu ici est le schéma de Hilbert $\Hilb(X,x,y)$ paramétrant
les sous-schémas fermés de~$X$ qui contiennent~$x$ et~$y$. L'étude infinitésimale de $\Hilb(X,x,y)$ est gouvernée non  par le faisceau $f^\star T_X$ mais par le faisceau normal de~$C$ dans~$X$.
La restriction de celui-ci à~$C_0$ \emph{est} sensible à l'ajout de dents au peigne (contrairement à $f^\star T_X$).}.

Cette technique est discutée dans~\cite{starrclay}.  L'adaptation à un corps non algébriquement clos ne présente qu'une difficulté: comme le peigne construit doit être
défini sur~$k$, lorsqu'on ajoute une dent au peigne sur~$\kbar$ il faut aussi ajouter ses conjuguées sous l'action de $\Gal(\kbar/k)$.
Or il se pourrait que la dent rencontre le manche au même point que l'une de ses conjuguées, auquel cas on n'obtient plus un peigne.  Un lemme de Kollár~\cite[Lemma~14]{kollsp}
montre que ce phénomène peut être évité.

\section{Variétés rationnellement connexes sur les corps finis et sur les corps pseudo-algébriquement clos}
\label{seckollsz}

Ce chapitre reprend les résultats établis en~2003 par Kollár et
Szabó~\cite{kollsz} au sujet des variétés séparablement
rationnellement connexes sur les corps finis.
La présentation que nous en donnons diffère légèrement de celle de Kollár et Szabó en
ceci que nous établissons d'abord un résultat sur les corps pseudo-algébriquement clos, duquel nous déduisons
le théorème principal de~\cite{kollsz} sur les corps finis, alors que Kollár et Szabó
travaillent sur un corps parfait arbitraire pour obtenir un énoncé dont ils tirent ensuite des conséquences
d'une part sur les corps parfaits pseudo-algébriquement clos et d'autre
part sur les corps finis.
Aucune perte de généralité ne résulte de cette façon de procéder et la présentation s'en trouve quelque peu allégée (notamment parce que les corps pseudo-algébriquement clos ont l'avantage d'être
fertiles; mais aussi parce que nous n'avons pas besoin de suivre dans toute la preuve la dépendance des constructions intermédiaires
par rapport aux entiers $\dim(X)$, $\deg(X)$ et $\deg(S)$ qui apparaissent dans~\cite[Th.~2]{kollsz}; et enfin parce que le recours
aux espaces de courbes stables de Kontsevich devient superflu).  Bien sûr, la preuve reste la même.

Dans la démonstration de Kollár et Szabó intervient un théorème «~de type Lefschetz~» pour les courbes rationnelles sur les variétés rationnellement connexes.
Au paragraphe~\ref{subsecapplicationgalinv} nous discutons une application de ce théorème au problème de Galois inverse.

\subsection{Énoncé du théorème principal; conséquences}

\medskip
Rappelons qu'un corps~$k$ est dit \emph{pseudo-algébriquement clos} si toute variété géométriquement intègre sur~$k$
admet un point rationnel (voir~\textsection\ref{subsubsecautrescorps}).

Soit~$X$ une variété projective, lisse et séparablement rationnellement connexe, sur~$k$.
La question~\ref{questionexisteample} admet une réponse positive si~$k$ est pseudo-algébriquement clos puisque tout corps pseudo-algébriquement clos est fertile.
Comme on l'a discuté au~§\ref{reqreqdirecte}, il n'est pas vrai en revanche que si~$k$ est fertile et si $S \subset X$ est un ensemble fini de points rationnels,
il existe sur~$X$ une courbe rationnelle passant par~$S$.
Le théorème principal de ce chapitre affirme qu'une telle courbe existe bien lorsque~$k$ est pseudo-algébriquement clos.  Les points de~$S$ n'ont même pas besoin d'être rationnels:

\bigskip
\begin{theoreme}[ (Kollár--Szabó~\cite{kollsz})]%
\label{thkspac}
Soit~$k$ un corps pseudo-algébriquement clos.
Soit~$X$ une variété projective et lisse sur~$k$, séparablement rationnellement connexe, de dimension~$\geq 2$.
Soit $S \subset X$ une sous-variété de dimension~$0$ lisse sur~$k$.
Alors il existe un morphisme $f \colon \P^1_k \rightarrow X$
tel que $(f^*T_X)(-\deg(S))$ soit ample et que l'inclusion $S \subset X$ se factorise par~$f$.
De plus, si $\dim(X)\geq 3$, on peut demander que~$f$ soit une immersion fermée.
\end{theoreme}

\bigskip
Si~$f$ est une immersion fermée, la condition «~l'inclusion $S \subset X$ se factorise par~$f$~» est bien sûr équivalente à la condition «~$S \subset f(\P^1_k)$~».
Du théorème~\ref{thkspac} résulte\footnote{Le corollaire~\ref{correqtrivpac} est immédiat si~$k$ est parfait puisque le support de tout $0$\nobreakdash-cycle sur~$X$ est alors lisse sur~$k$
(condition imposée à~$S$ dans le théorème~\ref{thkspac}).
Par ailleurs, l'assertion sur la $R$\nobreakdash-équivalence est
claire que~$k$ soit parfait ou non.  Pour établir la nullité de~$\CHz(X)$ en général, fixons un point rationnel $x \in X(k)$.
Il suffit de montrer que tout point fermé $z \in X$ est rationnellement équivalent à~$\deg(z)x$.  Pour cela, il suffit que sur la variété $X \otimes_k k(z)$,
les points rationnels~$x$ et~$z$ soient rationnellement équivalents.  Mais cela découle du théorème~\ref{thkspac} appliqué à la variété $X \otimes_k k(z)$ sur le corps~$k(z)$
(qui est encore pseudo-algébriquement clos d'après~\cite[11.2.5]{friedjarden}).}:

\bigskip
\begin{corollaire}
\label{correqtrivpac}
Soient~$k$ un corps pseudo-algébriquement clos
et~$X$ une variété projective, lisse et séparablement rationnellement connexe, sur~$k$.
Alors $\Card(X(k)/R)=1$ et $\CHz(X)=0$.
\end{corollaire}

\bigskip
Le corollaire~\ref{correqtrivpac} est à mettre en parallèle avec les questions~\ref{questionax} et~\ref{questionctqac}.

Du théorème~\ref{thkspac} nous déduirons au~§\ref{subsecdescorpspacauxfinis}
le même énoncé avec pour~$k$ un corps fini de cardinal assez grand devant~$\dim(X)$, $\deg(X)$ et~$\deg(S)$, où~$\deg(X)$ désigne le degré de~$X$ dans un plongement
projectif $X \subset \P^N_k$ donné:

\bigskip
\begin{corollaire}
\label{corkolszfini}
Il existe une application $\Phi \colon \N^3 \rightarrow \N$ telle que l'assertion suivante soit vérifiée.
Soit~$k$ un corps fini.
Soit $X \subset \P^N_k$ une variété projective et lisse sur~$k$, séparablement rationnellement connexe et de dimension~$\geq 2$.
Soit $S \subset X$ une sous-variété de dimension~$0$ lisse sur~$k$.
Si le cardinal de~$k$ est supérieur à $\Phi(\dim(X),\deg(X),\deg(S))$, il existe un morphisme $f \colon \P^1_k \rightarrow X$
tel que $(f^*T_X)(-\deg(S))$ soit ample et que l'inclusion $S \subset X$ se factorise par~$f$.
De plus, si $\dim(X)\geq 3$, on peut demander que~$f$ soit une immersion fermée.
\end{corollaire}

\bigskip
Sans hypothèse sur le cardinal du corps fini~$k$, on ne peut pas espérer trouver un morphisme $f \colon \P^1_k \rightarrow X$ qui soit une immersion fermée, même
si $\dim(X) \geq 3$.  En effet cela impliquerait que~$X$ admet au moins $\Card(k)+1$ points rationnels, or Swinnerton-Dyer~\cite{swinnertondyeruniversal} a donné un exemple
d'une surface cubique lisse~$V$ sur un corps fini~$k$ telle que $\Card(V(k))=1$.  En posant $X=V \times V$ on a alors $\Card(X(k))=1$ et $\dim(X) \geq 3$.
Il se pourrait néanmoins que privé de sa dernière phrase, l'énoncé du corollaire~\ref{corkolszfini}
soit vrai sans hypothèse sur le cardinal de~$k$, à condition toutefois de supposer que~$S$ se plonge dans~$\P^1_k$.
C'est une question ouverte même dans le cas où $S=\emptyset$.

Le corollaire~\ref{corkolszfini} implique tout de suite:

\bigskip
\begin{corollaire}
\label{corkolszreqfini}
Il existe une application $\Phi \colon \N^2 \rightarrow \N$ telle que l'assertion suivante soit vérifiée.
Soient~$k$ un corps fini et $X \subset \P^N_k$ une variété projective, lisse et séparablement rationnellement connexe, sur~$k$.
Si $\Card(k) \geq \Phi(\dim(X),\deg(X))$ alors $\Card(X(k)/R)=1$.
\end{corollaire}

\bigskip
La question de savoir si l'énoncé du corollaire~\ref{corkolszreqfini} reste vrai sans l'hypothèse sur le cardinal de~$k$ est ouverte.

On peut aussi déduire du corollaire~\ref{corkolszfini} l'égalité $\CHz(X)=0$ pour~$X$ projective, lisse et séparablement rationnellement connexe sur un corps fini,
et cette fois quel que soit le cardinal du corps de base
(cf.~\cite{kollsz}).
Mais ce résultat était déjà connu.  En effet, Kato et Saito~\cite{katosaito} prouvent que pour \emph{toute} variété projective, lisse et géométriquement
connexe sur un corps fini~$k$, le groupe $\CHz(X)$ est isomorphe à un quotient du groupe fondamental abélianisé $\pi_1^{\mathrm{ab}}(X \otimes_k \kbar)$; or Kollár a montré
(à l'aide du théorème de Graber, Harris, Starr et de Jong) que le groupe fondamental de toute variété projective, lisse et séparablement rationnellement connexe sur un corps algébriquement
clos est trivial (cf.~\cite[3.6]{debarrebourbaki}).

Le corollaire~\ref{corkolszfini} a également des conséquences sur les corps $p$\nobreakdash-adiques.
Si~$X$ est une variété projective et lisse sur un corps $p$\nobreakdash-adique~$k$, on dit que~$X$ a \emph{bonne réduction séparablement rationnellement connexe}
si, notant~$\Orond$ l'anneau des entiers de~$k$, il existe un $\Orond$\nobreakdash-schéma projectif et lisse
dont la fibre spéciale géométrique est une variété séparablement rationnellement connexe et dont la fibre générique est isomorphe à~$X$.
Cela entraîne que la variété~$X$ elle-même est rationnellement connexe (cf.~\cite[IV.3.11]{kollrational}).

Nous prouverons les corollaires~\ref{corkolszreqfiniloc} et~\ref{corkolszreqfinilocch} au~§\ref{kolfinisubseccortant}.

\bigskip
\begin{corollaire}
\label{corkolszreqfiniloc}
Il existe une application $\Phi \colon \N^2 \rightarrow \N$ telle que l'assertion suivante soit vérifiée.
Soient~$k$ un corps $p$\nobreakdash-adique et~$X$ une variété projective et lisse sur~$k$, telle qu'il existe
un modèle projectif $\Xrond \subset \P^N_{\Orond}$ de~$X$ dont la fibre spéciale soit une variété lisse et séparablement rationnellement connexe.
Si le cardinal du corps résiduel de~$k$ est $\geq \Phi(\dim(X),\deg(X))$, alors $\Card(X(k)/R)=1$.
\end{corollaire}

\bigskip
Ici encore, c'est une question ouverte de savoir si le corollaire~\ref{corkolszreqfiniloc} reste vrai sans l'hypothèse sur le cardinal du corps résiduel de~$k$.

Sur les corps $p$\nobreakdash-adiques, le corollaire sur le groupe de Chow des $0$\nobreakdash-cycles est cette fois entièrement nouveau.  L'énoncé ci-dessous
n'était connu que dans le cas des surfaces rationnelles (grâce à Bloch, Colliot-Thélène, Sansuc; cf.~\cite{ctkdeux}).

\bigskip
\begin{corollaire}
\label{corkolszreqfinilocch}
Soient~$k$ un corps $p$\nobreakdash-adique et~$X$ une variété projective et lisse sur~$k$, ayant bonne réduction séparablement rationnellement connexe.
Alors $\CHz(X)=0$.
\end{corollaire}

\bigskip
Des corollaires~\ref{corkolszreqfiniloc} et~\ref{corkolszreqfinilocch} résulte immédiatement:

\bigskip
\begin{corollaire}
\label{corkolszcdn}
Soient~$k$ un corps de nombres et~$X$ une variété projective, lisse et rationnellement connexe, sur~$k$.
Pour presque toute place~$v$ de~$k$, on a $\Card(X(k_v)/R)=1$ et $\CHz(X \otimes_k k_v)=0$.
\end{corollaire}

\bigskip
L'énoncé du corollaire~\ref{corkolszcdn} avait été conjecturé par Colliot-Thélène (voir notamment~\cite[3.1~(b)]{ctbordeaux}).

\subsection{Des corps pseudo-algébriquement clos aux corps finis}
\label{subsecdescorpspacauxfinis}

\medskip
Dans ce paragraphe nous déduisons le corollaire~\ref{corkolszfini} du théorème~\ref{thkspac}.

L'outil de départ est l'estimation donnée indépendamment par Lang et Weil~\cite{langweil} et par Nisnevi{\v{c}}~\cite{nisnevic}
du nombre de points rationnels d'une sous-variété fermée géométriquement irréductible~$X$
de $\P^n_k$, pour~$k$ fini, en fonction du cardinal de~$k$, du degré de~$X$ dans~$\P^n_k$ et de la dimension de~$X$.  Celle-ci a pour conséquence:

\bigskip
\begin{theoreme}[ (Lang--Weil--Nisnevi{\v{c}})]%
\label{thlangweil}
Il existe une application $\Phi \colon \N^3 \rightarrow \N$ telle que l'assertion suivante soit vérifiée.
Soient~$k$ un corps fini et~$X$ une variété quasi-projective et géométriquement irréductible sur~$k$.
Fixons une immersion $X \subset \P^n_k$.
Notons $\Xbarre$ l'adhérence de~$X$ dans $\P^n_k$ et posons $\partial X = \Xbarre \setminus X$.
Supposons enfin satisfaite l'inégalité
$\Card(k)\geq \Phi(n,\deg(\Xbarre),\deg(\partial X))$,
où $\deg(Z)$ désigne ici,
pour tout fermé $Z \subset \P^n_k$ (non nécessairement équidimensionnel),
la somme des degrés des composantes irréductibles de~$Z$.
Alors $X(k)\neq\emptyset$.
\end{theoreme}

\bigskip
Les corps pseudo-algébriquement clos auxquels on appliquera le théorème~\ref{thkspac} seront fournis par le lemme général suivant:

\bigskip
\begin{lemme}
\label{lemmepac}
Pour tout corps~$k$, il existe un corps pseudo-algébriquement clos~$K$ contenant~$k$ et dans lequel~$k$ est algébriquement fermé.
\end{lemme}

\bigskip
\begin{esquissedemo}
L'idée remonte à Merkurjev et consiste \emph{grosso modo} à construire~$K$ inductivement en partant de~$k$ et en remplaçant~$K$ par~$K(X)$ chaque fois que~$X$ est une variété
géométriquement intègre sur~$K$.
Cette construction est formulée de façon rigoureuse par exemple dans~\cite{ctcras} (preuve du théorème~2.1).
\end{esquissedemo}

\bigskip
Le lemme~\ref{lemmepac}, le théorème de Lang--Weil--Nisnevi{\v{c}} et le théorème~\ref{thkspac} entraînent de façon formelle que si~$k$ est un corps fini et si l'on se donne
$X$ et~$S$ vérifiant les hypothèses du corollaire~\ref{corkolszfini}, il existe un entier~$N$ tel que pour toute extension finie~$k'/k$ de degré~\mbox{$\geq N$}, il existe un morphisme
$f \colon \P^1_{k'} \rightarrow X \otimes_k k'$ tel que l'inclusion $S \otimes_k k' \subset X \otimes_k k'$ se factorise par~$f$ et que $(f^\star T_X)(-\deg(S))$ soit ample.
Malheureusement, cette information ne permet de démontrer aucun des corollaires~\ref{corkolszreqfini}, \ref{corkolszreqfiniloc}, \ref{corkolszreqfinilocch} et~\ref{corkolszcdn}.
Au minimum il faudrait pour cela que l'entier~$N$ ne dépende pas du choix de~$S$ (mais seulement de son degré).  Pour le corollaire~\ref{corkolszcdn} il est crucial que~$N$ ne dépende
même pas de la caractéristique de~$k$.  Le lemme suivant nous permettra de résoudre ces difficultés et d'établir le corollaire~\ref{corkolszfini} dans toute sa généralité.

\bigskip
\begin{lemme}
\label{monlemme}
Soient~$S$ un schéma de type fini sur~$\Z$ et $f \colon X \rightarrow S$ un
morphisme de schémas localement de type fini.
Supposons que pour tout corps pseudo-algébriquement clos~$L$, l'application $X(L) \rightarrow S(L)$ induite par~$f$ soit surjective.
Alors il existe un entier~$N$ tel que pour tout point $s \in S$ dont le corps résiduel est fini et de cardinal~$\geq N$, la fibre $f^{-1}(s)$ possède un point
rationnel.
\end{lemme}

\bigskip
\begin{demo}
Quitte à remplacer~$S$ tour à tour par chacune de ses composantes irréductibles (munies de la structure de sous-schéma fermé réduit), on peut supposer~$S$ intègre.
Notons~$\kappa(S)$ le corps des fonctions de~$S$.  D'après le lemme~\ref{lemmepac}, il existe un corps pseudo-algébriquement clos~$L$ contenant~$\kappa(S)$ et
tel que~$\kappa(S)$ soit algébriquement fermé dans~$L$.  Comme l'application $X(L) \rightarrow S(L)$ est surjective, la fibre générique de~$f$ admet un $L$\nobreakdash-point;
et comme~$f$ est localement de type fini, on en déduit qu'il existe une sous-extension de type fini $L_0/\kappa(S)$
de $L/\kappa(S)$ telle que la fibre générique de~$f$
admette un $L_0$\nobreakdash-point.  Le corps~$L_0$ est alors le corps des fonctions d'une variété intègre, affine et géométriquement irréductible sur~$\kappa(S)$, puisque~$\kappa(S)$ est
séparablement fermé dans~$L_0$.  D'où l'existence d'un ouvert affine dense $S^0 \subset S$ et d'un $S^0$\nobreakdash-schéma affine~$W^0$, de type fini et à fibres géométriquement irréductibles,
muni d'un $S$\nobreakdash-morphisme $W^0 \rightarrow X$.  Fixons une $S^0$\nobreakdash-immersion fermée $W^0 \subset \A^n_{S^0}$ pour un $n \geq 1$.
Quitte à rétrécir~$S^0$, on peut supposer 
que pour $s \in S^0$, les degrés dans~$\P^n_s$ de $\overline{W^0_s}$ et de $\partial(W^0_s)$ sont indépendants de~$s$
(avec les notations du théorème~\ref{thlangweil} relatives à $W^0_s \subset \A^n_s \subset \P^n_s$).
Comme les fibres de $W^0 \rightarrow S^0$ sont géométriquement irréductibles,
le théorème de Lang--Weil--Nisnevi{\v{c}} rappelé ci-dessus
entraîne alors l'existence d'un entier~$N$ tel que pour tout point $s \in S^0$ dont le corps résiduel est fini et de cardinal~$\geq N$,
la fibre de $W^0 \rightarrow S^0$ en~$s$ possède un point rationnel.  Compte tenu de l'existence d'un $S$\nobreakdash-morphisme $W^0 \rightarrow X$, il en résulte
que pour tout point $s \in S^0$ dont le corps résiduel est fini et de cardinal~$\geq N$, la fibre $f^{-1}(s)$ admet un point rationnel.
Pour conclure il ne reste plus qu'à établir la conclusion du lemme pour la restriction de~$f$ au-dessus de $S \setminus S^0$, ce qui est justiciable d'un raisonnement par récurrence
puisque $\dim(S \setminus S^0) < \dim(S)$.
\end{demo}

\bigskip
\begin{remarque}
La considération d'ultraproduits de corps finis permet de donner une seconde démonstration du lemme~\ref{monlemme} \emph{via} la théorie des modèles.

Lorsque~$f$ est un morphisme de type fini, cette approche ne présente aucune difficulté puisque si~$k$ désigne un corps,
la surjectivité de l'application $X(k) \to S(k)$ induite par~$f$ s'exprime alors comme une formule logique du premier ordre dans le corps~$k$.

L'argument suivant, issu d'une discussion avec Antoine Chambert-Loir pour laquelle je le remercie, permet de procéder sans supposer~$f$ de type fini.
Soit~$I$ l'ensemble des ouverts quasi-compacts de~$X$.  Pour $U \in I$, notons $I_U$ l'ensemble des $V \in I$ contenant~$U$.
Soit enfin~$\Brond$ l'ensemble des parties de $I \times \N$ de la forme
 $I_U \times C$ pour $U \in I$ et $C \subset \N$ cofini.
Toute intersection finie d'éléments de~$\Brond$ est non vide; il existe donc un ultrafiltre~$\Frond$ sur $I \times \N$ contenant~$\Brond$
(cf.~\cite[Corollary~7.5.3]{friedjarden}).
Supposons la conclusion du lemme~\ref{monlemme} en défaut.  Pour chaque entier~$N$ choisissons alors un corps fini $L_N$ de cardinal $\geq N$
et un point $s_N \in S(L_N)$ n'admettant pas d'antécédent dans $X(L_N)$.
Pour $U \in I$, notons $L_{U,N}=L_N$ et $s_{U,N}=s_N$.
Soit~$L$ le corps ultraproduit des $L_{U,N}$ selon~$\Frond$
puis $s \in S(L)$ l'ultraproduit des $s_{U,N}$.  Grâce au théorème de Lang--Weil--Nisnevi{\v{c}}
et compte tenu que
$I \times C = I_\emptyset \times C \in \Brond \subset \Frond$ pour toute partie $C \subset \N$ cofinie,
le corps~$L$ est pseudo-algébriquement clos
(cf.~\cite[Proposition~7.7.1]{friedjarden}).
D'après l'hypothèse du
lemme~\ref{monlemme}, le point~$s$ se relève donc dans~$X(L)$. En particulier il existe un $V \in I$
tel que~$s$ se relève dans~$V(L)$.
Rappelons d'autre part que quels que soient~$U$ et~$N$, le point $s_{U,N}$ ne se relève pas dans $U(L_{U,N})$; le point $s_{U,N}$ ne se relève donc dans $V(L_{U,N})$ pour aucun $(U,N) \in I_V \times \N$.
Comme $I_V \times \N \in \Brond \subset \Frond$, il s'ensuit
que~$s$ ne se relève pas dans $V(L)$ (cf.~\cite[Proposition~7.7.1]{friedjarden}), ce qui est absurde.
\end{remarque}

\bigskip
Démontrons maintenant le corollaire~\ref{corkolszfini}.
Fixons des entiers~$\delta$, $d$ et~$r$ (ce seront respectivement la dimension de~$X$, le degré de~$X$ et le degré de~$S$) et posons $n=\delta+d-1$.
Quels que soient l'entier~$N$ et le corps~$k$, toute sous-variété géométriquement irréductible
de~$\P^N_k$ de degré~$d$ et de dimension~$\delta$ est incluse dans un sous-espace projectif de~$\P^N_k$ de dimension au plus~$n$
(cf.~\cite[Corollary~18.12]{harrisfirstcourse}). C'est ce qui nous permettra ci-dessous de nous borner à considérer les sous-variétés de~$\P^n_k$ de degré~$d$,
plutôt que les sous-variétés de~$\P^N_k$ de degré~$d$ et de dimension~$\delta$, avec~$N$ variable.

Soit $H_1 \subset \Hilb(\P^n_\Z/\Z)$ l'ensemble des points $h \in \Hilb(\P^n_\Z/\Z)$ tels que la sous-variété fermée
correspondante de~$\P^n_h$ soit lisse (sur le corps résiduel de~$h$), de
degré~$d$ et séparablement rationnellement connexe.  L'ensemble~$H_1$ est un
ouvert de~$\Hilb(\P^n_\Z/\Z)$ (cf.~\cite[IV.3.11]{kollrational}).  Notons $\Xrond_1
\rightarrow H_1$ la famille universelle, de sorte que $\Xrond_1$ est un sous-schéma fermé de
$\P^n_{H_1}$ lisse sur~$H_1$.
Soit $H_2 \subset \Hilb(\Xrond_1/H_1)$ l'ouvert constitué des points~$h$ tels
que la sous-variété fermée correspondante de $\Xrond_1 \times_{H_1} h$ soit
lisse, de dimension~$0$ et de degré~$r$ (sur le corps résiduel de~$h$).
Posons $\Xrond_2 = \Xrond_1 \times_{H_1} H_2$ et
notons $\Srond_2 \rightarrow H_2$ la famille universelle, de sorte que~$\Srond_2$ est un sous-schéma
fermé de $\Xrond_2$, fini étale de degré~$r$ sur~$H_2$.
Notons $\Hom_{H_2}^0(\P^1_{H_2},\Xrond_2) \subset \Hom_{H_2}(\P^1_{H_2},\Xrond_2)$ le schéma de Hilbert des
$H_2$\nobreakdash-morphismes $f \colon \P^1_{H_2} \rightarrow \Xrond_2$ tels
que $(f^\star T_{\Xrond_2/H_2})(-r)$ soit ample relativement à~$H_2$.
Soit enfin~$M$ l'image réciproque, par le morphisme de composition
$$
\Hom_{H_2}(\Srond_2,\P^1_{H_2}) \times_{H_2} \Hom_{H_2}^0(\P^1_{H_2},\Xrond_2) \longrightarrow \Hom_{H_2}(\Srond_2,\Xrond_2)\rlap{\text{,}}
$$
de la section de $\Hom_{H_2}(\Srond_2,\Xrond_2) \rightarrow H_2$ correspondant à l'inclusion $\Srond_2 \subset \Xrond_2$.

Pour résumer le paragraphe ci-dessus, disons que l'on a défini, au-dessus de $\Spec(\Z)$, le schéma de modules~$H_2$ paramétrant
les données qui apparaissent dans le corollaire~\ref{corkolszfini}: une sous-variété projective lisse et séparablement rationnellement connexe
de degré~$d$ de l'espace projectif de dimension~$n$ (fibre de $\Xrond_2 \rightarrow H_2$), munie d'une sous-variété lisse de dimension~$0$ et de degré~$r$
(fibre de $\Srond_2 \rightarrow H_2$).  De plus on a défini un $H_2$\nobreakdash-schéma~$M$ dont la fibre au-dessus de chaque point de~$H_2$ est l'espace de modules
paramétrant les morphismes~$f$ que l'on cherche à construire.
Tous ces schémas sont localement de type fini sur~$\Z$.

Le théorème~\ref{thkspac} affirme exactement que pour tout corps pseudo-algébriquement clos~$L$, l'application induite $M(L) \rightarrow H_2(L)$ est surjective.
Si~$H_2$ est de type fini sur~$\Z$, le lemme~\ref{monlemme} appliqué à $M \rightarrow H_2$ entraîne donc immédiatement
la validité du corollaire~\ref{corkolszfini}: on définit $\Phi(\delta,d,r)$ comme étant
l'entier~$N$ donné par le lemme~\ref{monlemme}.
(Pour prouver la toute dernière assertion du corollaire~\ref{corkolszfini}
on raisonne de la même manière, à ceci près que l'on remplace~$M$ par l'ouvert de~$M$
obtenu en ajoutant la condition que le morphisme~$f$ est une immersion fermée; il s'agit bien d'un ouvert de~$M$ d'après~\cite[I.6.2]{kollrational}.)

Il suffit donc, pour conclure, de vérifier que le schéma~$H_2$ est de type fini sur~$\Z$; autrement dit, qu'il est quasi-compact,
ou encore, qu'il ne possède qu'un nombre fini de composantes connexes.
C'est seulement ici que le fait d'avoir fixé~$\delta$, $d$ et~$r$ va jouer un rôle.
Le morphisme $H_2 \rightarrow H_1$ est quasi-compact car~$H_2$ est inclus dans $\Hilb_r(\Xrond_1/H_1)$, qui est un $H_1$\nobreakdash-schéma projectif.
(Ici $\Hilb_r(\Xrond_1/H_1)$ désigne le lieu dans $\Hilb(\Xrond_1/H_1)$ des sous-variétés de polynôme de Hilbert égal à~$r$.)
Il~reste donc seulement à vérifier que le schéma~$H_1$ est quasi-compact.  Or cela résulte tout de suite de la théorie des formes de Chow (cf.~\cite[I.6.6.1]{kollrational}).
Ce qui compte ici est que~$H_1$ paramètre des sous-variétés fermées réduites de degré borné dans un espace projectif de dimension bornée.

\subsection{Preuve du théorème principal}

\medskip
Nous supposons dorénavant les hypothèses du théorème~\ref{thkspac} satisfaites.  En particulier, le corps~$k$ est  pseudo-algébriquement clos.

\subsubsection{Esquisse de l'argument}

Le principe de la preuve est le suivant.  Supposons, pour simplifier, que $\dim(X)\geq 3$ et que tous les points de~$S$ soient $k$\nobreakdash-rationnels.  Notons-les $S=\{\uplet{x_1}{x_n}\}$.
Nous voulons construire
une courbe rationnelle sur~$X$ passant par tous les~$x_i$.
On commence par appliquer à chaque~$x_i$ le théorème principal du chapitre précédent
(théorème~\ref{corpslocauxkollarth}), ce qui est licite puisque le corps~$k$ est fertile, étant pseudo-algébriquement clos.  On obtient, pour chaque~$i$,
une courbe rationnelle très libre $C_i \subset X$ passant par~$x_i$.  On aimerait voir les~$C_i$ comme les dents d'un peigne.  Pour cela il faut construire un manche,
c'est-à-dire trouver une courbe rationnelle sur~$X$ qui rencontre tous les~$C_i$.
Nous ne construirons une telle courbe qu'après avoir déformé les~$C_i$.
Plus précisément, comme les~$C_i$ sont des courbes rationnelles très libres, elles s'insèrent dans des familles qui d'une part sont
paramétrées par des variétés lisses et qui d'autre part balaient un ouvert dense de~$X$.
En d'autres termes, pour chaque~$i$ il existe une variété lisse et connexe~$H_i$ munie d'un point rationnel $h_i \in H_i(k)$
et un morphisme dominant $\phi_i \colon \P^1_{H_i} \rightarrow X$ dont la restriction à~$\P^1_{h_i}$ a pour image~$C_i$ et dont la restriction à~$\P^1_h$, pour tout $h \in H_i$,
a pour image une courbe rationnelle très libre contenant~$x_i$.
«~Fabriquer un manche quitte à déformer les~$C_i$~» revient alors
 à trouver une courbe rationnelle~$C_0$ sur~$X$ dont l'image réciproque par~$\phi_i$ contienne, pour tout~$i$, un point rationnel\footnote{En toute rigueur il faut aussi s'assurer que deux dents ne rencontrent jamais le manche au même point; sans cela on n'obtient pas un peigne.}.
Comme~$k$ est pseudo-algébriquement clos, il suffit pour cela que la variété $\phi_i^{-1}(C_0) \subset \P^1_{H_i}$ soit géométriquement irréductible sur~$k$.
Or $\P^1_{H_i}$ est lui-même géométriquement irréductible sur~$k$ puisque $H_i$ est lisse et connexe et contient un point rationnel (à savoir~$h_i$).
On est ainsi ramené au problème suivant: étant donnés une variété géométriquement irréductible~$Y$
et un morphisme dominant $Y \rightarrow X$ (ici~$\P^1_{H_i}$ et $\phi_i \colon \P^1_{H_i} \rightarrow X$), montrer que la restriction de~$Y$ au-dessus de «~presque toute~» courbe rationnelle sur~$X$ est encore
une variété géométriquement irréductible.  Ce problème, qui rappelle les théorèmes de Lefschetz sur les sections hyperplanes des variétés projectives
et qui est essentiellement de nature géométrique, a été résolu
par Kollár~\cite{kollarlefschetz}.  Le résultat de Kollár entraîne l'existence d'une famille
de courbes rationnelles très libres sur~$X$, paramétrée par une variété géométriquement irréductible~$H_0$, telle que pour tout~$i$, l'image réciproque par~$\phi_i$ de toute courbe
apparaissant dans cette famille soit une variété géométriquement irréductible.  Comme~$k$ est pseudo-algébriquement clos, on a $H_0(k)\neq\emptyset$.  D'où finalement l'existence
d'un peigne sur~$X$ dont le manche~$C_0$ est une courbe rationnelle très libre et dont les dents $\uplet{C_1}{C_n}$ sont des courbes rationnelles très libres telles que $x_i \in C_i$ pour tout~$i$.
La technique déjà employée dans la preuve du théorème~\ref{corpslocauxkollarth} montre alors que ce peigne se déforme en
une courbe rationnelle lisse $C \subset X$ contenant tous les~$x_i$ et vérifiant la condition d'amplitude qui apparaît dans le théorème~\ref{corpslocauxkollarth}.
(L'hypothèse $\dim(X)\geq 3$ sert à assurer que l'on peut choisir~$C$ lisse.  Si~$C$ n'est pas lisse, la condition «~$C$ contient les~$x_i$~» n'est pas tout à fait
celle qui nous intéresse.)

\subsubsection{La preuve proprement dite}
\label{subsecproprementdite}

Soit $s \in S$.  Comme le corps~$k(s)$ est une extension finie séparable
de~$k$, c'est un corps pseudo-algébriquement clos (par restriction des
scalaires à la Weil)
et en particulier c'est un corps fertile.  Par conséquent, d'après le théorème~\ref{corpslocauxkollarth}, le $k(s)$\nobreakdash-schéma
$\Hom^0(\P^1_{k(s)},X \otimes_k k(s)\;\!; 0 \mapsto s)$
qui paramètre les $k(s)$\nobreakdash-morphismes très libres de~$\P^1_{k(s)}$ vers $X \otimes_k k(s)$ envoyant~$0$ sur~$s$
admet un point rationnel.  Choisissons-en un et notons~$H_s$ un voisinage ouvert de ce point.  Quitte à rétrécir~$H_s$ on peut supposer que~$H_s$ est une variété connexe.
La variété~$H_s$ est alors géométriquement irréductible sur~$k(s)$ puisqu'elle est lisse et connexe et qu'elle admet un point rationnel sur~$k(s)$.

Notons $\phi_s \colon \P^1_{k(s)} \times H_s \rightarrow X \otimes_k k(s)$ le $k(s)$\nobreakdash-morphisme universel paramétré par~$H_s$.
C'est un morphisme dominant puisque sa restriction à $(\P^1_{k(s)} \setminus \{0\}) \times H_s$ est lisse (cf.~\cite[II.3.5.3]{kollrational}).
Le théorème de type Lefschetz que nous allons maintenant appliquer est le suivant:

\bigskip
\newcommand{\refkollef}{\cite[Theorem~6, Corollary~7 et~§28]{kollarlefschetz}}
\begin{theoreme}[ (Kollár~\refkollef)]%
\label{thkollarlefschetz}
Soient~$k$ un corps et~$X$ une variété projective, lisse et séparablement rationnellement connexe, sur~$k$.  Notons~$\kbar$ une clôture algébrique de~$k$
et $\Hom^0(\P^1_k,X)$ l'ouvert du schéma de Hilbert $\Hom(\P^1_k,X)$ constitué des morphismes très libres.
Alors il existe un ouvert géométriquement intègre $U \subset \Hom^0(\P^1_k,X)$ tel que
pour toute variété irréductible~$Y$ sur~$\kbar$ munie d'un morphisme dominant $\phi \colon Y \rightarrow X \otimes_k \kbar$, il existe un
ouvert dense $V \subset U$ vérifiant la propriété suivante: pour tout $v \in V(\kbar)$,
le produit fibré de $\phi \colon Y \rightarrow X \otimes_k \kbar$ avec $v \colon \P^1_\kbar \rightarrow X \otimes_k \kbar$ est irréductible.
\end{theoreme}

\bigskip
Soit $U \subset \Hom^0(\P^1_k,X)$ comme dans le théorème~\ref{thkollarlefschetz}.  Quitte à rétrécir~$U$, on peut supposer
que pour tout~$s$, la conclusion du théorème~\ref{thkollarlefschetz} 
pour $\phi=\phi_s \otimes_{k(s)} \kbar$
est satisfaite avec $V=U$.  Quitte à rétrécir~$U$ de nouveau, on peut supposer de plus
que l'image du morphisme universel $U \times \P^1_k \rightarrow X$ est disjointe de~$S$ (en effet le morphisme $U \times \P^1_k \rightarrow X$ est lisse~\cite[II.3.5.4]{kollrational}
et~$S$ est de codimension~$\geq 2$ dans~$X$).
Choisissons alors $u \in U(k)$ (un tel choix est possible car~$k$ est pseudo-algébriquement clos et~$U$ est géométriquement intègre). Pour $s \in S$, considérons le carré
\begin{equation}
\begin{aligned}
\label{lecarre}
\xymatrix@C=10ex{
D_s \ar[d]^{\psi_s} \ar[r] & \P^1_{k(s)} \times H_s \ar[d]^{\phi_s} \\
\P^1_{k(s)} \ar[r]^(.45){u \otimes_k k(s)} & X \otimes_k k(s)
}
\end{aligned}
\end{equation}
où~$D_s$ est défini comme le produit fibré de~$\phi_s$ et de $u \otimes_k k(s)$ (et où~$\psi_s$ est la seconde projection).

Comme~$\phi_s$ est lisse au-dessus du complémentaire de $S \otimes_k k(s)$ (cf.~\cite[II.3.5.3]{kollrational}),
le morphisme~$\psi_s$ est lisse.  En particulier il est dominant et~$D_s$ est une variété géométriquement intègre sur~$k(s)$
(elle est géométriquement irréductible par construction).
Comme~$k$ est pseudo-algébriquement clos, il en résulte que l'on peut choisir pour chaque $s \in S$
un $k(s)$\nobreakdash-point de~$D_s$ de telle façon que les images dans~$\P^1_k$ de tous ces points (par $D_s \xrightarrow{\psi_s} \P^1_{k(s)} \rightarrow \P^1_k$)
soient deux à deux distinctes et que le corps résiduel de chacune de ces images soit~$k(s)$.  Pour forcer cette dernière condition, il suffit de remarquer que trouver un
$k(s)$\nobreakdash-point
de~$D_s$ dont l'image par $D_s \xrightarrow{\psi_s} \P^1_{k(s)} \rightarrow \P^1_k$ ait pour corps résiduel~$k(s)$ équivaut, en notant $R_{k(s)/k}$ la restriction des scalaires à la Weil de~$k(s)$ à~$k$,
à trouver un point $k$\nobreakdash-rationnel de $R_{k(s)/k}D_s$ dont l'image dans $R_{k(s)/k}\P^1_{k(s)}$ appartienne à l'ouvert complémentaire de la réunion
des sous-variétés strictes $R_{k'/k}\P^1_{k'}$ où~$k'$ parcourt les sous-extensions de $k(s)/k$ distinctes de~$k(s)$.

Autrement dit, il existe $(h_s)_{s \in S} \in \prod_{s \in S}
H_s(k(s))$,  $(q_s)_{s \in S} \in \prod_{s \in S} \P^1(k(s))$ et
 $(m_s)_{s \in S} \in \prod_{s \in S} \P^1(k(s))$ 
 satisfaisant pour tout~$s \in S$
l'égalité $u(m_s)=\phi_s(q_s,h_s)$ dans $X(k(s))$
et
telles que les images dans~$\P^1_k$ des~$m_s$ soient deux à deux distinctes et de corps résiduel~$k(s)$;
 de plus, on a $q_s \neq 0$ pour tout~$s$ puisque $\phi_s(0,h_s)=s$ alors que l'image
de $U\times \P^1_k \rightarrow X$ est disjointe de~$S$.  Toutes ces données définissent un peigne (muni d'un morphisme vers~$X$), que nous allons maintenant
déformer en procédant exactement comme au~§\ref{chap2preuveproprementdite}.

Soit~$T$ une courbe lisse et connexe sur~$k$ munie d'un point rationnel $t \in T(k)$.
Les points~$m_s$ définissent une immersion fermée $S \subset \P^1_k$.
Soient~$\Crond$ la surface lisse obtenue en éclatant $t \times S$ dans $T \times \P^1_k$
et $\pi \colon \Crond \rightarrow T$ la composée de l'éclatement et de la première projection.
Pour $s \in S$, notons $R_s \subset \Crond$ le transformé strict de $T \times s$ et $E_s \subset \Crond$
la fibre de l'éclatement $\Crond \rightarrow T \times \P^1_k$ au-dessus de $t \times s$.
Posons enfin $R = \bigcup_{s \in S} R_s$ et notons $C_0 \subset \Crond$ le transformé
strict de $t \times \P^1_k$.
Sur~$E_s$ il y a deux points $k(s)$\nobreakdash-rationnels distingués: les points $E_s \cap C_0$ et $E_s \cap R_s$.
Fixons un $k(s)$\nobreakdash-isomorphisme $E_s \isoto \P^1_{k(s)}$ envoyant le premier sur~$q_s$ et le second sur~$0$.
Les morphismes $\phi_s(-,h_s) \colon \P^1_{k(s)} \rightarrow X \otimes_k k(s)$ et $u \colon \P^1_k \rightarrow X$
se recollent alors en un morphisme $f_t \colon \Crond_t \rightarrow \Xrond_t$ (où $\Xrond=T \times X$ et $\Xrond_t = t \times X = X$)
tel que $f_t(R|_{\Crond_t})=t \times S$.

\bigskip
\begin{proposition}
\label{kollpacproplisse}
Le $T$\nobreakdash-schéma $\Hom_T(\Crond,\Xrond; R \mapsto T \times S)$ est lisse au point~$[f_t]$.
\end{proposition}

\bigskip
\begin{demo}
D'après la proposition~\ref{criterelissite} il suffit de vérifier que
$$H^1(\Crond_t, f_t^\star T_X \otimes_{\Orond_{\Crond_t}} \Irond)=0$$
où $\Irond \subset \Orond_{\Crond_t}$ désigne le faisceau d'idéaux défini par $R|_{\Crond_t} \subset \Crond_t$.
Or cela résulte du lemme~\ref{lemmepeigne}, compte tenu que~$\Crond_t$ est un peigne et que les restrictions de $f_t^\star T_X$
au manche et aux dents sont des faisceaux amples.
\end{demo}

\bigskip
D'après la proposition~\ref{kollpacproplisse}, il existe une courbe $B \subset \Hom_T(\Crond,\Xrond; R \mapsto T \times S)$ lisse sur~$k$, connexe, passant par~$[f_t]$
et dominant~$T$.  Notons $\phi \colon \Crond \times_T B \rightarrow X$ le morphisme donné par le propriété universelle
de $\Hom_T(\Crond,\Xrond; R \mapsto T \times S)$
et posons $B^0 = B \times_T (T \setminus \{t\})$.
Comme $[f_t] \in B(k)$ et que le corps~$k$ est fertile (étant pseudo-algébriquement clos), l'ensemble $B(k)$ est dense dans~$B$.
Pour $b \in B^0(k)$, la restriction de~$\phi$ à $\Crond \times_T b$ est un morphisme
$\phi_b \colon \P^1_k \rightarrow X$.
L'inclusion $S \subset X$ se factorise par~$\phi_b$ puisque
$\phi_b(R \times_T b)=S$.
Il reste donc seulement à vérifier que l'on peut choisir~$b$ de sorte que le faisceau
$(\phi_b^\star T_X)(-\deg(S))$ soit ample, c'est-à-dire que 
$H^1(\P^1_k,(\phi_b^\star T_X)(-\deg(S)-2))=0$.  Par semi-continuité de la cohomologie, il suffit pour cela
qu'il existe une section $D \subset \Crond$ de $\pi \colon \Crond \rightarrow T$
telle que, si l'on pose $\Frond=\phi^\star T_X \otimes \Irond_{R \times_T B} \otimes \Irond_{D \times_T B}^{\otimes 2}$
(où $\Irond_{R \times_T B}$ et $\Irond_{D \times_T B}$ désignent respectivement les faisceaux (localement libres) d'idéaux
de $\Orond_{\Crond \times_T B}$ définis par $R \times_T B$ et $D \times_T B$),
 on ait $H^1(\Crond \times_T [f_t],\Frond|_{\Crond \times_T [f_t]})=0$.
Fixons un point $z \in \P^1(k)$ tel que $z \not\in S$ et prenons
pour~$D$ l'image réciproque, par l'éclatement $\Crond \rightarrow T \times \P^1$, de $T \times z$.
La restriction de $\Frond|_{\Crond \times_T [f_t]}$ au manche
du peigne $\Crond \times_T [f_t] = \Crond_t$ est alors isomorphe à une somme directe $\bigoplus_{i=1}^r \Orond(a_i)$ pour des $a_i \geq -1$
et sa restriction à chaque dent est isomorphe à une somme directe  $\bigoplus_{i=1}^r \Orond(a_i)$ pour des $a_i \geq 0$;
le lemme~\ref{lemmepeigne} permet donc de conclure.

Supposons enfin $\dim(X) \geq 3$ et prouvons la dernière assertion du théorème~\ref{thkspac}.
Comme $(\phi_b^\star T_X)(-\deg(S))$ est ample et que $\dim(X) \geq 3$,
le $k$\nobreakdash-schéma $$\Hom(\P^1_k,X\;\!; S \mapsto S)$$ est lisse en~$[\phi_b]$
(cf.~\cite[2.9]{debarrelivre})
et tout voisinage ouvert de~$[\phi_b]$ contient un ouvert non vide
paramétrant des immersions fermées (cf.~\cite[II.3.14.3]{kollrational}).
Par conséquent, compte tenu que~$k$ est fertile, tout voisinage ouvert de~$[\phi_b]$ contient un point rationnel correspondant à une immersion fermée $f \colon \P^1_k \rightarrow X$.
Si l'on choisit un voisinage ouvert suffisamment petit, le faisceau $(f^\star T_X)(-\deg(S))$ sera nécessairement ample.

\subsection{Preuves des corollaires~\ref{corkolszreqfiniloc} et~\ref{corkolszreqfinilocch}}
\label{kolfinisubseccortant}

\medskip
\begin{demo}[ du corollaire~\ref{corkolszreqfiniloc}]
Notons~$F$ le corps résiduel de~$\Orond$ et $\Phi \colon \N^3 \rightarrow N$ l'application donnée par le corollaire~\ref{corkolszfini}.
Soient $x,y \in X(k)$.  Comme~$\Xrond$ est propre sur~$\Orond$, les points~$x$ et~$y$ s'étendent en des sections $\xtilde, \ytilde \in \Xrond(\Orond)$.
Quitte à remplacer~$X$ par $X \times \P^1_k$, $\Xrond$ par $\Xrond \times \P^1_\Orond$ et~$x$, $y$ par $(x,0)$, $(y,1)$,
on peut supposer que $\xtilde \cap \ytilde = \emptyset$.
Si l'inégalité $\Card(F) \geq \Phi(\dim(X),\deg(X),2)$ est satisfaite,
il existe alors un morphisme très libre $f \colon \P^1_F \rightarrow \Xrond \otimes_\Orond F$
tel que $f(0)=\xtilde \otimes_\Orond F$ et $f(1)=\ytilde \otimes_\Orond F$.
D'où un point rationnel~$[f]$ de la fibre spéciale du $\Orond$\nobreakdash-schéma $\Hom_\Orond(\P^1_\Orond,\Xrond;0\mapsto \xtilde,1\mapsto \ytilde)$.
Ce $\Orond$\nobreakdash-schéma est lisse en~$[f]$ puisque~$f$ est très libre
(cf.~proposition~\ref{criterelissite}).
Par conséquent, et compte tenu que~$\Orond$ est un anneau local complet, on peut relever~$[f]$ en une section de
$\Hom_\Orond(\P^1_\Orond,\Xrond;0\mapsto \xtilde,1\mapsto \ytilde)$ sur $\Spec(\Orond)$; autrement dit, $f$ est la restriction à~$\P^1_F$
d'un morphisme $\P^1_\Orond \rightarrow \Xrond$ qui envoie~$0$ sur~$\xtilde$ et~$1$ sur~$\ytilde$.  La restriction de ce morphisme
au-dessus du point générique de~$\Spec(\Orond)$ fournit une courbe rationnelle sur~$X$ reliant~$x$ à~$y$.
Ainsi a-t-on prouvé que $\Card(X(k)/R)\leq 1$ dès que $\Card(F) \geq \Phi(\dim(X),\deg(X),2)$.  D'autre part, il résulte du théorème de Lang--Weil--Nisnevi{\v{c}}
que $\Xrond(F) \neq \emptyset$, et donc $X(k)\neq\emptyset$, dès que $\Card(F)$ est assez grand
devant~$\dim(X)$ et~$\deg(X)$; d'où le corollaire.
\end{demo}

\bigskip
\begin{demo}[ du corollaire~\ref{corkolszreqfinilocch} (due à Colliot-Thélène)]%
Soit $\Xrond \subset \P^N_{\Orond}$ un modèle projectif de~$X$, lisse sur~$\Orond$ et de fibre spéciale séparablement rationnellement connexe.
Pour toute extension finie $k'/k$, la flèche $X \otimes_k k' \rightarrow X$ induit
(par image directe) une application norme $A_0(X \otimes_k k') \rightarrow A_0(X)$ dont la composée avec la flèche de
restriction $A_0(X) \rightarrow A_0(X \otimes_k k')$ est la multiplication par $[k':k]$ sur~$A_0(X)$.  Pour que $A_0(X)=0$, il suffit donc
qu'il existe des extensions $k'/k$ de degrés premiers entre eux telles que $A_0(X \otimes_k k')=0$.

Comme~$F$ est fini, il existe pour chaque $n \geq 1$ une extension $F_n/F$ de degré~$n$.  Celle-ci se relève en une extension non ramifiée $k_n/k$ de degré~$n$.
D'après le corollaire~\ref{corkolszreqfiniloc}, on a $\Card(X(k_n)/R)=1$ pour tout~$n$ assez grand.
En particulier on a $\Card(X(k_m)/R)=1$ et $\Card(X(k_n)/R)=1$
pour des entiers~$m$ et~$n$ premiers entre eux.  Quitte à remplacer~$k$ par~$k_m$ et par~$k_{n}$
on peut donc supposer que $\Card(X(k)/R)=1$, et même que $\Card(X(k')/R)=1$ pour toute extension finie~$k'/k$.
Fixons alors $x \in X(k)$.  Pour établir que $A_0(X)=0$ il suffit de montrer que tout point fermé de~$X$ est rationnellement équivalent à un multiple de~$x$.
Pour cela il suffit que pour tout point fermé $y \in X$, les points rationnels~$x$ et~$y$ de $X \otimes_k k(y)$ soient rationnellement équivalents sur
$X \otimes_k k(y)$ (en effet la norme de~$k(y)$ à~$k$ du $0$\nobreakdash-cycle $y-x$ sur $X \otimes_k k(y)$ est le $0$\nobreakdash-cycle $y-\deg(y)x$ sur~$X$);
or ils sont même $R$\nobreakdash-équivalents puisque $\Card(X(k(y))/R)=1$.
\end{demo}

\subsection{Application du théorème de type Lefschetz au problème de Galois inverse}
\label{subsecapplicationgalinv}

\medskip
Tout groupe fini est-il le groupe de Galois d'une extension galoisienne de~$\Q$~?  L'une des approches permettant de répondre partiellement à cette question ouverte
consiste à chercher des extensions finies galoisiennes du corps~$\Q(t)$ de groupe de Galois donné.  Par spécialisation,
une extension finie galoisienne de~$\Q(t)$ de groupe~$G$ fournit, pour une infinité de valeurs de~$t$, une extension galoisienne de~$\Q$ ayant~$G$ pour groupe de Galois
 (théorème d'irréductibilité de Hilbert,
cf.~\cite{serretopics}).

Que tout groupe fini soit groupe de Galois sur~$\C(t)$ est une conséquence classique du théorème d'existence de Riemann.
On commence par construire un revêtement galoisien topologique (ramifié) de $\P^1(\C)$
ayant le groupe de Galois désiré; le théorème d'existence de Riemann permet d'en déduire un revêtement de~$\P^1_\C$ et, par suite, une extension finie de $\C(t)$.
Par une méthode analogue, fondée sur la géométrie analytique rigide, Harbater~\cite{harbater} a démontré que tout groupe fini est le groupe de Galois
d'une extension galoisienne de~$\Qp(t)$, et plus généralement de~$k(t)$ pour tout corps~$k$ complet pour une valuation discrète.  On a même:

\bigskip
\begin{theoreme}[ (Harbater~\cite{harbater}, Pop~\cite{pop})]%
\label{thharbater}
Soit~$k$ un corps fertile.  Tout groupe fini est le groupe de Galois d'une extension galoisienne de~$k(t)$.
\end{theoreme}

\bigskip
Colliot-Thélène~\cite{ctann} a donné une démonstration nouvelle du théorème~\ref{thharbater} lorsque~$k$ est de caractéristique nulle, s'appuyant sur la
technique de déformation de courbes rationnelles employée dans le chapitre précédent.  Le résultat principal de~\cite{ctann} est plus précis:
il affirme que pour tout groupe fini~$G$ et tout torseur~$T$ sur~$k$ sous~$G$ (c'est-à-dire: tout $k$\nobreakdash-schéma fini étale~$T$ muni
d'une action simplement transitive de~$G$), il existe une extension galoisienne de~$k(t)$ de groupe de Galois~$G$, régulière sur~$k$, qui en $t=0$ se spécialise en le $k$\nobreakdash-schéma~$T$
(c'est-à-dire: il existe une courbe~$C$ lisse et géométriquement connexe sur~$k$ et un revêtement (ramifié) $C \rightarrow \P^1_k$ galoisien de groupe~$G$,
dont la fibre au-dessus de $0 \in \P^1(k)$ soit isomorphe à~$T$).
En prenant $T=\coprod_{g \in G} \Spec(k)$ on retrouve le théorème~\ref{thharbater} en caractéristique nulle; l'assertion avec~$T$ arbitraire était inconnue
dans cette généralité  (problème dit de Beckmann--Black).

Ce théorème a été généralisé dans plusieurs directions.  Moret-Bailly~\cite{mbconstruction} l'a étendu
aux corps fertiles de caractéristique quelconque, par une technique de recollement formel dans l'esprit de la preuve par Harbater du théorème~\ref{thharbater}.
Kollár~\cite{kollarlefschetz} a établi un énoncé de type Lefschetz (une légère variante du théorème~\ref{thkollarlefschetz} ci-dessus) ayant pour corollaire
une généralisation du théorème de Colliot-Thélène au cas où~$G$ est un groupe algébrique linéaire non nécessairement fini (toujours sous l'hypothèse de caractéristique nulle).

Nous nous contentons ici d'expliquer comment un théorème de type Lefschetz pour les courbes rationnelles sur les variétés rationnellement connexes peut avoir pour conséquence
le théorème de Harbater--Pop.

La variante du théorème~\ref{thkollarlefschetz} dont nous avons besoin est la suivante:

\bigskip
\newcommand{\refkollefd}{\cite[Theorem~3]{kollarlefschetz}}
\begin{theoreme}[ (Kollár~\refkollefd)]%
\label{thkollarlefschetzd}
Soient~$k$ un corps et~$X$ une variété projective, lisse et séparablement rationnellement connexe, sur~$k$.  Soit $x \in X(k)$.
Notons~$\kbar$ une clôture algébrique de~$k$
et $\Hom^0(\P^1_k,X\;\!\!;0 \mapsto x)$ l'ouvert du schéma de Hilbert $\Hom(\P^1_k,X\;\!;0\mapsto x)$ constitué des morphismes très libres.
Alors il existe un ouvert géométriquement intègre $U \subset \Hom^0(\P^1_k,X\;\!;0\mapsto x)$ et une compactification $U \subset U'$ tels que~$U'$
possède un point rationnel lisse et que
pour toute variété irréductible~$Y$ sur~$\kbar$ munie d'un morphisme dominant $\phi \colon Y \rightarrow X \otimes_k \kbar$, il existe un
ouvert dense $V \subset U$ vérifiant la propriété suivante: pour tout $v \in V(\kbar)$,
le produit fibré de $\phi \colon Y \rightarrow X \otimes_k \kbar$ avec $v \colon \P^1_\kbar \rightarrow X \otimes_k \kbar$ est irréductible.
\end{theoreme}

\bigskip
Autrement dit il s'agit du même énoncé que le théorème~\ref{thkollarlefschetz}, à ceci près qu'on ne considère que les courbes passant par un point rationnel donné de~$X$
et que l'ouvert~$U$ apparaissant dans la conclusion admet une compactification possédant un point rationnel lisse.  En fait, Kollár déduit le théorème~\ref{thkollarlefschetz}
du théorème~\ref{thkollarlefschetzd} (en étendant les scalaires de~$k$ à~$k(X)$ et en prenant pour~$x$ le point générique).

\bigskip
\begin{demo}[ du théorème~\ref{thharbater} en caractéristique nulle]%
Soient~$k$ un corps fertile de caractéristique~$0$ et~$G$ un groupe fini.  Choisissons un plongement $G \subset \mathrm{GL}_N$ pour un $N \geq 1$
et notons $Q=\mathrm{GL}_N/G$ le quotient.  Comme~$k$ est de caractéristique~$0$, il existe, d'après Hironaka, une compactification lisse $Q \subset X$,
une variété projective et lisse~$Y$ sur~$k$
et un morphisme $\phi \colon Y \rightarrow X$ dont la restriction au-dessus de~$Q$ soit isomorphe à la projection $\mathrm{GL}_N \rightarrow \mathrm{GL}_N/G$.
Notons $x \in Q(k) \subset X(k)$ l'image de $1 \in \mathrm{GL}_N(k)$ par~$\phi$.
La variété~$X$ est (séparablement) unirationnelle puisque~$\mathrm{GL}_N$ est un ouvert de l'espace affine.  En particulier elle est (séparablement) rationnellement connexe
et on peut donc lui appliquer le théorème~\ref{thkollarlefschetzd}.  D'où~$U$, $U'$ et~$V$ comme dans la conclusion de ce théorème.
Comme~$k$ est fertile et que~$U'$ admet un point rationnel lisse, l'ensemble~$V(k)$ est non vide.
Soit $v \in V(k)$.  Soit $C \rightarrow \P^1_k$ le produit fibré de $v \colon \P^1_k \rightarrow X$ par $\phi \colon Y \rightarrow X$.
Par construction, $C$ est une courbe géométriquement irréductible.  Le revêtement ramifié $C \rightarrow \P^1_k$ est galoisien de groupe~$G$
puisque $\phi^{-1}(Q) \rightarrow Q$ est galoisien de groupe~$G$ et que l'image de $v \colon \P^1_k \rightarrow X$ rencontre~$Q$ (à savoir, en~$x$).
Le groupe~$G$ est donc bien un groupe de Galois sur~$k(t)$.
\end{demo}

\bigskip
Il ressort de la preuve que~$G$ est même un groupe de Galois régulier sur~$k(t)$.
D'autre part, en variant le choix du point $x \in Q(k)$ dans l'argument, on contrôle à volonté la fibre de  $C \rightarrow \P^1_k$ au-dessus de $0 \in \P^1_k$, ce qui
résout le problème de Beckmann--Black pour~$k$ fertile de caractéristique~$0$.

\section{Existence de points rationnels sur les corps finis: le point de vue motivique}
\label{secesnault}

\subsection{Introduction}

\medskip
Le but de ce chapitre est d'établir le théorème suivant:

\bigskip
\begin{theoreme}[ (Esnault~\cite{esnaultinvent})]%
\label{thesn}
Sur un corps fini, toute variété propre, lisse et rationnellement connexe par chaînes admet un point rationnel.
\end{theoreme}

\bigskip
Ce théorème répond par l'affirmative à la question~\ref{laquestionc1} dans le cas des corps finis.
Il est à noter que la question~\ref{laquestionc1} portait seulement sur les variétés
séparablement rationnellement connexes.  L'hypothèse de connexité rationnelle
par chaînes est bien plus faible; par exemple le théorème~\ref{thesn}
s'applique aussi à des variétés de type général.  Cependant, la question~\ref{laquestionc1} étendue à toutes les variétés
rationnellement connexes par chaînes admet une réponse négative\footnote{Notons $X \subset \P^3_K$ l'hypersurface
de Fermat $w^n + tx^n + t^2y^n + t^3z^n = 0$ sur le corps $K=k((t))$, où~$k$ est algébriquement clos de caractéristique $p>2$.
D'après Shioda~\cite{shiodaexample}, cette variété est unirationnelle dès que~$n$ divise $p+1$
(cf.~aussi~\cite[Ex.~2.5.1]{debarrelivre}).
Par ailleurs, si~$n$ est premier à~$6$, on peut vérifier, par un calcul de valuations,
que le degré sur~$K$ de tout point fermé de~$X$ est divisible par~$n$.
Pour de nombreuses valeurs de~$n$ et de~$p$, la variété projective~$X$ est donc lisse et unirationnelle et n'admet pas de point $K$\nobreakdash-rationnel,
ni même de point rationnel sur la clôture parfaite de~$K$, bien que~$K$ (resp.~la clôture parfaite de~$K$) soit un corps~$(C_1)$.}.

Nous dirons qu'une variété~$X$ est \emph{de Fano} si elle est propre, lisse, géométriquement irréductible et si le faisceau inversible $\det(T_X)$ est ample.
Toute variété de Fano est rationnellement connexe par chaînes d'après un théorème dû à Campana et, indépendamment, à Kollár, Miyaoka et
Mori (cf.~\cite[V.2.13]{kollrational}). Le théorème~\ref{thesn} entraîne donc:

\bigskip
\begin{corollaire}
\label{corlangmanin}
Toute variété de Fano sur un corps fini admet un point rationnel.
\end{corollaire}

\bigskip
L'énoncé du corollaire~\ref{corlangmanin} avait été conjecturé par Lang et Manin~\cite[§2.6]{maninfano}.

C'est une question ouverte de savoir si toute variété de Fano
est séparablement rationnellement connexe\footnote{Kollár~\cite{kollnonrational} a néanmoins construit des exemples de variétés de Fano \emph{singulières} (c'est-à-dire
de variétés~$X$ propres et normales telles que $-K_X$ soit $\Q$\nobreakdash-Cartier et ample)
qui ne sont pas séparablement rationnellement connexes.}.
Il est donc crucial, pour obtenir le corollaire~\ref{corlangmanin}, que le théorème~\ref{thesn} ne se limite pas aux variétés séparablement rationnellement connexes.
En ce sens, le théorème~\ref{thesn} est à rapprocher du théorème de Chevalley--Warning (cf.~§\ref{subseccorpsci}), qui lui aussi s'applique à des variétés
de Fano sans supposer qu'elles sont séparablement rationnellement connexes.  (La question de savoir si toute hypersurface lisse $X \subset \P^n$
de degré $d \leq n$ est séparablement rationnellement connexe est ouverte (cf.~\cite[§4.5]{derenthalkollar}).  Dans le cas des hypersurfaces de Fermat
elle a été étudiée par Conduché~\cite{conduche}.)

Que le théorème~\ref{thesn} s'applique à toute variété rationnellement connexe par chaînes signifie en même temps qu'il n'est pas question pour le démontrer
d'utiliser des techniques de déformation de courbes rationnelles comme aux deux chapitres précédents.
La géométrie de~$X$ ne jouera qu'un rôle mineur dans la démonstration.  Le théorème sera établi en étudiant la cohomologie de~$X$ (son «~motif~»).

\subsection{Preuve du théorème~\ref{thesn}}
\label{preuveesn}

\medskip
La preuve du théorème~\ref{thesn} s'effectue en plusieurs étapes indépendantes, que nous détaillons dans les paragraphes~\ref{subsecformuleslefschetz} à~\ref{subsecchowconiveau}.

Dans tout le~§\ref{preuveesn} nous désignons par~$k$ un corps fini de cardinal~$q$ et de caractéristique~$p$, par~$\kbar$ une clôture algébrique de~$k$,
par $\F_{q^n}$ le sous-corps de~$\kbar$ de cardinal $q^n$,
par~$X$ une variété propre, lisse et géométriquement connexe sur~$k$
(non nécessairement rationnellement connexe par chaînes) et par~$F$ l'endomorphisme
de Frobenius $F \colon X \rightarrow X$, c'est-à-dire le morphisme de variétés sur~$k$ qui induit l'identité sur l'espace topologique sous-jacent
et qui agit sur~$\Orond_X$ par $f \mapsto f^q$.

\subsubsection{Formules de Lefschetz}
\label{subsecformuleslefschetz}

Pour étudier les points rationnels des variétés algébriques sur le corps fini~$k$, on définit, suivant Weil, la \emph{fonction zêta} de~$X$ par la formule
$$
\zeta_X(t)=\exp\left(\sum_{s=1}^\infty N_s \frac{t^s}{s}\right)
$$
où $N_s=\Card(X(\F_{q^s}))$.  C'est une série formelle en~$t$, à coefficients rationnels (et même entiers).

En formulant ses célèbres «~conjectures de Weil~» (aujourd'hui démontrées grâce, notamment, aux travaux de Dwork, Grothendieck, Artin et Deligne), Weil donnait à la fonction zêta une
interprétation cohomologique.  Son intuition était la suivante.
Comme le corps~$k$ est fini, l'ensemble $X(k)$ des points rationnels de~$X$ admet une définition purement géométrique:
c'est l'ensemble des points de~$X(\kbar)$ fixes par~$F$.  Ainsi, pour étudier les points rationnels on peut oublier le corps de base et raisonner seulement sur la variété
$\Xbarre=X \otimes_k \kbar$ munie de l'endomorphisme~$F \otimes_k \kbar$ (que l'on notera encore~$F$).
Les points fixes de~$F$ agissant sur~$\Xbarre$ sont les points d'intersection, dans $\Xbarre \times \Xbarre$,
de la diagonale $\Delta \subset \Xbarre \times \Xbarre$ avec le graphe $\Gamma \subset \Xbarre \times \Xbarre$ de~$F$.  Or ces deux sous-variétés de $\Xbarre \times \Xbarre$
se rencontrent transversalement en tout point d'intersection (car la différentielle de~$F$ est nulle);
ainsi~$N_s$ n'est autre que le nombre d'intersection de~$\Delta$ avec~$\Gamma$.  Dans la situation classique de la géométrie différentielle, le calcul du nombre d'intersection
de deux sous-variétés de codimensions complémentaires d'une variété donnée est bien connu pour être un problème de nature cohomologique.  Poursuivant ce raisonnement jusqu'à son terme,
Weil s'était rendu compte que l'existence
d'une bonne théorie cohomologique pour les variétés sur~$k$ (qui est un corps de caractéristique $p>0$) impliquerait formellement, pour la fonction zêta,
d'admirables propriétés, dont on trouvera la liste dans~\cite{weil}.

À l'heure actuelle on dispose de deux bonnes théories cohomologiques pour les variétés définies sur un corps~$k$ de caractéristique~$p$: la cohomologie (étale) $\ell$\nobreakdash-adique
(dépendant du choix d'un nombre premier $\ell\neq p$), à coefficients dans~$\Ql$, et la cohomologie rigide,
à coefficients dans le corps des fractions de l'anneau des vecteurs de Witt de~$k$.
Dans le reste de ce chapitre on pourrait travailler avec l'une quelconque de ces deux théories.  Le
choix n'a aucune incidence sur la preuve, ce qui reflète son caractère «~motivique~».
Comme l'article originel~\cite{esnaultinvent} et le rapport~\cite{chambsembour} sont rédigés avec la cohomologie rigide, nous emploierons ci-dessous la cohomologie $\ell$\nobreakdash-adique.
Fixons donc une fois pour toutes un nombre premier~$\ell$ différent de~$p$.

Pour $i \geq 0$, Grothendieck a défini le groupe de cohomologie $\ell$\nobreakdash-adique $H^i(\Xbarre,\Ql)$.  C'est un $\Ql$\nobreakdash-espace vectoriel de dimension finie,
nul si $i> 2 \dim(X)$.
Le morphisme de Frobenius $F \colon \Xbarre \rightarrow \Xbarre$ induit un
endomorphisme de $H^i(\Xbarre,\Ql)$; on le note encore~$F$.  Voici maintenant
comment la fonction zêta de~$X$ s'exprime en termes de la cohomologie
$\ell$-adique de~$X$.
Notons $\fp{P}i(t)=\det(1-tF \;|\; H^i(\Xbarre,\Ql)) \in \Ql[t]$
le polynôme caractéristique de l'endomorphisme de $H^i(\Xbarre,\Ql)$ induit par~$F$.
D'après Grothendieck et Artin, on a alors l'égalité
\begin{equation}
\label{eqlefschetz1}
\zeta_X(t) = \frac{\fp{P}1(t) \fp{P}3(t) \dots \fp{P}{2n-1}(t)}{\fp{P}0(t)\fp{P}2(t)\dots \fp{P}{2n}(t)}
\end{equation}
dans le corps $\Ql((t))$, où $n=\dim(X)$.  En particulier~$\zeta_X(t)$ est une fraction rationnelle
et l'égalité ci-dessus vaut dans $\Ql(t)$.
(En fait, d'après Deligne, les polynômes $\fp{P}i(t)$ appartiennent à~$\Q(t)$ et l'égalité ci-dessus vaut donc dans $\Q((t))$, mais nous n'aurons pas besoin de ce résultat.)

En prenant dans~(\ref{eqlefschetz1}) les dérivées logarithmiques puis en évaluant en~$t=0$, on obtient la formule (équivalente)
\begin{equation}
\label{eqlefschetz2}
\Card X(k) = \sum_{i \geq 0} (-1)^i \;\mathrm{Tr}(F \;|\; H^i(\Xbarre,\Ql))\text{,}
\end{equation}
appelée \emph{formule des traces de Lefschetz}.
Comme $H^0(\Xbarre,\Ql)=\Ql$ (avec action triviale de~$F$), le terme $i=0$ de la somme ci-dessus est égal à~$1$.
Par conséquent:

\bigskip
\begin{theoreme}[ (Grothendieck, Artin)]%
\label{thlefschetz}
Supposons que pour tout $i>0$, les valeurs propres de l'endomorphisme de $H^i(\Xbarre,\Ql)$ induit par~$F$ soient des entiers algébriques divisibles par~$q$
(\emph{i.e.}~appartiennent à $q \Zbbarre$, où~$\Zbbarre$ désigne la fermeture intégrale de~$\Z$ dans une clôture algébrique de~$\Ql$).  Alors $\congru{\Card X(k)}{1}{q}$
et en particulier $X(k)\neq\emptyset$.
\end{theoreme}

\bigskip
Que les valeurs propres de Frobenius agissant sur $H^i(\Xbarre,\Ql)$ soient
toujours des entiers algébriques (divisibles ou non) était l'une des propriétés attendues par
Weil de toute bonne théorie cohomologique sur les corps finis.  Pour la cohomologie $\ell$\nobreakdash-adique elle fut établie par Deligne~\cite{deligneappendice}.

Pour démontrer le théorème~\ref{thesn}, nous prouverons que l'hypothèse du théorème~\ref{thlefschetz} est satisfaite si~$X$ est rationnellement connexe par chaînes.

\bigskip
\begin{remarques}(i) Que les valeurs propres de~$F$ soient divisibles par~$q$ est une condition non seulement suffisante,
mais aussi presque nécessaire,
pour que $\congru{\Card X(k)}{1}{q}$.  Plus précisément, comme les hypothèses du théorème~\ref{thesn} sont stables par extension finie des scalaires (étant même de nature
géométrique), si l'on prouve que $\congru{\Card X(k)}{1}{q}$ alors on aura même prouvé que $\congru{N_s}{1}{q^s}$ pour tout $s \geq 1$.  Or cette dernière condition est équivalente
à ce que les valeurs propres de~$F$ sur $H^i(\Xbarre,\Ql)$ soient divisibles par~$q$ pour tout~$i>0$.  Cela résulte de l'égalité
\begin{equation}
\label{egalitelefschns}
N_s = 1 + \sum_{i \geq 1} (-1)^i \;\mathrm{Tr}(F^s \;|\; H^i(\Xbarre,\Ql))
\end{equation}
(c'est-à-dire la formule~(\ref{eqlefschetz2}) appliquée à la variété $X \otimes_k \F_{q^s}$ sur $\F_{q^s}$)
et du lemme de Fatou (cf.~\cite[Lemma~8.3]{illusie}),
compte tenu de l'absence de simplifications dans la somme~(\ref{egalitelefschns}). (Par «~absence de simplifications~» on entend que les endomorphismes de $H^i(\Xbarre,\Ql)$ et de $H^j(\Xbarre,\Ql)$
induits par~$F$ n'ont aucune valeur propre en commun si $i\neq j$.  Il s'agit là d'un corollaire de l'hypothèse de Riemann pour~$X$, conjecturée par Weil et démontrée par Deligne.
Nous ne nous servirons pas de l'hypothèse de Riemann pour~$X$ dans la preuve du théorème~\ref{thesn}.)

(ii)
Il existe une autre formule de Lefschetz \emph{a priori} plus prometteuse pour établir le théorème~\ref{thesn}: la \emph{formule de Woods Hole}.
C'est un avatar de la formule de Lefschetz en cohomologie cristalline mais son énoncé est élémentaire: elle ne fait intervenir que des groupes de cohomologie cohérente.
Il s'agit de l'égalité
$$
\Card X(k) = \sum_{i \geq 0} (-1)^i \; \mathrm{Tr}(F \;|\; H^i(X,\Orond_X)) \text{,}
$$
qui est une égalité entre éléments de~$k$
(puisque $H^i(X,\Orond_X)$ est un $k$\nobreakdash-espace vectoriel) et qui ne renseigne donc sur l'entier~$\Card X(k)$ que modulo~$p$ (et non modulo~$q$).
Voir~\cite[6.13.2]{sga5groth} pour une démonstration.

Il résulte de la formule de Woods Hole que l'on a $X(k) \neq \emptyset$ dès que $H^i(X,\Orond_X)=0$ pour tout $i>0$.
(Plus généralement, l'inégalité entre polygone de Newton et polygone de Hodge (cf.~\cite[4.4(c)]{illusie}) signifie que
les valeurs propres
de~$F$ agissant sur $H^i(\Xbarre,\Ql)$ sont divisibles par des puissances de~$q$ d'autant plus élevées que davantage de groupes de cohomologie de Hodge de~$X$ s'annulent.
Cela permet d'ailleurs d'obtenir une congruence modulo~$q$ (et non uniquement modulo~$p$) à partir de l'hypothèse que les groupes $H^i(X,\Orond_X)$ s'annulent pour $i>0$.)

Il est bien vrai que $H^i(X,\Orond_X)=0$ pour tout $i>0$ si~$X$ est une variété propre, lisse et rationnellement connexe sur un corps de caractéristique~$0$
(cf.~\cite[p.~16]{chambsembour}).
Malheureusement, sans hypothèse sur la caractéristique, il existe des variétés~$X$ propres, lisses et rationnellement connexes (même unirationnelles)
telles que $H^i(X,\Orond_X)\neq 0$ (cf.~\cite[p.~17]{chambsembour}).  C'est une question ouverte de savoir s'il existe de tels exemples parmi les variétés de Fano.
\end{remarques}

\subsubsection{Valeurs propres de Frobenius et coniveau}
\label{subsecvpfrob}

Le but de ce paragraphe est de relier la condition de divisibilité sur les valeurs propres de Frobenius (l'hypothèse du théorème~\ref{thlefschetz})
à une condition de «~coniveau~» sur la cohomologie de~$X$.
Fixons un entier $i \geq 0$.
La \emph{filtration par le coniveau} sur $H^i(\Xbarre,\Ql)$ est la filtration décroissante de $H^i(\Xbarre,\Ql)$ définie par la formule
$$
N^rH^i(\Xbarre,\Ql) = \bigcup \Ker\left(H^i(\Xbarre,\Ql) \longrightarrow H^i(\Xbarre \setminus \Zbarre,\Ql)\right)
$$
où la réunion porte sur les fermés $Z \subset X$ de codimension~$\geq r$.  
Elle fut introduite par Grothendieck~\cite[III, §10]{grothbrauer} (cf.~\cite[p.~300]{grothhodge} pour un erratum;
cf.~aussi~\cite[p.~118]{illusie}).  Bien entendu les sous-espaces $N^rH^i(\Xbarre,\Ql) \subset H^i(\Xbarre,\Ql)$ sont
stables par~$F$.

\bigskip
\begin{theoreme}[ (Deligne, Esnault, Katz)]%
\label{thconiveau}
Pour tout $i \geq 0$ et tout $r \geq 0$, les valeurs propres de l'endomorphisme de $N^rH^i(\Xbarre,\Ql)$ induit par~$F$
sont des entiers algébriques divisibles par~$q^r$.
\end{theoreme}

\bigskip
\begin{demo}[ du théorème~\ref{thconiveau}]%
Le théorème~\ref{thconiveau} repose sur un dévissage et sur le théorème d'intégralité de Deligne~\cite{deligneappendice}.  Celui-ci affirme
que si~$U$ est une variété lisse et séparée sur~$k$ (non nécessairement propre), alors pour tout~$i$, les valeurs propres de l'endomorphisme
de $H^i(\Ubarre,\Ql)$ induit par~$F$ sont des entiers algébriques\footnote{En toute rigueur, cette affirmation n'apparaît pas explicitement dans~\cite{deligneappendice}
mais elle résulte tout de suite de \cite[Corollaire~5.3.3 (iii)]{deligneappendice} combiné avec la dualité de Poincaré.}.

Le formalisme de la cohomologie $\ell$\nobreakdash-adique comprend une notion de cohomologie à supports: si~$U$ est une variété sur~$k$ et si $A \subset U$ est un fermé de~$U$,
Grothendieck a défini les groupes  $H^i_\Abarre(\Ubarre,\Ql)$  de cohomologie de~$\Ubarre$ à supports dans~$\Abarre$.  (En topologie on parlerait de la cohomologie de~$\Ubarre$
relative à $\Ubarre \setminus \Abarre$.) Ce sont des $\Ql$\nobreakdash-espaces vectoriels de dimension finie.
Ils s'insèrent dans une suite exacte longue du triple
\begin{equation}
\label{selocalisation}
\xymatrix@C=1.702em{
H^i_{\Bbarre}(\Ubarre,\Ql) \ar[r] & H^i_{\Abarre}(\Ubarre,\Ql) \ar[r] & H^i_{\Abarre \setminus \Bbarre}(\Ubarre \setminus \Bbarre,\Ql) \ar[r] & H^{i+1}_{\Bbarre}(\Ubarre,\Ql) \ar[r] & \cdots
}
\end{equation}
chaque fois que l'on a des fermés $B \subset A \subset U$.  En outre, il y a une flèche d'oubli de support $H^i_\Abarre(\Ubarre,\Ql) \rightarrow H^i(\Ubarre,\Ql)$, qui  est un isomorphisme si $A=U$.
Tous ces espaces sont munis d'un endomorphisme induit par l'endomorphisme de Frobenius de~$U$; la suite~(\ref{selocalisation}) commute aux Frobenius.

La suite exacte longue de la paire (\emph{i.e.}~la suite~(\ref{selocalisation}) appliquée à $U=X$, $A=X$, $B=Z$)
montre que le sous-espace
$N^rH^i(\Xbarre,\Ql) \subset H^i(\Xbarre,\Ql)$ est engendré par les images des flèches d'oubli de support
$H^i_{\Zbarre}(\Xbarre,\Ql) \rightarrow H^i(\Xbarre,\Ql)$
pour $Z \subset X$ de codimension~$\geq r$.  Il suffit donc, pour établir le théorème, de vérifier que les valeurs propres de~$F$
agissant sur $H^i_{\Zbarre}(\Xbarre,\Ql)$ sont des entiers algébriques divisibles par~$q^r$ si $Z \subset X$ est un fermé de codimension~$\geq r$.

Si $Z \subset X$ est lisse sur~$k$ et purement de codimension~$r$, le théorème de pureté (analogue au théorème d'isomorphisme de Thom pour les variétés topologiques orientées)
affirme que la cohomologie de~$\Xbarre$ à supports dans~$\Zbarre$ en degré~$i$ s'identifie à la cohomologie de~$\Zbarre$ en degré $i-2r$, à ceci près que l'on doit
«~tordre~» le Frobenius pour obtenir cette identification.
Plus précisément, l'espace vectoriel $H^i_\Zbarre(\Xbarre,\Ql)$ muni de l'endomorphisme~$F$ est isomorphe (non canoniquement)
à l'espace vectoriel $H^{i-2r}(\Zbarre,\Ql)$ muni de l'endomorphisme $q^rF$.
Il s'ensuit que les valeurs propres de l'endomorphisme de
 $H^i_\Zbarre(\Xbarre,\Ql)$ induit par~$F$ sont des entiers algébriques divisibles par~$q^r$.

Si $Z \subset X$ n'est pas lisse, posons $Z_0=Z$ et
 notons $Z_1 \subset Z_0$ le lieu singulier de~$Z_0$, puis $Z_2 \subset Z_1$ le lieu singulier de~$Z_1$, etc.
Compte tenu des suites exactes
$$
\xymatrix{
H^i_{\Zbarre_{j+1}}(\Xbarre,\Ql) \ar[r] & H^i_{\Zbarre_j}(\Xbarre,\Ql) \ar[r] & H^i_{\Zbarre_{j} \setminus \Zbarre_{j+1}}(\Xbarre \setminus \Zbarre_{j+1},\Ql)
}
$$
et de la vacuité de~$Z_j$ pour~$j$ assez grand, il suffit pour conclure de montrer que les valeurs propres de~$F$ agissant sur
$H^i_{\Zbarre_{j} \setminus \Zbarre_{j+1}}(\Xbarre \setminus \Zbarre_{j+1},\Ql)$ sont des entiers algébriques divisibles par~$q^r$. Mais cela résulte, comme précédemment,
de la pureté et du théorème d'intégralité de Deligne (appliqué aux composantes connexes de $\Zbarre_{j} \setminus \Zbarre_{j+1}$, qui sont bien lisses; noter qu'ici il est important
de disposer du théorème d'intégralité de Deligne pour des variétés lisses qui ne sont pas propres).
\end{demo}

\bigskip
En combinant les théorèmes~\ref{thlefschetz} et~\ref{thconiveau} on obtient:

\bigskip
\begin{theoreme}
\label{thconivcong}
Supposons que pour tout $i>0$ on ait $N^1H^i(\Xbarre,\Ql)=H^i(\Xbarre,\Ql)$.  Alors $\congru{\Card X(k)}{1}{q}$ et en particulier $X(k)\neq \emptyset$.
\end{theoreme}

\subsubsection{Du groupe de Chow au coniveau}
\label{subsecchowconiveau}

Un argument connu sous le nom de \emph{décomposition de la diagonale} et utilisé par Bloch~\cite{blochdiag} pour étendre en caractéristique arbitraire un théorème
de Mumford sur le groupe de Chow des $0$\nobreakdash-cycles de certaines surfaces complexes permet d'établir:

\bigskip
\begin{lemme}
\label{lemmedecompdiag}
Soit~$K$ une clôture algébrique de~$k(X)$.
Si $CH_0(X \otimes_k K)=\Z$ alors $N^1H^i(\Xbarre,\Ql)=H^i(\Xbarre,\Ql)$ pour tout $i>0$.
\end{lemme}

\bigskip
\begin{remarque}
L'hypothèse du lemme équivaut à ce que la flèche degré induise un isomorphisme
$\CH_0(X \otimes_k K) \otimes_\Z \Q \isoto \Q$, d'après un théorème de Roitman
complété par Milne (cf.~\cite[Remarque~6.5]{chambsembour}).
\end{remarque}

\bigskip
\begin{demo}
Quitte à remplacer~$k$ par une extension finie, on peut supposer que $X(k)\neq\emptyset$.
Fixons $x \in X(k)$.  Notons $\Delta \subset X \times X$ la diagonale de~$X$
et définissons un cycle algébrique~$z$ de dimension $\dim(X)$ sur $X \times X$ par la formule
 $z = \Delta - (x \times X)$.
La restriction~$z'$ de~$z$ à $X \times \eta$, où~$\eta$ désigne le point générique de~$X$, est un $0$\nobreakdash-cycle de degré~$0$
sur la $k(X)$\nobreakdash-variété $X \times \eta = X \otimes_k k(X)$.
D'après l'hypothèse, sa classe dans $\CH_0(X \otimes_k k(X))$ s'annule dans $\CH_0(X \otimes_k K)$ et donc dans $\CH_0(X \otimes_k L)$ pour une certaine extension finie $L/k(X)$.
Compte tenu de l'existence d'une application norme $\CH_0(X \otimes_k L) \rightarrow \CH_0(X \otimes_k k(X))$ dont la composée avec la flèche
de restriction $\CH_0(X \otimes_k k(X)) \rightarrow \CH_0(X \otimes_k L)$ est égale à la multiplication par $[L:k(X)]$ sur $\CH_0(X \otimes_k k(X))$,
il en résulte qu'un multiple non nul de~$z'$ est rationnellement équivalent à~$0$ sur $X \otimes_k k(X)$.  D'où l'existence d'un entier $N \geq 1$
et d'un ouvert dense $U \subset X$ tels que~$Nz$ soit rationnellement équivalent à~$0$ sur~$X \times U$.
On a alors, dans $CH_{\dim(X)}(X \times X) \otimes_\Z \Q$, la relation
\begin{equation}
\label{decompdiag}
[\Delta] = [x \times X] + \frac{1}{N}[w]
\end{equation}
où~$w$ est un cycle sur $X \times X$ dont le support ne rencontre pas $X \times \eta$: c'est la décomposition de la diagonale.

Le groupe $CH_{\dim(X)}(X \times X) \otimes_\Z \Q$ agit sur $H^i(\Xbarre,\Ql)$ par $c(\alpha)=q_\star(\cl(c) \smile p^\star \alpha)$
pour $c \in \CH_{\dim(X)}(X \times X) \otimes_\Z \Q$ et $\alpha \in H^i(\Xbarre,\Ql)$, où $p,q\colon X\times X \rightarrow X$ sont les deux projections, où $\smile$ désigne le cup-produit
et où $\cl(c) \in H^{2\dim(X)}(\Xbarre \times \Xbarre,\Ql)$ est la classe de cohomologie du cycle~$c$.
Le cycle $\Delta$ agit sur $H^i(\Xbarre,\Ql)$ par l'identité; le cycle $x \times X$ agit sur $H^i(\Xbarre,\Ql)$ par~$0$ pour $i>0$
puisque
$$
\cl(x \times X) \smile p^\star \alpha = p^\star \cl(x) \smile p^\star \alpha = p^\star (\cl(x) \smile \alpha) \in p^\star H^{2\dim(X)+i}(\Xbarre,\Ql)
$$
et $H^{2\dim(X)+i}(\Xbarre,\Ql)=0$ pour $i>0$; enfin, le cycle~$w$ agit sur $H^i(\Xbarre,\Ql)$ par un endomorphisme dont l'image
est incluse dans $N^1H^i(\Xbarre,\Ql)$.  La décomposition~(\ref{decompdiag}) montre donc que l'identité de $H^i(\Xbarre,\Ql)$ est à valeurs dans $N^1H^i(\Xbarre,\Ql)$ pour tout~$i>0$.
\end{demo}

\bigskip
En emboîtant le théorème~\ref{thconivcong} avec le lemme~\ref{lemmedecompdiag} on obtient:

\bigskip
\begin{theoreme}
\label{thfinalch}
Soit~$K$ une clôture algébrique de~$k(X)$. Si $\CH_0(X \otimes_k K)=\Z$ alors $\congru{\Card X(k)}{1}{q}$ et en particulier $X(k)\neq\emptyset$.
\end{theoreme}

\bigskip
Une variété rationnellement connexe par chaînes vérifie de façon évidente l'hypothèse du théorème~\ref{thfinalch} puisque tous les points rationnels de $X \otimes_k K$ sont
même $R$\nobreakdash-équivalents.  Le théorème~\ref{thesn} est donc prouvé.

\subsection{Au-delà du théorème~\ref{thesn}}
\label{subsecapres}

\medskip
Soient~$k$ un corps fini de cardinal~$q$ et~$X$ une variété projective sur~$k$.

Le théorème~\ref{thfinalch} affirme que si~$X$ est lisse et rationnellement connexe par chaînes, alors $\congru{\Card X(k)}{1}{q}$.

Le théorème de Chevalley--Warning amélioré par Ax~\cite[p.~260]{axzeroes} affirme que
si~$X$ est l'intersection d'hypersurfaces de~$\P^n_k$ de degrés $\uplet{d_1}{d_r}$ avec $\sum_{i=1}^r d_i \leq n$, alors $\congru{\Card X(k)}{1}{q}$.

Dans l'énoncé (et dans la preuve) du théorème de Chevalley--Warning (et de l'amélioration due à Ax), la lissité de~$X$ ne joue aucun rôle.
Le théorème~\ref{thfinalch} n'est donc pas optimal et il est naturel de chercher à généraliser sa démonstration de façon
à retrouver les résultats de Chevalley--Warning et d'Ax.

Malheureusement, si~$X$ n'est pas lisse, la preuve présentée au~§\ref{thesn} s'effondre:
d'une part la décomposition de la diagonale~(\ref{decompdiag}) ne fournit plus l'information cohomologique du lemme~\ref{lemmedecompdiag}, faute de disposer
d'une application classe de cycle $\CH_{\dim(X)}(X \times X) \rightarrow H^{2\dim(X)}(\Xbarre \times \Xbarre,\Ql)$;
d'autre part, du théorème d'intégralité de Deligne (qui reste vrai sans hypothèse de lissité, cf.~\cite[1.2]{hodgeweightanalogies}) il n'est plus possible de déduire
la divisibilité des valeurs propres de Frobenius agissant sur $N^1H^i(\Xbarre,\Ql)$,
faute de pouvoir appliquer le théorème de pureté
(et de fait, le théorème~\ref{thconiveau} est faux pour les variétés singulières; cf.~\cite[0.5]{deligneesnaultapp}).

En outre, pour~$X$ singulière, il peut y avoir des «~simplifications~» dans la formule des traces de Lefschetz (cf.~remarque~(i) après le théorème~\ref{thlefschetz}),
de sorte que la congruence $\congru{\Card X(k)}{1}{q}$ pourrait être satisfaite sans que les valeurs propres de Frobenius agissant sur $H^i(\Xbarre,\Ql)$ ne soient
divisibles par~$q$ pour $i>0$.  Cependant, dans la situation considérée par Chevalley, Warning et Ax, la philosophie des motifs fournit des indications\footnote{Il existe une étroite analogie, découverte par Grothendieck et Deligne,
entre la divisibilité par~$q^\kappa$
des valeurs propres de Frobenius agissant sur $H^i(\Xbarre,\Ql)$ (lorsque~$k$ est fini) et la trivialité des~$\kappa$ premiers gradués de la filtration de Hodge
sur $H^i(X(\C),\C)$ (lorsque $k=\C$). (Au moins si~$X$ est lisse,
ces deux propriétés devraient rendre compte de la «~divisibilité~» du motif effectif
$H^i(X)$ par le motif de Tate $\Q(-\kappa)$.)  Or il est bien vrai que si~$X$ est une intersection d'hypersurfaces de~$\P^n_\C$ de degrés $\uplet{d_1}{d_r}$ avec
$\sum_{i=1}^r d_i \leq n$, le premier gradué de la filtration de Hodge sur $H^i(X(\C),\C)$ est nul si $i>0$ (cf.~\cite{esnaultnorisrinivas}).
Pour plus de détails
on pourra consulter~\cite{delignedimca} ou~\cite[§1]{hodgeweightanalogies}.}
qui permettent de s'attendre
à ce que les valeurs propres de Frobenius
agissant sur $H^i(\Xbarre,\Ql)$ soient quand même divisibles par~$q$ pour tout $i>0$; cela a d'ailleurs été vérifié \emph{a posteriori}, c'est-à-dire en se servant des résultats
de Chevalley--Warning et d'Ax, par Esnault et Katz~\cite{esnaultkatz}.  Bloch, Esnault et Levine~\cite{blochesnaultlevine} ont ainsi réussi à adapter
la preuve du théorème~\ref{thfinalch}, en établissant une variante de la décomposition de la diagonale, pour retrouver \emph{via} le théorème~\ref{thlefschetz} la congruence
$\congru{\Card X(k)}{1}{q}$ si~$X$ est une hypersurface de~$\P^n_k$ de degré~$d$ avec $d \leq n$.
Mais leur résultat ne couvre pas le cas plus général des intersections d'hypersurfaces de degrés $\uplet{d_1}{d_r}$ avec $\sum_{i=1}^r d_i \leq n$
et ne permet pas non plus de retrouver les énoncés plus précis d'Ax et Katz~\cite{katzsurax} donnant des congruences modulo des puissances supérieures de~$q$.

Plutôt que d'essayer d'adapter la méthode de preuve du théorème~\ref{thesn}, on peut aussi chercher un énoncé qui généralise à la fois le théorème~\ref{thesn}
et le théorème de Chevalley--Warning.
Si $\sum_{i=1}^r d_i \leq n$,
l'intersection de~$r$ hypersurfaces génériques de~$\P^n_k$ de degrés $\uplet{d_1}{d_r}$
est une variété lisse et de Fano, donc rationnellement connexe par chaînes.
Par spécialisation, il s'ensuit\footnote{La connexité rationnelle par chaînes
est préservée par spécialisation équidimensionnelle (cf.~\cite[IV.3.5.2]{kollrational}).
Ainsi, si $X \subset \P^n_k$ désigne l'intersection de~$r$
hypersurfaces de degrés $\uplet{d_1}{d_r}$ avec $\sum_{i=1}^r d_i \leq n$, la connexité rationnelle par chaînes de~$X$ s'ensuit de façon immédiate si $\dim(X)=n-r$.
Sans hypothèse sur la dimension de~$X$, on remarque qu'étant donnés deux points~$P$
et~$Q$ de~$X$, l'intersection de~$r$ hypersurfaces de~$\P^n_k$
de degrés $\uplet{d_1}{d_r}$,
contenant~$P$ et~$Q$ et génériques pour ces propriétés, est de dimension~$n-r$. Cette variété est rationnellement connexe par chaînes et
se spécialise en une sous-variété rationnellement connexe par chaînes de~$X$ contenant~$P$ et~$Q$, ce qui démontre la connexité rationnelle par chaînes de~$X$ puisque~$P$ et~$Q$ sont arbitraires.
Je~remercie Jason Starr pour cet argument.} que toute intersection dans~$\P^n_k$
d'hypersurfaces de degrés $\uplet{d_1}{d_r}$ avec $\sum_{i=1}^r d_i \leq n$ est rationnellement connexe par chaînes (bien que non nécessairement lisse).
Ainsi, le premier espoir venant à l'esprit
en direction d'une généralisation commune serait que sur un corps
fini, toute variété projective, géométriquement irréductible et
rationnellement connexe par chaînes (mais non nécessairement lisse) admette un point rationnel.
Cet espoir est malheureusement mis en défaut par un exemple de variété normale, projective, géométriquement irréductible, rationnellement connexe par chaînes et dépourvue de points rationnels,
sur un corps fini, construit par Kollár
(cf.~\cite[Remark~3.4]{fakhruddinrajan}).  Néanmoins, Fakhruddin et Rajan~\cite{fakhruddinrajan} se sont
rendu compte que l'existence d'une déformation lisse et rationnellement connexe par chaînes de~$X$ suffit à assurer que $X(k)\neq\emptyset$.  Plus précisément, ils ont établi
le théorème suivant, qui généralise bien à la fois le théorème~\ref{thesn} et le théorème de Chevalley--Warning revu par Ax:

\bigskip
\newcommand{\citefr}{\cite[Corollary~1.2]{fakhruddinrajan}}
\begin{theoreme}[ (Fakhruddin--Rajan~\citefr)]%
\label{thfakhruddinrajan}
Soit~$k$ un corps fini de cardinal~$q$.
Soit $f \colon Y \rightarrow S$ un morphisme propre et surjectif entre variétés irréductibles et lisses sur~$k$.
Notons~$Z$ la fibre générique de~$f$ et~$K$ une clôture algébrique de~$k(Y)$.
Si $\CH_0(Z \otimes_{k(S)} K)=\Z$ alors $\congru{\Card f^{-1}(s)(k)}{1}{q}$ pour tout $s \in S(k)$.
\end{theoreme}

\bigskip
Le théorème~\ref{thfakhruddinrajan} est à comparer avec le théorème de
Hogadi et Xu~\cite{hogadixu} que nous signalions à la fin du~§\ref{subsubsecautrescorps},
 selon lequel sur un corps de caractéristique nulle,
toute variété obtenue en faisant dégénérer une variété rationnellement connexe
contient une sous-variété géométriquement irréductible (et même une sous-variété
rationnellement connexe).
Fakhruddin et Rajan montrent que cet énoncé vaut aussi sur un corps fini,
un point rationnel étant en effet \emph{a fortiori} une sous-variété
géométriquement irréductible (et même rationnellement connexe).

Esnault a étendu le théorème~\ref{thfakhruddinrajan} à une situation d'inégale caractéristique:

\bigskip
\begin{theoreme}[ (Esnault~\cite{esnaultann})]%
\label{thesnann}
Soit~$R$ un anneau de valuation discrète complet à corps résiduel fini.  Notons~$K$ le corps des fractions de~$R$ et~$k$ son corps résiduel.
Soit~$Y$ une variété projective, lisse et géométriquement connexe sur~$K$ admettant un modèle régulier~$\Yrond$ projectif sur~$R$.
Notons~$\Ybarre$ la variété déduite de~$Y$ par extension des scalaires de~$K$ à une clôture algébrique de~$K$ et fixons un nombre premier~$\ell$ inversible dans~$K$.
Si $N^1H^i(\Ybarre,\Ql)=H^i(\Ybarre,\Ql)$ pour tout $i>0$, alors $\congru{\Card \Yrond(k)}{1}{q}$, où~$q$ désigne le cardinal de~$k$.
\end{theoreme}

\bigskip
Lorsque la base est à la fois de dimension~$1$ et d'égale caractéristique,
le théorème~\ref{thesnann} améliore le théorème~\ref{thfakhruddinrajan}
puisque l'hypothèse sur le groupe de Chow y est remplacée par la condition que la cohomologie $\ell$\nobreakdash-adique
de la fibre générique géométrique est de coniveau~$\geq 1$.
Ces deux hypothèses devraient être équivalentes (cela résulterait d'une généralisation de la conjecture de Bloch, cf.~\cite[§3]{jannsenmotives}) mais l'hypothèse sur le coniveau
est plus facile à vérifier que l'hypothèse sur le groupe de Chow. Par exemple, si~$K$ est de caractéristique~$0$, elle est satisfaite pour toute surface~$Y$ telle que $H^1(Y,\Orond_Y)=0$
et $H^2(Y,\Orond_Y)=0$ (cf.~\cite[§1]{esnaultann}).

Tout récemment, Berthelot, Esnault et Rülling~\cite{berthelotesnaultruelling} ont réussi
à démontrer que le théorème~\ref{thesnann} reste valable si l'on remplace l'hypothèse sur le coniveau par la condition
que $H^i(Y,\Orond_Y)=0$ pour tout $i>0$.
D'après la théorie de Hodge mixte de Deligne,
cette condition est impliquée par l'hypothèse sur le coniveau.
Il résulterait de la conjecture de Hodge généralisée par Grothendieck~\cite{grothhodge} que les deux sont équivalentes.
Des exemples de variétés de type général~$Y$ sur des extensions finies non ramifiées de~$\Qp$ vérifiant $H^i(Y,\Orond_Y)=0$ pour $i>0$ mais pour lesquelles on ne sait pas vérifier que $N^1H^i(\Ybarre,\Ql)=H^i(\Ybarre,\Ql)$ sont donnés
dans~\cite[\textsection9]{berthelotesnaultruelling}.

Signalons enfin que Blickle et Esnault~\cite{blickleesnault} ont étendu le théorème~\ref{thfinalch} à une classe de variétés non nécessairement lisses (mais dont les singularités
sont «~Witt-rationnelles~», cf.~\cite[Definition~2.3]{blickleesnault}).

Dans une toute autre direction, Kahn~\cite{kahnmotives} a donné une nouvelle preuve des théorèmes~\ref{thesn} et~\ref{thfinalch}.
Plus précisément, la théorie des motifs birationnels, due à Kahn et Sujatha, entraîne que
si~$X$ est une variété projective et lisse sur un corps~$k$
telle que
$\CH_0(X \otimes_k K)=\Z$, où~$K$ est une clôture algébrique de~$k(X)$,
alors
le motif de Chow de~$X$ (à coefficients rationnels) admet une décomposition de la forme
$$
\mathfrak{h}(X) \simeq \mathbf{1} \oplus (\mathbf{L}\otimes M)
$$
où $\mathbf{L}$ désigne le motif de Lefschetz et~$M$ un motif de Chow effectif
(cf.~\cite[\textsection7.3]{kahnmotives}).  Lorsque~$k$ est un corps fini de cardinal~$q$, cette décomposition se traduit par l'égalité
$\Card X(k) = 1 + qm$, où $m \in \Q$ est la trace de l'endomorphisme de Frobenius sur le motif~$M$.
Or, d'après Kahn~\cite[Theorem~8.1]{kahnmotives}, on a $m \in \Z$; ainsi retrouve-t-on bien la congruence $\congru{\Card X(k)}{1}{q}$.
Des considérations analogues permettent d'établir:

\bigskip
\newcommand{\refkahn}{\cite[Corollary~9.6]{kahnmotives}}
\begin{theoreme}[ (Kahn~\refkahn)]%
\label{thkahn}
Soient~$k$ un corps fini de cardinal~$q$ et $X$, $Y$ des variétés projectives et lisses sur~$k$.
Supposons qu'il existe des variétés projectives lisses rationnellement connexes par chaînes~$X'$ et~$Y'$ sur~$k$ et des applications rationnelles dominantes
$X \times X' \dashrightarrow Y$ et $Y \times Y' \dashrightarrow X$.  Alors $\congru{\Card X(k)}{\Card Y(k)}{q}$.
\end{theoreme}

\bigskip
En prenant $X'=Y=\Spec(k)$ et $Y'=X$, on retrouve le théorème~\ref{thesn}.  En
prenant $X'=Y'=\Spec(k)$, on retrouve, au moins dans le cas projectif, un théorème d'Ekedahl revu par Chambert-Loir (cf.~\cite{chambsembour}).
Nous renvoyons à~\cite{kahnmotives} pour d'autres énoncés de ce type.

\addcontentsline{toc}{section}{Bibliographie}
\bibliographystyle{rcarithbibstyle}
\bibliography{rcarith}
\end{document}